\documentclass[11pt]{article}

\usepackage[margin=1in]{geometry}
\usepackage{amsmath,amssymb,amsthm}
\usepackage[authoryear,round]{natbib}
\usepackage[colorlinks=true,linkcolor=blue,citecolor=blue,urlcolor=blue]{hyperref}
\usepackage{cleveref}

\usepackage{algorithm}
\usepackage{algorithmic}

\newcommand{\sD}{{\mathsf{D}}}
\newcommand{\sW}{{\mathsf{W}}}
\newcommand{\sS}{{\mathsf{S}}}
\newcommand{\sC}{{\mathsf{C}}}
\newcommand{\neff}{{\mathsf{n_{\rm eff}}}}

\newcommand{\rmJ}{{\mathrm{J}}}
\theoremstyle{plain}
\newtheorem{assumption}{Assumption}
\newtheorem{definition}{Definition}
\newtheorem{lemma}{Lemma}
\newtheorem{theorem}{Theorem}
\crefname{assumption}{assumption}{assumptions}
\Crefname{assumption}{Assumption}{Assumptions}
\Crefname{definition}{Definition}{Definitions}

\newtheorem{remark}{Remark}

\title{Optimal Neural Estimation for High-Dimensional Nonparametric Regression under Multisource Covariate Shift}

\author{%
  \begin{minipage}{0.96\textwidth}
  \centering
  Chenghao Xu\textsuperscript{1},
  Haotian Xu\textsuperscript{2},
  Shivam Kumar\textsuperscript{3},
  and
  Daren Wang\textsuperscript{1}\\[0.85em]
  \textsuperscript{1}Department of Mathematics, University of California, San Diego\\
  \textsuperscript{2}Department of Mathematics and Statistics, Auburn University\\
  \textsuperscript{3}Booth School of Business, University of Chicago
  \end{minipage}%
}
\date{September 26, 2026}

\begin{document}
\maketitle

\begin{abstract}
We study high-dimensional nonparametric regression under multisource covariate shift, where $\mathrm J$ sources of training data share a common regression function but differ in their covariate distributions, sample sizes, and noise variances. Our goal is to combine these sources to estimate the regression function accurately under a target covariate distribution. We propose an oracle weighted neural estimator and establish its minimax optimality up to logarithmic factors through matching upper and lower bounds. The resulting rate is $\mathsf n_{\rm eff}^{-2p/(2p+K)}$, where $p$ and $K$ denote the effective smoothness and dimension of the regression function, and the effective sample size is
\begin{align}
\mathsf n_{\rm eff}
=\sum_{j=1}^{\mathrm J}
\frac{n_j}{(1+\sigma_{\rho,j})\sigma_{\epsilon,j}^{2}}.
\nonumber
\end{align}
Here, $n_j$ is the sample size of source $j$, $\sigma_{\epsilon,j}^{2}$ is its noise variance, and $\sigma_{\rho,j}$ controls the deviation between the covariate distribution of source $j$ and the target covariate distribution. We also develop an adaptive procedure that requires no knowledge of source noise variances or covariate-shift parameters. Under suitable conditions, this procedure attains the oracle minimax rate up to logarithmic factors. 
\end{abstract}

\section{Introduction}\label{sec:introduction}

Artificial intelligence plays an increasingly prominent role in prediction, and vast amounts of data are available from diverse sources. In many applications, the relationship between the covariates and the mean response remains the same across sources, while the covariate distributions can vary substantially. For example, surveys may sample demographic groups in different proportions, and online advertising systems may collect data under different ad-selection policies. We study this form of covariate shift across multiple data sources and ask how to combine their observations to improve prediction for a target population of interest.

To formalize this setting, suppose we observe data from $\mathrm{J}$ sources. For source $j$, the observations follow
\begin{align}
Y_{ij}&=f^*(X_{ij})+\epsilon_{ij}, \qquad X_{ij}\sim\mathcal{P}_j, \qquad i=1,\ldots,n_j,\quad j=1,\ldots,\mathrm{J},
\label{eq:multisource-model}
\end{align}
where $\mathbb{E}[\epsilon_{ij}\mid X_{ij}]=0$ and $\sigma_{\epsilon,j}^2=\operatorname{Var}(\epsilon_{ij})$ denotes the noise variance in source $j$. The sources share a common regression function $f^*$, while both the covariate distribution $\mathcal{P}_j$ and the noise variance $\sigma_{\epsilon,j}^2$ may differ across sources. Our goal is to estimate $f^*$ using the source observations, with estimation error measured under a target covariate distribution $\mathcal{Q}$, which may differ from each $\mathcal{P}_j$.

This framework captures several common applications. In online advertising, click data are collected under different ad-selection policies, each inducing a source distribution over user characteristics, displayed advertisements, and display contexts \citep{agarwal2017}.  When the click probability given these covariates remains unchanged across policies, historical observations can be combined to predict responses under the new policy. The same framework applies to surveys that sample demographic groups in different proportions and experiments conducted under different input designs, provided that changes in sampling or design preserve the conditional mean response \citep{shimodaira2000,sugiyama2007}.

\subsection{Literature review}
Early work on covariate shift developed importance-weighting methods that adjust the contribution of source observations to better represent the target population \citep{shimodaira2000,sugiyama2007}. Subsequent studies developed data-driven weighting procedures \citep{huang2006kmm,bickel2009} and examined how variable or extreme weights affect estimation accuracy, motivating theoretical error bounds and truncation methods \citep{cortes2010importance,ionides2008}.

Kernel and local smoothing methods have provided sharp statistical guarantees for nonparametric regression under covariate shift. \citet{ma2023} develop kernel ridge regression estimators that attain nearly minimax rates under bounded or potentially unbounded importance weights. \citet{wang2026} uses pseudo-labels on unlabeled target observations to select the regularization parameter in kernel ridge regression, achieving nearly minimax rates while adapting to the target distribution and covariate shift. For Nadaraya--Watson kernel regression, \citet{pathak2022} characterize how source--target overlap determines the minimax estimation rate. \citet{deb_phase} study kernel-weighted local polynomial regression when the target covariates concentrate near a low-dimensional manifold. Despite these advances, conventional kernel and local smoothing methods remain subject to the curse of dimensionality over general smooth function classes, with convergence rates deteriorating as the number of covariates increases. This motivates methods that exploit additional structure in high-dimensional regression functions.

Neural network methods can adapt to low-dimensional structure in high-dimensional nonparametric regression problems. In the standard setting without covariate shift, \citet{schmidthieber2020} and \citet{kohlerlanger2021} establish nearly minimax rates for regression functions formed by hierarchical compositions of simpler functions. These results provide a foundation for studying whether the advantages of neural estimation can be retained under covariate shift.

Neural methods have also been developed for regression under covariate shift. \citet{feng2024} study deep quantile regression with estimated importance weights, while \citet{xu2025} identify conditions under which unweighted neural regression retains accurate target prediction. \citet{kengne2026} develop penalized neural estimators that adapt to unknown smoothness and accommodate dependent observations. \citet{tilt2026} incorporate unlabeled target observations into the training objective and establish guarantees for sparse ReLU networks. These studies focus on a single source. Optimal prediction rates and adaptive procedures for high-dimensional nonparametric regression under multisource covariate shift remain largely unexplored in the literature.

\subsection{Contributions}
Our main contributions are as follows.
\begingroup
\setlength{\leftmargini}{0.5em}
\begin{itemize}
\item \textbf{A growing number of heterogeneous sources.}
Our framework allows the number of sources to increase with the total sample size, with sources differing in their sample sizes, covariate distributions, and noise variances. Our proposed approach combines labeled source observations to improve target prediction.

\item \textbf{Minimax rates in the oracle setting.}
Define the effective sample size as
\begin{align}
\mathsf{n}_{\rm eff}
=\sum_{j=1}^{\mathrm J}
\frac{n_j}{(1+\sigma_{\rho,j})\sigma_{\epsilon,j}^{2}},
\label{eq:effective-sample-size}
\end{align}
where, for source $j$, $n_j$ is the sample size, $\sigma_{\epsilon,j}^{2}$ is the noise variance, and $\sigma_{\rho,j}$ measures the scale of covariate shift. See \Cref{ass:holder-shift} for details of this parameter. For regression functions in the hierarchical H\"older class defined in \Cref{def:hierarchical}, our oracle estimator achieves a prediction error of order $\mathsf{n}_{\rm eff}^{-2p/(2p+K)}$ up to logarithmic factors, where $p$ and $K$ denote the effective smoothness and effective dimension, respectively, of the high-dimensional regression function. We also establish a new matching lower bound in \Cref{thm:holder-covariate-shift-lower-bound}, showing that the proposed oracle procedure is minimax optimal up to log factors.

\item \textbf{An adaptive procedure.}
We develop a new data-driven neural estimator that combines estimated importance weights with source weighting and validation-based tuning. The procedure adapts to unknown source noise variances and covariate-shift parameters. Under suitable assumptions, it recovers the oracle minimax rate up to logarithmic factors.

\end{itemize}
\endgroup

\section{Oracle estimation}\label{sec:hierarchical-covariate-shift}
In this section, we consider an oracle procedure that uses the true density ratios and source noise variances to combine observations across sources. The oracle analysis establishes the minimax rate up to logarithmic factors and guides the development of the adaptive procedure in \Cref{sec:multisource-adaptive}. For each source $j=1,\ldots,\mathrm J$, we observe a sample
\begin{align}
\mathcal D_j&=\{(X_{ij},Y_{ij}):i=1,\ldots,n_j\}
\nonumber
\end{align}
from the regression model \eqref{eq:multisource-model}. The observations are identically distributed within each source, and the covariates are mutually independent across all samples. We take $[0,1]^d$ as the common covariate domain for the source and target distributions.
Our goal is to develop a neural estimator $\widehat f$ that controls the target prediction error
\begin{align}
\mathbb E\|\widehat f-f^*\|_{\mathcal Q}^{2}
=
\mathbb E\int_{[0,1]^d}
\left( \widehat f(x)-f^*(x) \right)^{2}\,d\mathcal Q(x).
\nonumber
\end{align}

As in \Cref{sec:introduction}, let $\mathcal P_j$ denote the covariate distribution of source $j$ and $\mathcal Q$ the target covariate distribution. We assume $\mathcal Q\ll\mathcal P_j$ for each source and define the density ratios by
\begin{align}
\rho_j&=\frac{d\mathcal Q}{d\mathcal P_j},
\qquad j=1,\ldots,\mathrm J.
\nonumber
\end{align}

We impose the following conditions on the response noise and density ratios.
\begin{assumption}[\textbf{Noise}]\label{ass:noise}
The errors $\{\epsilon_{ij}:1\le i\le n_j,\ 1\le j\le\mathrm J\}$ are mutually independent and independent of all covariates $\{X_{ij}:1\le i\le n_j,\ 1\le j\le\mathrm J\}$. Write $\sigma_{\epsilon,j}^2=\operatorname{Var}(\epsilon_{ij})$. There is an absolute constant $C_\epsilon\ge1$ such that, for every $j=1,\ldots,\mathrm J$ and $i=1,\ldots,n_j$,
\begin{align}
\mathbb E[\epsilon_{ij}]=0,
\qquad
\mathbb E\exp(t\epsilon_{ij})
\le\exp\left(\frac{C_\epsilon\sigma_{\epsilon,j}^2t^2}{2}\right)
\qquad\text{for every }t\in\mathbb R.
\nonumber
\end{align}
\end{assumption}

\begin{assumption}[\textbf{Density-ratio}]
\label{ass:holder-shift}
For each $j$, the density ratio $\rho_j$ satisfies
\begin{align}
\mathbb E_{X\sim\mathcal P_j}
\left[\exp\{t(\rho_j(X)-1)\}\right]
\le\exp\left(\frac{\sigma_{\rho,j}^2t^2}{2}\right)
\qquad\text{for every }t\in\mathbb R,\quad j=1,\ldots,\mathrm J,
\nonumber
\end{align}
for some parameter $0\le\sigma_{\rho,j}<n_j$.
\end{assumption}

\Cref{ass:noise} is standard in the literature, and we allow the noise variance to vary across sources. \Cref{ass:holder-shift} controls density-ratio tails through a sub-Gaussian condition, without imposing a uniform upper bound. Bounded ratios, often considered in single-source analyses \citep{cortes2010importance,ma2023,feng2024}, satisfy this tail condition. Since $\mathbb E_{X\sim\mathcal P_j}[\rho_j(X)]=1$, the condition
controls the centered version of $\rho_j$. 
If $\sigma_{\rho,j}=0$, then $\rho_j(X)=1$
$\mathcal P_j$-almost surely, and hence $\mathcal P_j=\mathcal Q$.

We next impose structure on the regression function. To accommodate high-dimensional settings, we consider a hierarchical H\"older class, in which the regression function is a composition of smooth functions with a small number of inputs. We first define H\"older smoothness.

\begin{definition}
\label{def:holder}
Let $p\ge1$, $R>0$, $q=\lceil p\rceil-1$, and $s=p-q\in(0,1]$.
Write $\mathbb N_0=\{0,1,\ldots\}$. For a multi-index $\alpha=(\alpha_1,\ldots,\alpha_K)\in\mathbb N_0^K$, let $|\alpha|=\sum_{j=1}^K\alpha_j$ and let $\partial^\alpha u$ denote the corresponding partial derivative, with $\partial^0u=u$.
A function $u$ defined on $[0,1]^K$ or $[-R,R]^K$ is called $(p,R)$-smooth
if it has continuous partial derivatives through order $q$ and satisfies
\begin{align}
\sum_{|\alpha|\le q}\|\partial^\alpha u\|_\infty
+\max_{|\alpha|=q}\sup_{z\ne w}
\frac{|\partial^\alpha u(z)-\partial^\alpha u(w)|}
{|z-w|_2^s}\le R,
\nonumber
\end{align}
where the supremums are taken over the domain of $u$, and $|\cdot|_2$ denotes the Euclidean norm.
\end{definition}

\begin{definition}
\label{def:hierarchical}
Let $L\ge1$ and $1\le K\le d$ be integers, and let $p\ge1$ and $R>0$.
Define the hierarchical composition class recursively.
For $\ell=1$, let
\begin{align}\nonumber
\mathcal H(1,p,K,R)
=\Bigl\{h:[0,1]^d\to[-R,R]:\;
h(x)=u(x_{\pi(1)},\ldots,x_{\pi(K)}),
u:[0,1]^K\to[-R,R]&
\\\nonumber
\text{ is $(p,R)$-smooth},
\, \pi:\{1,\ldots,K\}\to\{1,\ldots,d\}&\Bigr\}.
\end{align}
For $2\le\ell\le L$, let
\begin{align}\nonumber
\mathcal H(\ell,p,K,R)
=\Bigl\{h:[0,1]^d\to[-R,R]:\;
h(x)=u(h_1(x),\ldots,h_K(x)),
u:[-R,R]^K\to[-R,R] &
\\\nonumber
\text{ is $(p,R)$-smooth},
\; h_j\in\mathcal H(\ell-1,p,K,R),
\quad j=1,\ldots,K & \, \Bigr\}.
\end{align}
\end{definition}
The hierarchical class in \Cref{def:hierarchical} has effective smoothness $p$ and effective dimension $K$. Deep neural estimators can adapt to this structure and avoid the curse of dimensionality \citep{schmidthieber2020,kohlerlanger2021}. Classical methods such as wavelet estimators \citep{donoho1995wavelet} may not adapt to this structure. For example, \citet[Section~5]{schmidthieber2020} gives models in which the rate of wavelet series estimators depends on the ambient dimension, while the rate of deep neural estimators depends on the effective dimension.

\begin{assumption}[Hierarchical composition]\label{ass:hierarchical-composition}
For $M>0$, the regression function $f^*$ in \eqref{eq:multisource-model} satisfies
\begin{align}
f^*\in\mathcal H(L,p,K,R),
\qquad\sup_{x\in[0,1]^d}|f^*(x)|\le M.
\nonumber
\end{align}
\end{assumption}
Under \Cref{ass:hierarchical-composition}, the regression function has effective smoothness $p$ and effective dimension $K$. These parameters determine the prediction-rate exponent.

\subsection{Oracle estimator}
We construct the oracle estimator by weighted least squares over a class of deep ReLU networks. We first define the network class.

\begin{definition}[ReLU network class]\label{def: nn class}
Let $\mathsf D\ge2$, $\mathsf W\ge1$, and $\mathsf S\ge1$ be integers, and let $\mathsf C_A\ge1$ and $\mathsf C_M>0$. A network of depth $\mathsf D$ has $\mathsf D-1$ hidden layers, each of width $\mathsf W$, and a scalar output. The network function is given by
\begin{align}
f_\theta(x)
&=
A_{\mathsf D}\,
r\Bigl(
A_{\mathsf D-1}\,
r\bigl(\cdots r(A_1x+b_1)\cdots\bigr)
+b_{\mathsf D-1}
\Bigr)
+b_{\mathsf D},
\nonumber
\end{align}
where $r(u)=\max\{u,0\}$ acts coordinatewise and $\theta$ collects
the weight matrices and biases of compatible dimensions.
We define
\begin{align}
\mathcal F(\mathsf D,\mathsf W,\mathsf S,\mathsf C_A,\mathsf C_M)
=\Bigl\{f_\theta:\;&
\max_{1\le a\le\mathsf D}
\{|A_a|_{\max},|b_a|_\infty\}\le\mathsf C_A,\nonumber\\
&\sum_{a=1}^{\mathsf D}
\bigl(|A_a|_0+|b_a|_0\bigr)\le\mathsf S,\qquad
\sup_{x\in[0,1]^d}|f_\theta(x)|\le\mathsf C_M
\Bigr\},
\nonumber
\end{align}
where $|A|_{\max}$ is the largest absolute matrix entry and $|\cdot|_0$ counts nonzero entries.
\end{definition}


We define the oracle estimator as the minimizer of the weighted least-squares criterion over the network class in \Cref{def: nn class},
\begin{align}
\widehat f
&\in
\operatorname*{arg\,min}_{
f\in\mathcal F(\mathsf D,\mathsf W,\mathsf S,\mathsf C_A,\mathsf C_M)}
\sum_{j=1}^{\mathrm J}\sum_{i=1}^{n_j}
w_{j,\tau}(X_{ij})
\{Y_{ij}-f(X_{ij})\}^{2},
\label{eq:oracle-estimator}
\end{align}
where, for positive thresholds $\tau_j$, we choose the weights
\begin{align}
w_{j,\tau}(x)
&=
\frac{1}{\sigma_{\epsilon,j}^{2}}
\min\left\{\frac{\rho_j(x)}{\tau_j},1\right\},
\qquad j=1,\ldots,\mathrm J.
\label{eq:oracle-generic-weights}
\end{align}
The inverse noise variances $\sigma_{\epsilon,j}^{-2}$ give greater weight to observations from sources with less noise.
The density ratios $\rho_j$ adjust for differences between the source and target covariate distributions. The weights in \eqref{eq:oracle-generic-weights} can be written as
\begin{align}
w_{j,\tau}(x)
=\frac{\min\{\rho_j(x),\tau_j\}}
{\sigma_{\epsilon,j}^{2}\tau_j}.
\nonumber
\end{align}
Thus, each density ratio exceeding $\tau_j$ is replaced by $\tau_j$, limiting the influence of observations with large density ratios. The amount removed from the density ratio is $\rho_j(x)-\min\{\rho_j(x),\tau_j\}$. Under \Cref{ass:holder-shift}, the thresholds defined in \eqref{eq:oracle-thresholds} satisfy
\begin{align}
\mathbb E_{X\sim\mathcal P_j}
\left[\rho_j(X)-\min\{\rho_j(X),\tau_j\}\right]
=o(1).
\label{eq:oracle-density-ratio-reduction}
\end{align}
Thus, the expected reduction in each density ratio tends to zero as the total sample size increases.

\subsection{Statistical guarantees}\label{sec:oracle-conditions}

We next provide statistical guarantees for the oracle estimator in \eqref{eq:oracle-estimator}.  

\begin{theorem}[Oracle prediction bound]
\label{thm:estimated-covariate-shift-rate}
Suppose \Cref{ass:noise,ass:holder-shift,ass:hierarchical-composition} hold and $\sigma_{\epsilon,j}^2\ge c_\epsilon$ for every $j=1,\ldots,\mathrm J$, where $c_\epsilon>0$ is fixed. Assume that the effective sample size defined in \eqref{eq:effective-sample-size} satisfies

\begin{align}\label{eq:oracle-effective-sample-size}
    \mathsf{n}_{\rm eff} \ge
    \sqrt{12\log\left(e\sum_{k=1}^{\mathrm J}n_k\right)}.
\end{align}
Fix $\mathsf C_M\ge M$. There exist constants $\mathsf C_A\ge1$ and $C_1,C_2,C_3,C_4>0$, depending only on $d$, $L$, $p$, $K$, $R$, $M$, $c_\epsilon$, $C_\epsilon$, $\mathsf C_M$, such that
the following holds. Let $\widehat f$ be a minimizer in \eqref{eq:oracle-estimator}, using the weights in \eqref{eq:oracle-generic-weights} with
\begin{align}
\tau_j=(1+\sigma_{\rho,j})
\sqrt{12\log\left(e\sum_{k=1}^{\mathrm J}n_k\right)},
\qquad j=1,\ldots,\mathrm J,
\label{eq:oracle-thresholds}
\end{align}
and the network parameters
\begin{align}
\sD&=\left\lceil
C_1\log\left(e\sum_{k=1}^{\mathrm J}n_k\right)
\right\rceil,
 \sW=\left\lceil C_2
\left\{
\frac{\mathsf{n}_{\rm eff}}
{\sqrt{12\log\left(e\sum_{k=1}^{\mathrm J}n_k\right)}}
\right\}^{\frac{K}{(2p+K)}}
\right\rceil,
 \sS=\left\lceil
C_3\sW\sD\right\rceil.
\label{eq:multisource-sparsity}
\end{align}
Then
\begin{align}
\mathbb E\|\widehat f-f^*\|_{\mathcal Q}^2
\le C_4\,
\left(
\sum_{j=1}^{\mathrm J}
\frac{n_j}{(1+\sigma_{\rho,j})\sigma_{\epsilon,j}^2}
\right)^{-2p/(2p+K)}
\left\{\log\left(e\sum_{k=1}^{\mathrm J}n_k\right)\right\}
^{3+p/(2p+K)}.
\nonumber
\end{align}
\end{theorem}

Condition \eqref{eq:oracle-effective-sample-size} requires the effective sample size to be at least $\sqrt{12\log(e\sum_{k=1}^{\mathrm J}n_k)}$, which grows slowly with the total sample size. It holds for sufficiently large total sample sizes whenever the covariate-shift parameters and noise variances are uniformly bounded above. The thresholds in \eqref{eq:oracle-thresholds} ensure that the expected reduction in each density ratio tends to zero, as stated in \eqref{eq:oracle-density-ratio-reduction}.

\Cref{thm:estimated-covariate-shift-rate} shows that source sample sizes, covariate shifts, and noise variances determine prediction accuracy through the effective sample size $\mathsf n_{\rm eff}$ defined in \eqref{eq:effective-sample-size}. Each source contributes in proportion to its sample size $n_j$ and inversely to its noise variance $\sigma_{\epsilon,j}^2$ and covariate shift factor $1+\sigma_{\rho,j}$. The rate exponent depends on the effective smoothness $p$ and effective dimension $K$. We next establish a matching minimax lower bound up to logarithmic factors.

\begin{theorem}
\label{thm:holder-covariate-shift-lower-bound}
Fix integers $L\ge1$ and $1\le K\le d$, and constants $p\ge1$, $R,M>0$, and $C_\epsilon\ge1$. Fix the number of sources $\mathrm J$ and, for each source $j$, the sample size $n_j$, noise variance $\sigma_{\epsilon,j}^2>0$, and covariate shift parameter $0\le\sigma_{\rho,j}<n_j$. Suppose that the effective sample size satisfies \eqref{eq:oracle-effective-sample-size}. There exists a constant $C_1>0$, depending only on $L,p,K,d,R,M$, such that
\begin{align}
\inf_{\widetilde f}
\sup_{\{\mathcal P_j\}_{j=1}^{\mathrm J},
      \mathcal Q,f^*,\epsilon}
\mathbb E\|\widetilde f-f^*\|_{\mathcal Q}^2
\ge C_1\,
\left(
\sum_{j=1}^{\mathrm J}
\frac{n_j}{(1+\sigma_{\rho,j})\sigma_{\epsilon,j}^2}
\right)^{-2p/(2p+K)}.
\nonumber
\end{align}
The supremum is over source and target distributions, regression functions, and noise distributions satisfying \Cref{ass:noise,ass:holder-shift,ass:hierarchical-composition} with the specified parameters. The infimum is over all estimators based on the source samples.
\end{theorem}
Together, \Cref{thm:estimated-covariate-shift-rate,thm:holder-covariate-shift-lower-bound} show that the oracle estimator is minimax optimal up to logarithmic factors when the conditions of both theorems hold and the ambient dimension and class parameters are fixed.

\section{Adaptive estimation}\label{sec:multisource-adaptive}
The oracle estimator requires knowledge of the density ratios $\rho_j$, the covariate shift parameters $\sigma_{\rho,j}$, and the noise variances $\sigma_{\epsilon,j}^2$. In this section, we propose an adaptive procedure that, under suitable assumptions, attains the oracle prediction rate. Recall that $\rho_j=d\mathcal Q/d\mathcal P_j$ is the density ratio of the target distribution to source $j$. We first assume that separate estimators $\widehat\rho_j$ of the density ratios $\rho_j$, $j=1,\ldots,\mathrm J$, are available and satisfy the following condition.

\begin{assumption}[\textbf{Density-ratio estimates}]
\label{ass:pooled-adaptive-source-ratios}
The collection $(\widehat\rho_j)_{j=1}^{\mathrm J}$ is independent of the labeled samples $(\mathcal D_j)_{j=1}^{\mathrm J}$ and satisfies
\begin{align}
\mathbb E\|\widehat\rho_j-\rho_j\|_{\mathcal P_j}^2 \le\varepsilon_j^2,
\qquad j=1,\ldots,\mathrm J, \label{eq:pooled-adaptive-source-errors}
\end{align}
for deterministic constants $\varepsilon_j>0$.
\end{assumption}

Density-ratio estimation has an extensive literature providing consistency and optimal convergence rates under suitable regularity conditions. \citet{kpotufe2017} establishes nearly optimal rates governed by the smoothness of the density ratio and the intrinsic dimension of the data. Using deep neural networks, \citet{zheng2026} establish nearly minimax mean-square error rates under smoothness and boundedness conditions and show that neural estimators can exploit low-dimensional manifold structure.

\begin{remark}
\label{rem:density-ratio-estimation}
The density-ratio estimators $\widehat\rho_j$ can be constructed by further splitting the available data. Alternatively, they can be constructed from additional unlabeled source and target data, consisting of covariates $X$ without responses $Y$. Let $m_j$ denote the number of unlabeled observations from source $j$. Under suitable regularity conditions and with sufficiently many unlabeled target observations, the density-ratio estimation error satisfies, up to logarithmic factors,
\begin{align}
\varepsilon_j^2\lesssim m_j^{-2\alpha/(2\alpha+d_\rho)},
\qquad j=1,\ldots,\mathrm J,
\nonumber
\end{align}
where $\alpha$ is the smoothness of the log-density ratio and $d_\rho$ is the effective dimension governing its estimation \citep{zheng2026}. When the covariates are supported on a low-dimensional manifold, $d_\rho$ reflects an embedding dimension determined by the manifold geometry and can be much smaller than the ambient dimension $d$, and potentially smaller than the effective regression dimension $K$.
\end{remark}

\subsection{Adaptive Estimation}\label{sec:multisource-adaptive-algorithm}
Let the total number of labeled observations be
\begin{align}
N=\sum_{j=1}^{\mathrm J}n_j.
\label{eq:multisource-total-sample-size}
\end{align}
For simplicity, suppose each $n_j\ge6$ is divisible by three. We split source $j$'s observations into three groups with index sets
\begin{align}
\mathcal I_{j,0} =\{1,\ldots,n_j/3\},  \quad 
\mathcal I_{j,1} =\{n_j/3+1,\ldots,2n_j/3\},   \quad 
\mathcal I_{j,2} =\{2n_j/3+1,\ldots,n_j\}. \nonumber
\end{align}
\Cref{alg:multisource-adaptive} uses these three groups as follows. Step 1 uses the first group to estimate response variances and construct the source weights and pooled density-ratio estimate. Step 2 uses the second group to fit candidate estimators with different weight thresholds. Step 3 uses the third group to compare all candidates using the same uncapped weights and select the candidate with the smallest validation loss.

For the candidate estimators in Step 2, we use the network class in \Cref{def: nn class}. Choose sufficiently large fixed constants $\sC_A,C_\sD,C_\sW,C_\sS$. For $N_{\widehat\nu}$ computed in Step 2, define the network depth by
\begin{align}
\sD=\lceil C_\sD\log(eN)\rceil,
\label{eq:adaptive-network-depth}
\end{align}
the width for each truncation level $\lambda$ by
\begin{align}
\sW_\lambda=\max\left\{d,
\left\lceil C_\sW
\left(\frac{N_{\widehat\nu}}{\lambda\{\log(eN)\}^3}\right)^{K/(2p+K)}
\right\rceil\right\},
\label{eq:adaptive-network-width}
\end{align}
and the sparsity by
\begin{align}
\sS_\lambda=\lceil C_\sS\sW_\lambda\log(eN)\rceil.
\label{eq:adaptive-network-sparsity}
\end{align}
The candidate network class is
\begin{align}
\mathcal F_\lambda=\mathcal F(\sD,\sW_\lambda,\sS_\lambda,\sC_A,M).
\label{eq:adaptive-network-class}
\end{align}

\begin{algorithm}[ht]
\caption{Adaptive estimation under multisource covariate shift}
\label{alg:multisource-adaptive}
\small
\setlength{\abovedisplayskip}{4pt}
\setlength{\belowdisplayskip}{4pt}
\setlength{\abovedisplayshortskip}{3pt}
\setlength{\belowdisplayshortskip}{3pt}
\textbf{Input:} Labeled samples $\{\mathcal D_j\}_{j=1}^{\mathrm J}$ and density-ratio estimates $\{\widehat\rho_j\}_{j=1}^{\mathrm J}$.
\par\smallskip
\textbf{Initialization:} Split each source into $\mathcal I_{j,0},\mathcal I_{j,1},\mathcal I_{j,2}$ as defined in \Cref{sec:multisource-adaptive-algorithm}.
\par\medskip
{\bf Step 1.} For each $j=1,\ldots,\mathrm J$, compute
\begin{align}
\overline Y_{j,0}
=\frac{3}{n_j}\sum_{i\in\mathcal I_{j,0}}Y_{ij},
\qquad
\widehat V_j=\frac{1}{n_j/3-1}\sum_{i\in\mathcal I_{j,0}}
(Y_{ij}-\overline Y_{j,0})^2.
\label{eq:adaptive-sample-variance-definition}
\end{align}
Set
\begin{align}
\widehat\nu_j=\max\{N^{-1},\widehat V_j\},
\qquad
\widehat\alpha_j=\frac{n_j/\widehat\nu_j}{\sum_{k=1}^{\mathrm J}n_k/\widehat\nu_k},
\qquad
\delta_N=\min\left\{\frac1{\log(eN)},\min_{1\le j\le\mathrm J}\widehat\alpha_j\right\}.
\nonumber
\end{align}
Compute
\begin{align}
\widehat\rho(x)=\left(\sum_{j=1}^{\mathrm J}
\frac{\widehat\alpha_j}{\max\{\delta_N,\widehat\rho_j(x)\}}\right)^{-1},
\qquad
\widehat r(x)=\max\{1,\widehat\rho(x)\}.
\nonumber
\end{align}
\par\medskip
{\bf Step 2.} Compute
\begin{align}
N_{\widehat\nu}=\sum_{j=1}^{\mathrm J}\frac{n_j}{3\widehat\nu_j}.
\nonumber
\end{align}
If $N_{\widehat\nu}<1$, return $\widehat f=0$. Otherwise, set
\begin{align}
\Lambda=\{2^\ell:\ell=0,\ldots,\lfloor\log_2N_{\widehat\nu}\rfloor\}.
\nonumber
\end{align}
For each $\lambda\in\Lambda$, use $\mathcal F_\lambda$ in \eqref{eq:adaptive-network-class} and choose a measurable minimizer
\begin{align}
\widehat f_\lambda\in\operatorname*{arg\,min}_{f\in\mathcal F_\lambda}
\sum_{j=1}^{\mathrm J}\sum_{i\in\mathcal I_{j,1}}
\frac{\min\{\widehat r(X_{ij}),\lambda\}}{\widehat\nu_j}
\bigl(Y_{ij}-f(X_{ij})\bigr)^2.
\label{eq:adaptive-training-estimator}
\end{align}
\par\medskip
{\bf Step 3.} Select
\begin{align}
\widehat\lambda\in\operatorname*{arg\,min}_{\lambda\in\Lambda}
\sum_{j=1}^{\mathrm J}\sum_{i\in\mathcal I_{j,2}}
\frac{\widehat r(X_{ij})}{\widehat\nu_j}
\bigl(Y_{ij}-\widehat f_\lambda(X_{ij})\bigr)^2,
\label{eq:adaptive-validation-selector}
\end{align}
breaking ties by choosing the smallest $\lambda$.
\par\medskip
\textbf{Output:} $\widehat f=\widehat f_{\widehat\lambda}$.
\end{algorithm}

In Step 1 of \Cref{alg:multisource-adaptive}, the law of large numbers implies that $\widehat\nu_j$ and $\widehat\alpha_j$ approach their population counterparts $\nu_j^*$ and $\alpha_j^*$, respectively, defined by
\begin{align}
\nu_j^*=\operatorname{Var}(Y_{ij}),
\qquad
\alpha_j^*=\frac{n_j/\nu_j^*}
{\sum_{k=1}^{\mathrm J}n_k/\nu_k^*},
\qquad j=1,\ldots,\mathrm J.
\label{eq:adaptive-population-source-weights}
\end{align}

Under \Cref{ass:noise,ass:hierarchical-composition}, and $\sigma_{\epsilon,j}^{2}\ge c_\epsilon$ as required in \eqref{eq:multisource-adaptive-sample-size},  it follows  that 
\begin{align}
\sigma_{\epsilon,j}^{2}
\le\nu_j^*
\le\left(1+\frac{M^{2}}{c_\epsilon}\right)\sigma_{\epsilon,j}^{2}.
\nonumber
\end{align}
Thus, $\widehat\nu_j^{-1}$ provides an estimate of $\sigma_{\epsilon,j}^{-2}$ up to a fixed multiplicative factor and sampling error, without first estimating $f^*$.
The population weights $\alpha_j^*$ in \eqref{eq:adaptive-population-source-weights} define the combined source distribution
\begin{align}
\overline{\mathcal P}=\sum_{j=1}^{\mathrm J}\alpha_j^*\mathcal P_j.
\label{eq:adaptive-source-mixture}
\end{align}
The density ratio of $\mathcal Q$ with respect to $\overline{\mathcal P}$ is
\begin{align}
\overline\rho
=\frac{d\mathcal Q}{d\overline{\mathcal P}}
=\left(\sum_{j=1}^{\mathrm J}\frac{\alpha_j^*}{\rho_j}\right)^{-1}
\qquad\mathcal Q\text{-almost surely}.
\label{eq:adaptive-mixture-density-ratio}
\end{align}

Step 1 of \Cref{alg:multisource-adaptive} constructs $\widehat\rho$ by substituting $\widehat\alpha_j$ for $\alpha_j^*$ and $\max\{\delta_N,\widehat\rho_j\}$ for $\rho_j$ in \eqref{eq:adaptive-mixture-density-ratio}. The positive floor $\delta_N$ prevents division by zero.

Step 2 of \Cref{alg:multisource-adaptive} uses $\widehat\rho$ to adjust for the difference between $\overline{\mathcal P}$ and $\mathcal Q$. Specifically, it weights each squared residual in \eqref{eq:adaptive-training-estimator} by
\begin{align}
\frac{\min\{\max\{1,\widehat\rho(X_{ij})\},\lambda\}}
{\widehat\nu_j},
\nonumber
\end{align}
where $\lambda$ limits the contribution of large values of $\widehat\rho(X_{ij})$.
Replacing $\widehat\nu_j$ by $\nu_j^*$ and $\widehat\rho$ by $\overline\rho$ in \eqref{eq:adaptive-training-estimator} and taking expectations over the training observations gives
\begin{align}
L_\lambda(f)
=\sum_{j=1}^{\mathrm J}\frac{n_j}{3\nu_j^*}
\mathbb E_j\!\left[
\min\{\max\{1,\overline\rho(X)\},\lambda\}
\{Y-f(X)\}^2
\right],
\label{eq:adaptive-population-criterion}
\end{align}
where $\mathbb E_j$ denotes expectation for an observation $(X,Y)$ from source $j$.
Moreover,
\begin{align}
\frac{L_\lambda(f)-L_\lambda(f^*)}
{\sum_{j=1}^{\mathrm J}n_j/(3\nu_j^*)}
&=\sum_{j=1}^{\mathrm J}\alpha_j^*
\mathbb E_{X\sim\mathcal P_j}\!\left[
\min\{\max\{1,\overline\rho(X)\},\lambda\}
\{f(X)-f^*(X)\}^2
\right]
\nonumber\\
&=\mathbb E_{X\sim\overline{\mathcal P}}\!\left[
\min\{\max\{1,\overline\rho(X)\},\lambda\}
\{f(X)-f^*(X)\}^2
\right],
\nonumber
\end{align}
where the first equality uses the regression model \eqref{eq:multisource-model}, the centered independent errors in \Cref{ass:noise}, and the definition of $\alpha_j^*$; and the second equality uses \eqref{eq:adaptive-source-mixture}.
Whenever $1\le\overline\rho(X)\le\lambda$, the weighting factor $\min\{\max\{1,\overline\rho(X)\},\lambda\}$ equals $\overline\rho(X)$. Weighting by $\overline\rho$ gives
\begin{align}
\mathbb E_{X\sim\overline{\mathcal P}}\!\left[
\overline\rho(X)\{f(X)-f^*(X)\}^2
\right]
=\|f-f^*\|_{\mathcal Q}^2,
\label{eq:adaptive-target-risk-identity}
\end{align}
where the equality follows from \eqref{eq:adaptive-mixture-density-ratio}. The identity \eqref{eq:adaptive-target-risk-identity} explains the connection to target prediction. Our analysis does not require $\overline\rho$ to be uniformly bounded above or bounded away from zero.

Step 3 of \Cref{alg:multisource-adaptive} evaluates every candidate using the same coefficients $\widehat r(X_{ij})/\widehat\nu_j$ in \eqref{eq:adaptive-validation-selector}. These coefficients do not depend on $\lambda$, so all candidates are compared using the same validation criterion.
Replacing $\widehat\nu_j$ by $\nu_j^*$ and $\widehat r$ by $\max\{1,\overline\rho\}$ in \eqref{eq:adaptive-validation-selector} and taking expectations over the $n_j/3$ validation observations from each source gives
\begin{align}
L_{\mathrm{val}}(f)
=\sum_{j=1}^{\mathrm J}\frac{n_j}{3\nu_j^*}
\mathbb E_j\!\left[
\max\{1,\overline\rho(X)\}\{Y-f(X)\}^2
\right].
\label{eq:adaptive-population-validation-criterion}
\end{align}
Moreover,
\begin{align}
\frac{L_{\mathrm{val}}(f)-L_{\mathrm{val}}(f^*)}
{\sum_{j=1}^{\mathrm J}n_j/(3\nu_j^*)}
&=\mathbb E_{X\sim\overline{\mathcal P}}\!\left[
\max\{1,\overline\rho(X)\}\{f(X)-f^*(X)\}^2
\right]
 \ge\|f-f^*\|_{\mathcal Q}^2,
\label{eq:adaptive-validation-risk-control}
\end{align}
where the equality uses the regression model \eqref{eq:multisource-model}, \Cref{ass:noise}, and the definition of $\overline{\mathcal P}$ in \eqref{eq:adaptive-source-mixture}; and the inequality uses $\max\{1,\overline\rho\}\ge\overline\rho$ and \eqref{eq:adaptive-target-risk-identity}.

Since $L_{\mathrm{val}}(f^*)$ is constant with respect to $\lambda$, minimizing $L_{\mathrm{val}}(\widehat f_\lambda)$ in \eqref{eq:adaptive-population-validation-criterion} also minimizes $L_{\mathrm{val}}(\widehat f_\lambda)-L_{\mathrm{val}}(f^*)$, which controls $\|\widehat f_\lambda-f^*\|_{\mathcal Q}^{2}$ by \eqref{eq:adaptive-validation-risk-control}.

\begin{remark}
Another approach is to pool all source observations directly and estimate the density ratio of the target distribution to the mixture
$
\sum_{j=1}^{\mathrm J}\frac{n_j}{N}\mathcal P_j,
$
where $N=\sum_{j=1}^{\mathrm J}n_j$ is the total number of labeled observations. However, this approach does not account for differences in the noise variances $\sigma_{\epsilon,j}^2$ across sources and can lead to suboptimal dependence on these variances in the prediction rate. On the other hand, our proposed adaptive procedure accounts for these differences through the weights $\widehat\nu_j^{-1}$ defined in Step 1 of \Cref{alg:multisource-adaptive}.
\end{remark}

\begin{theorem}
\label{thm:multisource-adaptive-rate}
Suppose \Cref{ass:noise,ass:holder-shift,ass:hierarchical-composition,ass:pooled-adaptive-source-ratios} hold. Assume that each $n_j$ is divisible by three and that
\begin{align}
\sigma_{\epsilon,j}^2\ge c_\epsilon,
\qquad n_j\ge\max\{6,C_0\log(eN)\},
\qquad j=1,\ldots,\mathrm J,
\label{eq:multisource-adaptive-sample-size}
\end{align}
where $c_\epsilon>0$ is fixed and $C_0>0$ is sufficiently large, depending only on $M,c_\epsilon,C_\epsilon$. For sufficiently large fixed constants $C_A,C_D,C_W,C_S$ and sufficiently large $N$, the adaptive estimator $\widehat f$ constructed in \Cref{sec:multisource-adaptive-algorithm} satisfies
\begin{align}
\mathbb E\|\widehat f-f^*\|_{\mathcal Q}^2
&\le C_1\mathsf n_{\rm eff}^{-2p/(2p+K)}
\{\log(eN)\}^{7p/(2p+K)}
+C_2\sum_{j=1}^{\mathrm J}
\frac{\varepsilon_j^2}{\alpha_j^*}+C_3N^{-3},
\nonumber
\end{align}
where $\mathsf n_{\rm eff}$ is defined in \eqref{eq:effective-sample-size}, $\alpha_j^*$ is the population source weight defined in \eqref{eq:adaptive-population-source-weights}, and $C_1,C_2,C_3>0$ depend only on $d,L,p,K,R,M,c_\epsilon,C_\epsilon$.
\end{theorem}

To briefly explain the rate in \Cref{thm:multisource-adaptive-rate}, consider the simplified setting where $\mathrm J$ is fixed and the population source weights $\alpha_j^*$ in \eqref{eq:adaptive-population-source-weights} are bounded away from zero. In this setting, the bound simplifies, up to logarithmic factors, to
\begin{align}
\mathbb E\|\widehat f-f^*\|_{\mathcal Q}^2
\lesssim \mathsf n_{\rm eff}^{-2p/(2p+K)}
+\sum_{j=1}^{\mathrm J}\varepsilon_j^2.
\nonumber
\end{align}

As discussed in \Cref{rem:density-ratio-estimation}, larger unlabeled samples or a smaller effective dimension $d_\rho$ can yield
\begin{align}
\sum_{j=1}^{\mathrm J}\varepsilon_j^2
=o\!\left(\mathsf n_{\rm eff}^{-2p/(2p+K)}\right).
\nonumber
\end{align}
In this case, the density-ratio estimation term is negligible relative to the regression estimation term, and the adaptive estimator attains the minimax rate up to logarithmic factors, matching the lower bound in \Cref{thm:holder-covariate-shift-lower-bound}.

\clearpage
\bibliographystyle{plainnat}
\bibliography{citations}

\clearpage
\appendix

\section{Technical results  for \Cref{sec:hierarchical-covariate-shift} }

Recall the total sample size $N$ defined in
\eqref{eq:multisource-total-sample-size}.
Let $\sD\ge2$, $\sW\ge d$, and $\sS\ge\max\{\sW,\sD\}$ be deterministic
integers, and let $\sC_A\ge1$ and $\sC_M\ge M$ be constants. Write
\begin{align}
\mathcal F = \mathcal F(\sD,\sW,\sS,\sC_A,\sC_M).
\label{eq:localization-network-class}
\end{align}

Set
\begin{align}
\neff
=\sum_{j=1}^{\rmJ}
\frac{n_j}{(1+\sigma_{\rho,j})\sigma_{\epsilon,j}^2}.
\label{eq:multisource-effective-sample-size}
\end{align}
Throughout this section we set
\begin{align}
w_j(x)
=\frac{1}{\sigma_{\epsilon,j}^2}
\left(
\frac{\rho_j(x)}
{(1+\sigma_{\rho,j})\sqrt{12\log(eN)}}
\wedge1\right).
\label{eq:estimated-weights}
\end{align}
For the weights in \eqref{eq:estimated-weights}, define the finite measure
\begin{align}
\nu(A)
=\frac{\sqrt{12\log(eN)}}{\neff}
\sum_{j=1}^{\rmJ}n_j\int_A w_j(x)\,d\mathcal P_j(x)
\qquad\text{for every Borel set }A\subseteq[0,1]^d.
\label{eq:multisource-population-measure}
\end{align}
Write $\|h\|_\nu^2=\int_{[0,1]^d}h(x)^2\,d\nu(x)$.

\subsection{Bounded-weight regression}
\begin{lemma}
\label{lem:weighted-quadratic-localization}
Suppose \Cref{ass:holder-shift,ass:hierarchical-composition}, and \eqref{eq:oracle-effective-sample-size}
hold and $\sigma_{\epsilon,j}^2\ge c_\epsilon>0$ for every
$j=1,\ldots,\rmJ$.
Let $\mathcal F$ be defined in \eqref{eq:localization-network-class},
and let $\neff$ and $\nu$ be defined in
\eqref{eq:multisource-effective-sample-size} and
\eqref{eq:multisource-population-measure}.
For every bounded measurable $f:[0,1]^d\to\mathbb R$, define
\begin{align}
\Phi(f)
=\frac{\sqrt{12\log(eN)}}{\neff}
\sum_{j=1}^{\rmJ}\sum_{i=1}^{n_j}
 w_j(X_{ij})\bigl(f(X_{ij})-f^*(X_{ij})\bigr)^2
-\|f-f^*\|_\nu^2.
\label{eq:weighted-quadratic-process}
\end{align}
There exists a constant $C_1>0$, depending only on $\sC_M$ and
$c_\epsilon$, such that, for every $0<\tau\le1$,
\begin{align}
&\mathbb E\left[
\sup_{f,g\in\mathcal F}
\left\{
|\Phi(f)-\Phi(g)|-\tau\|f-g\|_\nu^2
\right\}\right]
 \le
\frac{C_1\sS\sqrt{12\log(eN)}}{\tau\neff}
\left\{\sD\log\bigl((\sW+1)\sC_A\bigr)
+\log(eN)\right\}.
\nonumber
\end{align}
\end{lemma}

\begin{proof}
{\bf Step 1.}
Fix $0<\tau\le1$ and $f,g\in\mathcal F$. For each source $j$, define
\begin{align}
U_{j,f,g}(x)
=\frac{\sqrt{12\log(eN)}}{\neff}
 w_j(x)(f(x)-g(x))(f(x)+g(x)-2f^*(x)).
\label{eq:weighted-quadratic-increment}
\end{align}
The definitions \eqref{eq:weighted-quadratic-process},
\eqref{eq:weighted-quadratic-increment}, and
\eqref{eq:multisource-population-measure} give
\begin{align}
\Phi(f)-\Phi(g)
=\sum_{j=1}^{\rmJ}\sum_{i=1}^{n_j}
\left\{U_{j,f,g}(X_{ij})-\mathbb E U_{j,f,g}(X_{ij})\right\}.
\label{eq:weighted-quadratic-centered-sum}
\end{align}
By \eqref{eq:estimated-weights} and $\sigma_{\epsilon,j}^2\ge c_\epsilon$,
\begin{align}
0\le w_j(x)\le\sigma_{\epsilon,j}^{-2}\le c_\epsilon^{-1}.
\label{eq:weighted-quadratic-weight-bound}
\end{align}
Consequently,
\begin{align}
\left|U_{j,f,g}(X_{ij})-\mathbb E U_{j,f,g}(X_{ij})\right|
&\le |U_{j,f,g}(X_{ij})|+\mathbb E|U_{j,f,g}(X_{ij})|
\nonumber\\
&\le\frac{16\sC_M^2\sqrt{12\log(eN)}}{c_\epsilon\neff},
\label{eq:weighted-quadratic-centered-envelope}
\end{align}
where the first inequality uses the triangle inequality;
and the second inequality uses \eqref{eq:weighted-quadratic-increment},
\eqref{eq:weighted-quadratic-weight-bound}, $|f-g|\le2\sC_M$,
and $|f+g-2f^*|\le4\sC_M$.
Moreover,
\begin{align}
&\sum_{j=1}^{\rmJ}\sum_{i=1}^{n_j}
\operatorname{Var}\bigl(U_{j,f,g}(X_{ij})\bigr)
 \le\sum_{j=1}^{\rmJ}\sum_{i=1}^{n_j}
\mathbb E\bigl[U_{j,f,g}(X_{ij})^2\bigr]
\nonumber\\
&\quad\le16\sC_M^2
\left(\frac{\sqrt{12\log(eN)}}{\neff}\right)^2
\sum_{j=1}^{\rmJ}\sum_{i=1}^{n_j}
\mathbb E\left[w_j(X_{ij})^2(f(X_{ij})-g(X_{ij}))^2\right]
\nonumber\\
&\quad\le\frac{16\sC_M^2}{c_\epsilon}
\left(\frac{\sqrt{12\log(eN)}}{\neff}\right)^2
\sum_{j=1}^{\rmJ}\sum_{i=1}^{n_j}
\mathbb E\left[w_j(X_{ij})(f(X_{ij})-g(X_{ij}))^2\right]
\nonumber\\
&\quad=\frac{16\sC_M^2\sqrt{12\log(eN)}}{c_\epsilon\neff}
\|f-g\|_\nu^2,
\label{eq:weighted-quadratic-second-moment}
\end{align}
where  the second inequality uses \eqref{eq:weighted-quadratic-increment}
and $|f+g-2f^*|\le4\sC_M$;
the third inequality uses \eqref{eq:weighted-quadratic-weight-bound};
and the equality uses \eqref{eq:multisource-population-measure}
and the identical covariate distributions within each source $j$.
Applying \Cref{lem:bernstein-moment-bound} to the independent centered
summands in \eqref{eq:weighted-quadratic-centered-sum}, using
\eqref{eq:weighted-quadratic-centered-envelope} and
\eqref{eq:weighted-quadratic-second-moment}, gives a constant
$C_2>0$, depending only on $\sC_M$ and $c_\epsilon$, such that,
for every $t>0$,
\begin{align}
\mathbb P\left(
|\Phi(f)-\Phi(g)|
>C_2\|f-g\|_\nu
\sqrt{\frac{\sqrt{12\log(eN)}\,t}{\neff}}
+\frac{C_2\sqrt{12\log(eN)}\,t}{\neff}
\right)\le2e^{-t}.
\label{eq:weighted-quadratic-bernstein-tail}
\end{align}
For every $t>0$,
\begin{align}
&C_2\|f-g\|_\nu
\sqrt{\frac{\sqrt{12\log(eN)}\,t}{\neff}}
+\frac{C_2\sqrt{12\log(eN)}\,t}{\neff}
\nonumber\\
&\quad\le\tau\|f-g\|_\nu^2
+\left(\frac{C_2^2}{4\tau}+C_2\right)
\frac{\sqrt{12\log(eN)}\,t}{\neff}
\nonumber\\
&\quad\le\tau\|f-g\|_\nu^2
+\frac{C_3\sqrt{12\log(eN)}\,t}{\tau\neff},
\label{eq:weighted-quadratic-absorption}
\end{align}
where the first inequality uses Young's inequality;
and the second inequality uses $\tau\le1$, with a new constant
$C_3>0$ depending only on $\sC_M$ and $c_\epsilon$.
Combining \eqref{eq:weighted-quadratic-bernstein-tail} and
\eqref{eq:weighted-quadratic-absorption} yields
\begin{align}
\mathbb P\left(
|\Phi(f)-\Phi(g)|-\tau\|f-g\|_\nu^2
>\frac{C_3\sqrt{12\log(eN)}\,t}{\tau\neff}
\right)\le2e^{-t}.
\label{eq:weighted-quadratic-offset-tail}
\end{align}

{\bf Step 2.}
For every $f,g\in\mathcal F$ and $x\in[0,1]^d$,
\begin{align}
&\left|(f(x)-f^*(x))^2-(g(x)-f^*(x))^2\right|
\nonumber\\
&\quad=|f(x)-g(x)|\,|f(x)+g(x)-2f^*(x)|
\nonumber\\
&\quad\le4\sC_M|f(x)-g(x)|
\nonumber\\
&\quad\le4\sC_M\|f-g\|_\infty,
\label{eq:weighted-quadratic-pointwise-continuity}
\end{align}
where  
the first inequality uses $|f(x)|,|g(x)|\le \sC_M$ and
$|f^*(x)|\le M\le \sC_M$;
and the second inequality uses the definition of the uniform norm.
Moreover,
\begin{align}
&\int_{[0,1]^d}
\left|(f(x)-f^*(x))^2-(g(x)-f^*(x))^2\right|\,d\nu(x)
\nonumber\\
&\quad\le4\sC_M\|f-g\|_\infty\nu([0,1]^d)
\nonumber\\
&\quad\le4\sC_M\|f-g\|_\infty,
\label{eq:weighted-quadratic-population-continuity}
\end{align}
where the first inequality uses \eqref{eq:weighted-quadratic-pointwise-continuity};
and the second inequality uses \eqref{eq:multisource-measure-domination}
and $\mathcal Q([0,1]^d)=1$.
Simultaneously for all $f,g\in\mathcal F$,
\begin{align}
|\Phi(f)-\Phi(g)|
&\le\frac{\sqrt{12\log(eN)}}{\neff}
\sum_{j=1}^{\rmJ}\sum_{i=1}^{n_j}w_j(X_{ij})
\left|(f(X_{ij})-f^*(X_{ij}))^2-(g(X_{ij})-f^*(X_{ij}))^2\right|
\nonumber\\
&\quad+\int_{[0,1]^d}
\left|(f(x)-f^*(x))^2-(g(x)-f^*(x))^2\right|\,d\nu(x)
\nonumber\\
&\le4\sC_M\left\{
\frac{\sqrt{12\log(eN)}}{\neff}
\sum_{j=1}^{\rmJ}\sum_{i=1}^{n_j}w_j(X_{ij})+1
\right\}\|f-g\|_\infty
\nonumber\\
&\le4\sC_M\left\{
\frac{N\sqrt{12\log(eN)}}{c_\epsilon\neff}+1
\right\}\|f-g\|_\infty,
\label{eq:weighted-quadratic-continuity}
\end{align}
where the first inequality uses \eqref{eq:weighted-quadratic-process}
and the triangle inequality;
the second inequality uses \eqref{eq:weighted-quadratic-pointwise-continuity},
\eqref{eq:weighted-quadratic-population-continuity}, and $w_j\ge0$;
and the third inequality uses \eqref{eq:weighted-quadratic-weight-bound}
and \eqref{eq:multisource-total-sample-size}.

{\bf Step 3.}
By \Cref{lem:network-covering}, $\mathcal F$ is totally bounded
in the uniform norm, and there is a universal constant $C_4>0$
such that
\begin{align}
\log\mathcal N\left(N^{-1},\mathcal F,\|\cdot\|_\infty\right)
\le C_4S\left\{\sD\log\bigl((\sW+1)\sC_A\bigr)
+\log(4N)\right\}.
\label{eq:weighted-quadratic-entropy}
\end{align}
Moreover, $0\in\mathcal F$, so the class is nonempty.
Note that \begin{align}
\frac{1}{N}
&\le\frac{1}{c_\epsilon\neff}
 \le\frac{\sqrt{12\log(eN)}}{c_\epsilon\neff},
\label{eq:weighted-quadratic-discretization-scale}
\end{align}
where the first inequality uses \eqref{eq:multisource-effective-sample-size},
\eqref{eq:multisource-total-sample-size},
$1+\sigma_{\rho,j}\ge1$, and $\sigma_{\epsilon,j}^2\ge c_\epsilon$;
and the second inequality uses $\sqrt{12\log(eN)}\ge1$.
Therefore \begin{align}
&\mathbb E\left[
\sup_{f,g\in\mathcal F}
\left\{
|\Phi(f)-\Phi(g)|-\tau\|f-g\|_\nu^2
\right\}\right]
\nonumber\\
&\quad\le\frac{C_3\sqrt{12\log(eN)}}{\tau\neff}
\left\{2\log\mathcal N\left(N^{-1},\mathcal F,\|\cdot\|_\infty\right)
+\log2+1\right\}
 +8\sC_M\left\{
\frac{\sqrt{12\log(eN)}}{c_\epsilon\neff}
+\frac{1+\tau}{N}\right\}
\nonumber\\
&\quad\le\frac{C_5\sS\sqrt{12\log(eN)}}{\tau\neff}
\left\{\sD\log\bigl((\sW+1)\sC_A\bigr)
+\log(eN)\right\},
\nonumber
\end{align}
where the first inequality applies \Cref{lem:localized-supremum-bound}
using \eqref{eq:weighted-quadratic-offset-tail},
\eqref{eq:weighted-quadratic-continuity}, and
\eqref{eq:multisource-measure-domination}, with a covering radius equal
to $N^{-1}$;
and the second inequality uses \eqref{eq:weighted-quadratic-entropy},
\eqref{eq:weighted-quadratic-discretization-scale}, $\tau\le1$,
$S\ge1$, and
$\log(4N)\le2\log(eN)$.
This proves the lemma.
\end{proof}

\begin{lemma}
\label{lem:weighted-noise-localization}
Suppose \Cref{ass:noise,ass:holder-shift,}, and \eqref{eq:oracle-effective-sample-size}
hold and $\sigma_{\epsilon,j}^2\ge c_\epsilon>0$ for every
$j=1,\ldots,\rmJ$.
Let $\mathcal F$ be defined in \eqref{eq:localization-network-class},
and let $\neff$ and $\nu$ be defined in
\eqref{eq:multisource-effective-sample-size} and
\eqref{eq:multisource-population-measure}.
For every bounded measurable $f:[0,1]^d\to\mathbb R$, define
\begin{align}
\Psi(f)
=\frac{\sqrt{12\log(eN)}}{\neff}
\sum_{j=1}^{\rmJ}\sum_{i=1}^{n_j}
 w_j(X_{ij})\epsilon_{ij}f(X_{ij}).
\label{eq:weighted-noise-process}
\end{align}
There exists a constant $C_1>0$, depending only on $C_\epsilon$,
$\sC_M$, and $c_\epsilon$, such that, for every $0<\tau\le1$,
\begin{align}
&\mathbb E\left[
\sup_{f,g\in\mathcal F}
\left\{
|\Psi(f)-\Psi(g)|-\tau\|f-g\|_\nu^2
\right\}\right]
 \le
\frac{C_1\sS\sqrt{12\log(eN)}}{\tau\neff}
\left\{\sD\log\bigl((\sW+1)\sC_A\bigr)
+\log(eN)\right\}.
\nonumber
\end{align}
\end{lemma}

\begin{proof}
{\bf Step 1.}
Fix $0<\tau\le1$ and $f,g\in\mathcal F$. Define
\begin{align}
V_{ij}
=\frac{\sqrt{12\log(eN)}}{\neff}
 w_j(X_{ij})\epsilon_{ij}(f(X_{ij})-g(X_{ij})).
\label{eq:weighted-noise-increment}
\end{align}
The variables $V_{ij}$ are independent and centered by the sampling
model and \Cref{ass:noise}. By \eqref{eq:weighted-noise-process}
and \eqref{eq:weighted-noise-increment},
\begin{align}
\Psi(f)-\Psi(g)=\sum_{j=1}^{\rmJ}\sum_{i=1}^{n_j}V_{ij}.
\label{eq:weighted-noise-centered-sum}
\end{align}
The weights satisfy
\begin{align}
\sigma_{\epsilon,j}^2w_j(x)^2
&\le w_j(x)
 \le c_\epsilon^{-1},
\label{eq:weighted-noise-weight-bound}
\end{align}
where the first inequality uses $0\le\sigma_{\epsilon,j}^2w_j\le1$
from \eqref{eq:estimated-weights}; and the second inequality uses
$w_j\le\sigma_{\epsilon,j}^{-2}$ and $\sigma_{\epsilon,j}^2\ge c_\epsilon$.
For fixed $x\in[0,1]^d$ and $\lambda\in\mathbb R$,
\begin{align}
\mathbb E\left[e^{\lambda V_{ij}}\mid X_{ij}=x\right]
&=\mathbb E\exp\left\{
\lambda\frac{\sqrt{12\log(eN)}}{\neff}
w_j(x)(f(x)-g(x))\epsilon_{ij}\right\}
\nonumber\\
&\le\exp\left\{
\frac{C_\epsilon\sigma_{\epsilon,j}^2\lambda^2}{2}
\left(\frac{\sqrt{12\log(eN)}}{\neff}\right)^2
w_j(x)^2(f(x)-g(x))^2\right\},
\label{eq:weighted-noise-conditional-mgf}
\end{align}
where the equality uses \eqref{eq:weighted-noise-increment} and
the independence of $\epsilon_{ij}$ and $X_{ij}$;
and the inequality uses \Cref{ass:noise} with
$t=\lambda\sqrt{12\log(eN)}\,w_j(x)(f(x)-g(x))/\neff$.
For every $\lambda\in\mathbb R$ satisfying
\begin{align}
|\lambda|\le
\frac{\sqrt{c_\epsilon}\,\neff}
{2\sC_M\sqrt{C_\epsilon}\sqrt{12\log(eN)}},
\label{eq:weighted-noise-mgf-range}
\end{align}
we have
\begin{align}
\mathbb E e^{\lambda V_{ij}}
&=\mathbb E\left[\mathbb E\left[e^{\lambda V_{ij}}\mid X_{ij}\right]\right]
\nonumber\\
&\le\mathbb E\exp\left\{
\frac{C_\epsilon\sigma_{\epsilon,j}^2\lambda^2}{2}
\left(\frac{\sqrt{12\log(eN)}}{\neff}\right)^2
 w_j(X_{ij})^2(f(X_{ij})-g(X_{ij}))^2\right\}
\nonumber\\
&\le1+C_\epsilon\sigma_{\epsilon,j}^2\lambda^2
\left(\frac{\sqrt{12\log(eN)}}{\neff}\right)^2
\mathbb E\left[w_j(X_{ij})^2(f(X_{ij})-g(X_{ij}))^2\right]
\nonumber\\
&\le\exp\left\{
C_\epsilon\lambda^2
\left(\frac{\sqrt{12\log(eN)}}{\neff}\right)^2
\mathbb E\left[w_j(X_{ij})(f(X_{ij})-g(X_{ij}))^2\right]
\right\},
\label{eq:weighted-noise-individual-mgf}
\end{align}
where the equality uses the tower property of conditional expectation;
the first inequality uses \eqref{eq:weighted-noise-conditional-mgf};
the second inequality uses $e^u\le1+2u$ for $0\le u\le1$,
with the exponent at most $1/2$ by \eqref{eq:weighted-noise-mgf-range},
\eqref{eq:weighted-noise-weight-bound}, and $|f-g|\le2\sC_M$;
and the third inequality uses \eqref{eq:weighted-noise-weight-bound}
and $1+u\le e^u$.
Consequently,
\begin{align}
\log\mathbb E\exp\left\{\lambda\sum_{j=1}^{\rmJ}
\sum_{i=1}^{n_j}V_{ij}\right\}
&=\sum_{j=1}^{\rmJ}\sum_{i=1}^{n_j}
\log\mathbb E e^{\lambda V_{ij}}
\nonumber\\
&\le C_\epsilon\lambda^2
\left(\frac{\sqrt{12\log(eN)}}{\neff}\right)^2
\sum_{j=1}^{\rmJ}\sum_{i=1}^{n_j}
\mathbb E\left[w_j(X_{ij})(f(X_{ij})-g(X_{ij}))^2\right]
\nonumber\\
&=\frac{C_\epsilon\lambda^2\sqrt{12\log(eN)}}{\neff}
\|f-g\|_\nu^2,
\label{eq:weighted-noise-sum-mgf}
\end{align}
where the first equality uses independence;
the inequality uses \eqref{eq:weighted-noise-individual-mgf};
and the second equality uses \eqref{eq:multisource-population-measure}
and the identical covariate distributions within each source.
Applying Chernoff's bound to \eqref{eq:weighted-noise-sum-mgf}
for positive and negative $\lambda$ in \eqref{eq:weighted-noise-mgf-range},
and using \eqref{eq:weighted-noise-centered-sum}, gives a constant
$C_2>0$, depending only on $C_\epsilon$, $\sC_M$, and $c_\epsilon$,
such that, for every $t>0$,
\begin{align}
\mathbb P\left(
|\Psi(f)-\Psi(g)|
>C_2\|f-g\|_\nu
\sqrt{\frac{\sqrt{12\log(eN)}\,t}{\neff}}
+\frac{C_2\sqrt{12\log(eN)}\,t}{\neff}
\right)\le2e^{-t}.
\label{eq:weighted-noise-bernstein-tail}
\end{align}
For every $t>0$,
\begin{align}
&C_2\|f-g\|_\nu
\sqrt{\frac{\sqrt{12\log(eN)}\,t}{\neff}}
+\frac{C_2\sqrt{12\log(eN)}\,t}{\neff}
\nonumber\\
&\quad\le\tau\|f-g\|_\nu^2
+\left(\frac{C_2^2}{4\tau}+C_2\right)
\frac{\sqrt{12\log(eN)}\,t}{\neff}
\nonumber\\
&\quad\le\tau\|f-g\|_\nu^2
+\frac{C_3\sqrt{12\log(eN)}\,t}{\tau\neff},
\label{eq:weighted-noise-absorption}
\end{align}
where the first inequality uses Young's inequality;
and the second inequality uses $\tau\le1$, with a new constant
$C_3>0$ depending only on $C_\epsilon$, $\sC_M$, and $c_\epsilon$.
Combining \eqref{eq:weighted-noise-bernstein-tail} and
\eqref{eq:weighted-noise-absorption} yields
\begin{align}
\mathbb P\left(
|\Psi(f)-\Psi(g)|-\tau\|f-g\|_\nu^2
>\frac{C_3\sqrt{12\log(eN)}\,t}{\tau\neff}
\right)\le2e^{-t}.
\label{eq:weighted-noise-offset-tail}
\end{align}

{\bf Step 2.}
Simultaneously for all $f,g\in\mathcal F$,
\begin{align}
|\Psi(f)-\Psi(g)|
&=\left|\frac{\sqrt{12\log(eN)}}{\neff}
\sum_{j=1}^{\rmJ}\sum_{i=1}^{n_j}
 w_j(X_{ij})\epsilon_{ij}(f(X_{ij})-g(X_{ij}))\right|
\nonumber\\
&\le\frac{\sqrt{12\log(eN)}}{\neff}
\sum_{j=1}^{\rmJ}\sum_{i=1}^{n_j}
 w_j(X_{ij})|\epsilon_{ij}|\,|f(X_{ij})-g(X_{ij})|
\nonumber\\
&\le\left\{\frac{\sqrt{12\log(eN)}}{\neff}
\sum_{j=1}^{\rmJ}\sum_{i=1}^{n_j}
 w_j(X_{ij})|\epsilon_{ij}|\right\}\|f-g\|_\infty,
\label{eq:weighted-noise-continuity}
\end{align}
where the equality uses \eqref{eq:weighted-noise-process};
the first inequality uses the triangle inequality and $w_j\ge0$;
and the second inequality uses the uniform norm.
Moreover,
\begin{align}
&\mathbb E\left[\frac{\sqrt{12\log(eN)}}{\neff}
\sum_{j=1}^{\rmJ}\sum_{i=1}^{n_j}
 w_j(X_{ij})|\epsilon_{ij}|\right]
 =\frac{\sqrt{12\log(eN)}}{\neff}
\sum_{j=1}^{\rmJ}\sum_{i=1}^{n_j}
\mathbb E[w_j(X_{ij})]\mathbb E[|\epsilon_{ij}|]
\nonumber\\
&\quad\le\frac{\sqrt{12\log(eN)}}{\neff}
\sum_{j=1}^{\rmJ}\sum_{i=1}^{n_j}
\sigma_{\epsilon,j}\mathbb E[w_j(X_{ij})]
 \le\frac{\sqrt{12\log(eN)}}{\neff}
\sum_{j=1}^{\rmJ}\frac{n_j}{\sigma_{\epsilon,j}}
 \le\frac{N\sqrt{12\log(eN)}}{\sqrt{c_\epsilon}\,\neff},
\label{eq:weighted-noise-expected-coefficient}
\end{align}
where  
the first inequality uses  
$\mathbb E[\epsilon_{ij}^2]=\sigma_{\epsilon,j}^2$ from \Cref{ass:noise};
the second inequality uses $w_j\le\sigma_{\epsilon,j}^{-2}$
from \eqref{eq:estimated-weights};
and the third inequality uses $\sigma_{\epsilon,j}^2\ge c_\epsilon$
and \eqref{eq:multisource-total-sample-size}.

{\bf Step 3.}
By \Cref{lem:network-covering}, there is a deterministic
$N^{-1}$-net
$\mathcal F_0\subset\mathcal F$ with cardinality
\begin{align}
\log| \mathcal F_0| \le C_4\sS\left\{\sD\log\bigl((\sW+1)\sC_A\bigr)
+\log(4N)\right\},
\label{eq:weighted-noise-entropy}
\end{align}
where   $C_4>0$ is  a universal constant.

For each ordered pair $(f,g)\in\mathcal F_0\times\mathcal F_0$,
consider the random variable
\begin{align}
|\Psi(f)-\Psi(g)|-\tau\|f-g\|_\nu^2.
\nonumber
\end{align}
There are $|\mathcal F_0|^2$ such random variables, and each satisfies
\eqref{eq:weighted-noise-offset-tail}. The pairs with $f=g$ give zero,
so their maximum is nonnegative. Applying
\Cref{lem:finite-maximum-exponential-tail} therefore gives
\begin{align}
&\mathbb E\left[
\max_{f,g\in\mathcal F_0}
\left\{|\Psi(f)-\Psi(g)|-\tau\|f-g\|_\nu^2\right\}
\right]
 \le\frac{C_3\sqrt{12\log(eN)}}{\tau\neff}
\{2\log|\mathcal F_0|+\log2+1\}.
\label{eq:weighted-noise-net-maximum}
\end{align}
For each $f,g\in\mathcal F$, choose $f_0,g_0\in\mathcal F_0$
with $\|f-f_0\|_\infty,\|g-g_0\|_\infty
\le N^{-1}$. Then
\begin{align}
\left|\|f-g\|_\nu^2-\|f_0-g_0\|_\nu^2\right|
&\le\bigl(\|f-g\|_\nu+\|f_0-g_0\|_\nu\bigr)
\bigl(\|f-f_0\|_\nu+\|g-g_0\|_\nu\bigr)
\nonumber\\
&\le4\sC_M\left(\|f-f_0\|_\infty+\|g-g_0\|_\infty\right)
 \le\frac{8\sC_M}{N},
\label{eq:weighted-noise-penalty-comparison}
\end{align}
where 
the second inequality uses the uniform bound $\sC_M$ and
\eqref{eq:multisource-measure-domination}, which implies $\nu([0,1]^d)\le1$;
and the third inequality uses the net approximation.
Therefore, \eqref{eq:weighted-noise-penalty-comparison} gives
\begin{align}
-\tau\|f-g\|_\nu^2
\le-\tau\|f_0-g_0\|_\nu^2
+\frac{8\tau \sC_M}{N}.
\label{eq:weighted-noise-negative-penalty}
\end{align}

Simultaneously for all $f,g\in\mathcal F$,
\begin{align}
|\Psi(f)-\Psi(g)|
&\le|\Psi(f_0)-\Psi(g_0)|
+|\Psi(f)-\Psi(f_0)|+|\Psi(g)-\Psi(g_0)|
\nonumber\\
&\le|\Psi(f_0)-\Psi(g_0)|
+\frac{2\sqrt{12\log(eN)}}
{\neff N}
\sum_{j=1}^{\rmJ}\sum_{i=1}^{n_j}
 w_j(X_{ij})|\epsilon_{ij}|,
\label{eq:weighted-noise-process-comparison}
\end{align}
where the second inequality applies \eqref{eq:weighted-noise-continuity}
to $(f,f_0)$ and $(g,g_0)$ and uses
$\|f-f_0\|_\infty+\|g-g_0\|_\infty
\le2/N$.

Therefore,
\begin{align}
&|\Psi(f)-\Psi(g)|-\tau\|f-g\|_\nu^2
\nonumber\\
&\quad\le|\Psi(f_0)-\Psi(g_0)|-\tau\|f_0-g_0\|_\nu^2
 +
\frac{2\sqrt{12\log(eN)}}
{\neff N}
\sum_{j=1}^{\rmJ}\sum_{i=1}^{n_j}
 w_j(X_{ij})|\epsilon_{ij}|
+\frac{8\tau \sC_M}{N}
\nonumber\\
&\quad\le\max_{u,v\in\mathcal F_0}
\left\{|\Psi(u)-\Psi(v)|-\tau\|u-v\|_\nu^2\right\}
 +
\frac{2\sqrt{12\log(eN)}}
{\neff N}
\sum_{j=1}^{\rmJ}\sum_{i=1}^{n_j}
 w_j(X_{ij})|\epsilon_{ij}|
+\frac{8\tau \sC_M}{N},
\label{eq:weighted-noise-net-comparison}
\end{align}
where the first inequality follows from \eqref{eq:weighted-noise-process-comparison}
and \eqref{eq:weighted-noise-negative-penalty}.
Taking the supremum and then expectations gives
\begin{align}
&\mathbb E\left[
\sup_{f,g\in\mathcal F}
\left\{|\Psi(f)-\Psi(g)|-\tau\|f-g\|_\nu^2\right\}
\right]
\nonumber\\
&\quad\le\frac{C_3\sqrt{12\log(eN)}}{\tau\neff}
\{2\log|\mathcal F_0|+\log2+1\}
+\frac{2\sqrt{12\log(eN)}}{\sqrt{c_\epsilon}\,\neff}
+\frac{8\tau \sC_M}{N}
\nonumber\\
&\quad\le\frac{C_5\sS\sqrt{12\log(eN)}}{\tau\neff}
\left\{\sD\log\bigl((\sW+1)\sC_A\bigr)
+\log(eN)\right\},
\nonumber
\end{align}
where the first inequality uses
\eqref{eq:weighted-noise-net-comparison},
\eqref{eq:weighted-noise-net-maximum}, and
\eqref{eq:weighted-noise-expected-coefficient};
and the second inequality uses \eqref{eq:weighted-noise-entropy},
\eqref{eq:weighted-quadratic-discretization-scale}, $\tau\le1$,
$S\ge1$, and
$\log(4N)\le2\log(eN)$,
with a new constant $C_5>0$ depending only on $C_\epsilon$, $\sC_M$,
and $c_\epsilon$.
This proves the lemma.
\end{proof}
 
\subsection{Estimated importance weights}

\begin{lemma}
\label{lem:subgaussian-truncation-bound}
Suppose \Cref{ass:holder-shift} holds. For every $j=1,\ldots,\rmJ$,
\begin{align}
\mathbb E_{X\sim\mathcal P_j}
\left[
\left\{\rho_j(X)-(1+\sigma_{\rho,j})
\sqrt{12\log(eN)}\right\}_+
\right]
\le N^{-3}.
\label{eq:multisource-clipping-bound}
\end{align}
If $\sigma_{\rho,j}=0$, the expectation in \eqref{eq:multisource-clipping-bound} is zero.
\end{lemma}

\begin{proof}
{\bf Step 1.}
Fix $j\in\{1,\ldots,\rmJ\}$. If $\sigma_{\rho,j}=0$,
\Cref{ass:holder-shift} implies
$\mathbb P_{X\sim\mathcal P_j}(\rho_j(X)=1)=1$.
Since $\sqrt{12\log(eN)}\ge1$, the expectation in
\eqref{eq:multisource-clipping-bound} is zero.
Suppose henceforth that $\sigma_{\rho,j}>0$.
For $u>0$, set $t=u/\sigma_{\rho,j}^2$. Then
\begin{align}
\mathbb P_{X\sim\mathcal P_j}\bigl(\rho_j(X)>1+u\bigr)
&\le e^{-tu}\mathbb E_{X\sim\mathcal P_j}
\left[\exp\{t(\rho_j(X)-1)\}\right]
\nonumber\\
&\le\exp\left\{-tu+\frac{\sigma_{\rho,j}^2t^2}{2}\right\}
 =\exp\left\{-\frac{u^2}{2\sigma_{\rho,j}^2}\right\},
\label{eq:subgaussian-ratio-tail}
\end{align}
where the first inequality uses Markov's inequality; the second
inequality uses \Cref{ass:holder-shift}; and the equality uses
$t=u/\sigma_{\rho,j}^2$.

{\bf Step 2.}
The clipping loss satisfies
\begin{align}
&\mathbb E_{X\sim\mathcal P_j}
\left[\left\{\rho_j(X)-(1+\sigma_{\rho,j})
\sqrt{12\log(eN)}\right\}_+\right]
\nonumber\\
&\quad\le\mathbb E_{X\sim\mathcal P_j}
\left[\left\{\rho_j(X)-1-\sigma_{\rho,j}
\sqrt{12\log(eN)}\right\}_+\right]
 =\sigma_{\rho,j}\int_{\sqrt{12\log(eN)}}^\infty
\mathbb P_{X\sim\mathcal P_j}
\bigl(\rho_j(X)>1+\sigma_{\rho,j}v\bigr)\,dv
\nonumber\\
&\quad\le\sigma_{\rho,j}\int_{\sqrt{12\log(eN)}}^\infty
 e^{-v^2/2}\,dv
 \le\frac{\sigma_{\rho,j}}{\sqrt{12\log(eN)}}
\int_{\sqrt{12\log(eN)}}^\infty v e^{-v^2/2}\,dv
 =\frac{\sigma_{\rho,j}}{\sqrt{12\log(eN)}}(eN)^{-6}
\nonumber\\
&\quad\le N(eN)^{-6}
 \le N^{-3},
\nonumber
\end{align}
where the first inequality uses $\sqrt{12\log(eN)}\ge1$;
the first equality uses the tail integral formula of expectation;
the second inequality uses \eqref{eq:subgaussian-ratio-tail};
the third inequality uses $v\ge\sqrt{12\log(eN)}$ on the integration domain;
the second equality evaluates the integral;
the fourth inequality uses $\sqrt{12\log(eN)}\ge1$ and
$\sigma_{\rho,j}<n_j\le N$ from \Cref{ass:holder-shift} and
\eqref{eq:multisource-total-sample-size};
and the fifth inequality uses $N\ge2$.
\end{proof}

\begin{lemma}
\label{lem:multisource-target-risk-comparison}
Suppose \Cref{ass:holder-shift,ass:hierarchical-composition}, and \eqref{eq:oracle-effective-sample-size}
hold, and let $\sC_M\ge M$. For every measurable
$f:[0,1]^d\to\mathbb R$ with $\|f\|_\infty\le \sC_M$,
\begin{align}
\|f-f^*\|_\nu^2
\le\|f-f^*\|_{\mathcal Q}^2
\le\|f-f^*\|_\nu^2
+4\sC_M^2N^{-3}.
\nonumber
\end{align}
\end{lemma}

\begin{proof}
{\bf Step 1.}
For every Borel set $A\subseteq[0,1]^d$,
\begin{align}
\nu(A)
&=\frac{1}{\neff}\sum_{j=1}^{\rmJ}
\frac{n_j}{(1+\sigma_{\rho,j})\sigma_{\epsilon,j}^2}
\int_A\left\{\rho_j(x)\wedge
(1+\sigma_{\rho,j})\sqrt{12\log(eN)}\right\}\,d\mathcal P_j(x)
\nonumber\\
&\le\frac{1}{\neff}\sum_{j=1}^{\rmJ}
\frac{n_j}{(1+\sigma_{\rho,j})\sigma_{\epsilon,j}^2}
\int_A\rho_j(x)\,d\mathcal P_j(x)
\nonumber\\
&=\mathcal Q(A),
\label{eq:multisource-measure-domination}
\end{align}
where the first equality uses \eqref{eq:multisource-population-measure}
and \eqref{eq:estimated-weights}; the inequality uses $a\wedge b\le a$;
and the second equality uses $\rho_j=d\mathcal Q/d\mathcal P_j$
and \eqref{eq:multisource-effective-sample-size}.
By \eqref{eq:multisource-measure-domination}, integration of the
nonnegative function $(f-f^*)^2$ gives
\begin{align}
\|f-f^*\|_\nu^2\le\|f-f^*\|_{\mathcal Q}^2.
\nonumber
\end{align}

{\bf Step 2.}
The normalization gives
\begin{align}
1
 =\frac{1}{\neff}\sum_{j=1}^{\rmJ}
\frac{n_j}{(1+\sigma_{\rho,j})\sigma_{\epsilon,j}^2}
 =\frac{1}{\neff}\sum_{j=1}^{\rmJ}
\frac{n_j}{(1+\sigma_{\rho,j})\sigma_{\epsilon,j}^2}
\mathbb E_{X\sim\mathcal P_j}[\rho_j(X)],
\label{eq:multisource-mass-normalization}
\end{align}
where the first equality uses \eqref{eq:multisource-effective-sample-size};
and the second equality uses
$\mathbb E_{X\sim\mathcal P_j}[\rho_j(X)]=1$.
Moreover,
\begin{align}
\nu([0,1]^d)
&=\frac{\sqrt{12\log(eN)}}{\neff}
\sum_{j=1}^{\rmJ}n_j
\mathbb E_{X\sim\mathcal P_j}[w_j(X)]
\nonumber\\
&=\frac{1}{\neff}\sum_{j=1}^{\rmJ}
\frac{n_j}{(1+\sigma_{\rho,j})\sigma_{\epsilon,j}^2}
\mathbb E_{X\sim\mathcal P_j}
\left[\rho_j(X)\wedge
(1+\sigma_{\rho,j})\sqrt{12\log(eN)}\right],
\label{eq:multisource-clipped-mass}
\end{align}
where the first equality uses \eqref{eq:multisource-population-measure};
and the second equality uses \eqref{eq:estimated-weights}
and $b(a/b\wedge1)=a\wedge b$ for $b>0$.
Therefore,
\begin{align}
1-\nu([0,1]^d)
&=\frac{1}{\neff}\sum_{j=1}^{\rmJ}
\frac{n_j}{(1+\sigma_{\rho,j})\sigma_{\epsilon,j}^2}
\mathbb E_{X\sim\mathcal P_j}
\left[\rho_j(X)-\left\{\rho_j(X)\wedge
(1+\sigma_{\rho,j})\sqrt{12\log(eN)}\right\}\right]
\nonumber\\
&=\frac{1}{\neff}\sum_{j=1}^{\rmJ}
\frac{n_j}{(1+\sigma_{\rho,j})\sigma_{\epsilon,j}^2}
\mathbb E_{X\sim\mathcal P_j}
\left[\left\{\rho_j(X)-(1+\sigma_{\rho,j})
\sqrt{12\log(eN)}\right\}_+\right]
\nonumber\\
&\le\frac{1}{\neff}\sum_{j=1}^{\rmJ}
\frac{n_j}{(1+\sigma_{\rho,j})\sigma_{\epsilon,j}^2}
N^{-3}
=N^{-3},
\label{eq:multisource-missing-mass}
\end{align}
where the first equality subtracts \eqref{eq:multisource-clipped-mass}
from \eqref{eq:multisource-mass-normalization} and uses linearity of expectation;
the second equality uses $a-(a\wedge b)=(a-b)_+$;
the inequality uses \eqref{eq:multisource-clipping-bound};
and the third equality uses \eqref{eq:multisource-effective-sample-size}.
Consequently,
\begin{align}
\|f-f^*\|_{\mathcal Q}^2-\|f-f^*\|_\nu^2
&=\int_{[0,1]^d}(f(x)-f^*(x))^2\,d(\mathcal Q-\nu)(x)
\nonumber\\
&\le4\sC_M^2\{1-\nu([0,1]^d)\}
 \le4\sC_M^2N^{-3},
\nonumber
\end{align}
where the equality uses the definitions of the squared norms;
the first inequality uses \eqref{eq:multisource-measure-domination}
and $|f-f^*|\le \sC_M+M\le2\sC_M$;
and the second inequality uses \eqref{eq:multisource-missing-mass}.
\end{proof}

\begin{lemma}
\label{lem:capped-weighted-oracle-inequality}
Suppose \Cref{ass:noise,ass:holder-shift,ass:hierarchical-composition}, and \eqref{eq:oracle-effective-sample-size}
hold and $\sigma_{\epsilon,j}^2\ge c_\epsilon>0$ for every
$j=1,\ldots,\rmJ$.
Let $\sD\ge2$, $\sW\ge d$, and $\sS\ge\max\{\sW,\sD\}$ be deterministic
integers, and let $\sC_A\ge1$ and $\sC_M\ge M$.
Let $\mathcal F$ be defined in \eqref{eq:localization-network-class},
let $w_j$ be defined in \eqref{eq:estimated-weights}, and let
$\neff$ be defined in \eqref{eq:multisource-effective-sample-size}.
Let $\widehat f$ be a measurable minimizer satisfying
\begin{align}
\widehat f\in\operatorname*{arg\,min}_{f\in\mathcal F}
\sum_{j=1}^{\rmJ}\sum_{i=1}^{n_j}
 w_j(X_{ij})\bigl(Y_{ij}-f(X_{ij})\bigr)^2.
\label{eq:holder-weighted-erm}
\end{align}
There exists a constant $C_1>0$, depending only on $C_\epsilon$,
$\sC_M$, and $c_\epsilon$, such that
\begin{align}
\mathbb E\|\widehat f-f^*\|_{\mathcal Q}^2
&\le3\inf_{f\in\mathcal F}\|f-f^*\|_{\mathcal Q}^2
\nonumber\\
&\quad+\frac{C_1\sS\sqrt{12\log(eN)}}{\neff}
\left\{\sD\log\bigl((\sW+1)\sC_A\bigr)
+\log(eN)\right\}
 +4\sC_M^2N^{-3}.
\nonumber
\end{align}
\end{lemma}

\begin{proof}
{\bf Step 1.}
Fix a deterministic $f_0\in\mathcal F$.
The model \eqref{eq:multisource-model} gives, for every $f\in\mathcal F$,
\begin{align}
(Y_{ij}-f(X_{ij}))^2-(Y_{ij}-f^*(X_{ij}))^2
=(f(X_{ij})-f^*(X_{ij}))^2
-2\epsilon_{ij}(f(X_{ij})-f^*(X_{ij})).
\label{eq:capped-loss-expansion}
\end{align}
For every deterministic $f\in\mathcal F$, the population excess loss satisfies
\begin{align}
&\frac{\sqrt{12\log(eN)}}{\neff}
\sum_{j=1}^{\rmJ}\sum_{i=1}^{n_j}
\mathbb E\left[w_j(X_{ij})
\left\{(Y_{ij}-f(X_{ij}))^2-(Y_{ij}-f^*(X_{ij}))^2\right\}\right]
\nonumber\\
&\quad=\frac{\sqrt{12\log(eN)}}{\neff}
\sum_{j=1}^{\rmJ}\sum_{i=1}^{n_j}
\mathbb E\left[w_j(X_{ij})(f(X_{ij})-f^*(X_{ij}))^2\right]
\nonumber\\
&\quad=\|f-f^*\|_\nu^2,
\nonumber
\end{align}
where the first equality uses \eqref{eq:capped-loss-expansion}
and the centered errors independent of the covariates in \Cref{ass:noise};
and the second equality uses \eqref{eq:multisource-population-measure}.
The minimizing property \eqref{eq:holder-weighted-erm} gives
\begin{align}
\sum_{j=1}^{\rmJ}\sum_{i=1}^{n_j}
 w_j(X_{ij})(Y_{ij}-\widehat f(X_{ij}))^2
\le\sum_{j=1}^{\rmJ}\sum_{i=1}^{n_j}
 w_j(X_{ij})(Y_{ij}-f_0(X_{ij}))^2.
\label{eq:capped-minimizing-property}
\end{align}
Subtracting the loss at $f^*$ from both sides of
\eqref{eq:capped-minimizing-property}, applying \eqref{eq:capped-loss-expansion},
and multiplying by $\sqrt{12\log(eN)}/\neff$ yields
\begin{align}
&\frac{\sqrt{12\log(eN)}}{\neff}
\sum_{j=1}^{\rmJ}\sum_{i=1}^{n_j}
 w_j(X_{ij})\left\{
(\widehat f(X_{ij})-f^*(X_{ij}))^2
-(f_0(X_{ij})-f^*(X_{ij}))^2\right\}
\nonumber\\
&\quad\le\frac{2\sqrt{12\log(eN)}}{\neff}
\sum_{j=1}^{\rmJ}\sum_{i=1}^{n_j}
 w_j(X_{ij})\epsilon_{ij}(\widehat f(X_{ij})-f_0(X_{ij})).
\label{eq:capped-basic-empirical-inequality}
\end{align}
Using \eqref{eq:weighted-quadratic-process} and
\eqref{eq:weighted-noise-process} to rewrite
\eqref{eq:capped-basic-empirical-inequality}, we obtain
\begin{align}
\|\widehat f-f^*\|_\nu^2
\le\|f_0-f^*\|_\nu^2
+\Phi(f_0)-\Phi(\widehat f)
+2\{\Psi(\widehat f)-\Psi(f_0)\}.
\label{eq:capped-basic-inequality}
\end{align}
{\bf Step 2.}
The quadratic increment satisfies
\begin{align}
&\mathbb E\left[
\Phi(f_0)-\Phi(\widehat f)-\frac{1}{12}\|\widehat f-f_0\|_\nu^2
\right]
 \le\mathbb E\left[
\sup_{f,g\in\mathcal F}
\left\{|\Phi(f)-\Phi(g)|-\frac{1}{12}\|f-g\|_\nu^2\right\}
\right]
\nonumber\\
&\quad\le\frac{C_2\sS\sqrt{12\log(eN)}}{\neff}
\left\{\sD\log\bigl((\sW+1)\sC_A\bigr)
+\log(eN)\right\},
\label{eq:capped-quadratic-increment-bound}
\end{align}
where  
and the second inequality applies \Cref{lem:weighted-quadratic-localization}
with $\tau=1/12$, for a constant $C_2>0$ depending only on $\sC_M$ and
$c_\epsilon$.
Similarly,
\begin{align}
&\mathbb E\left[
\Psi(\widehat f)-\Psi(f_0)-\frac{1}{12}\|\widehat f-f_0\|_\nu^2
\right]
 \le\mathbb E\left[
\sup_{f,g\in\mathcal F}
\left\{|\Psi(f)-\Psi(g)|-\frac{1}{12}\|f-g\|_\nu^2\right\}
\right]
\nonumber\\
&\quad\le\frac{C_3\sS\sqrt{12\log(eN)}}{\neff}
\left\{\sD\log\bigl((\sW+1)\sC_A\bigr)
+\log(eN)\right\},
\label{eq:capped-noise-increment-bound}
\end{align}
where  
  the second inequality applies \Cref{lem:weighted-noise-localization}
with $\tau=1/12$, for a constant $C_3>0$ depending only on
$C_\epsilon$, $\sC_M$, and $c_\epsilon$.
Combining the two increments gives
\begin{align}
&\mathbb E\left[
\Phi(f_0)-\Phi(\widehat f)
+2\{\Psi(\widehat f)-\Psi(f_0)\}
-\frac14\|\widehat f-f_0\|_\nu^2\right]
\nonumber\\
&\quad=\mathbb E\left[
\Phi(f_0)-\Phi(\widehat f)
-\frac1{12}\|\widehat f-f_0\|_\nu^2\right]
\nonumber\\
&\qquad\quad+2\mathbb E\left[
\Psi(\widehat f)-\Psi(f_0)
-\frac1{12}\|\widehat f-f_0\|_\nu^2\right]
\nonumber\\
&\quad\le\frac{(C_2+2C_3)\sS\sqrt{12\log(eN)}}{\neff}
\left\{\sD\log\bigl((\sW+1)\sC_A\bigr)
+\log(eN)\right\},
\label{eq:capped-combined-increment-bound}
\end{align}
where the equality uses linearity of expectation and $1/12+2/12=1/4$;
and the inequality applies \eqref{eq:capped-quadratic-increment-bound}
once and \eqref{eq:capped-noise-increment-bound} twice.
Set $C_4=C_2+2C_3$, which depends only on $C_\epsilon$, $\sC_M$,
and $c_\epsilon$.
The squared distance satisfies
\begin{align}
\|\widehat f-f_0\|_\nu^2
&=\| (\widehat f-f^*)-(f_0-f^*)\|_\nu^2
 \le2\|\widehat f-f^*\|_\nu^2+2\|f_0-f^*\|_\nu^2. 
\label{eq:capped-comparator-distance}
\end{align}
Consequently,
\begin{align}
&\mathbb E\|\widehat f-f^*\|_\nu^2
\nonumber\\
 &\le\|f_0-f^*\|_\nu^2
+\frac14\mathbb E\|\widehat f-f_0\|_\nu^2
 +\frac{C_4\sS\sqrt{12\log(eN)}}{\neff}
\left\{\sD\log\bigl((\sW+1)\sC_A\bigr)
+\log(eN)\right\}
 \nonumber\\
 \le &\frac12\mathbb E\|\widehat f-f^*\|_\nu^2
+\frac32\|f_0-f^*\|_\nu^2
 +\frac{C_4\sS\sqrt{12\log(eN)}}{\neff}
\left\{\sD\log\bigl((\sW+1)\sC_A\bigr)
+\log(eN)\right\},
\label{eq:capped-risk-absorption}
\end{align}
where the first inequality takes expectations in
\eqref{eq:capped-basic-inequality} and applies
\eqref{eq:capped-combined-increment-bound} with $C_4=C_2+2C_3$;
and the second inequality uses \eqref{eq:capped-comparator-distance}.
Rearranging \eqref{eq:capped-risk-absorption} gives
\begin{align}
\mathbb E\|\widehat f-f^*\|_\nu^2
&\le3\|f_0-f^*\|_\nu^2
\nonumber\\
&\quad+\frac{C_5\sS\sqrt{12\log(eN)}}{\neff}
\left\{\sD\log\bigl((\sW+1)\sC_A\bigr)
+\log(eN)\right\},
\label{eq:capped-weighted-risk-bound}
\end{align}
where $C_5>0$ is a new constant depending only on $C_\epsilon$,
$\sC_M$, and $c_\epsilon$.
The bound \eqref{eq:capped-weighted-risk-bound} holds for every
deterministic $f_0\in\mathcal F$, with $C_5$ independent of $f_0$.

{\bf Step 3.}
The target risk satisfies
\begin{align}
&\mathbb E\|\widehat f-f^*\|_{\mathcal Q}^2
\nonumber \\
&\le\mathbb E\|\widehat f-f^*\|_\nu^2
+4\sC_M^2N^{-3}
\nonumber\\
&\le3\|f_0-f^*\|_\nu^2
 +\frac{C_5\sS\sqrt{12\log(eN)}}{\neff}
\left\{\sD\log\bigl((\sW+1)\sC_A\bigr)
+\log(eN)\right\}
 +4\sC_M^2N^{-3}
\nonumber\\
&\le3\|f_0-f^*\|_{\mathcal Q}^2
 +\frac{C_5\sS\sqrt{12\log(eN)}}{\neff}
\left\{\sD\log\bigl((\sW+1)\sC_A\bigr)
+\log(eN)\right\}
 +4\sC_M^2N^{-3},
\label{eq:capped-target-risk-bound}
\end{align}
where the first inequality applies \Cref{lem:multisource-target-risk-comparison}
to $\widehat f\in\mathcal F$;
the second inequality uses \eqref{eq:capped-weighted-risk-bound};
and the third inequality applies \Cref{lem:multisource-target-risk-comparison}.
Since \eqref{eq:capped-target-risk-bound} holds for every deterministic
$f_0\in\mathcal F$ and $C_5$ is independent of $f_0$, taking the infimum
over $f_0\in\mathcal F$ proves the lemma with $C_1=C_5$.
\end{proof}

\subsection{Proof of \Cref{thm:estimated-covariate-shift-rate}}

\begin{proof}[Proof of \Cref{thm:estimated-covariate-shift-rate}]
{\bf Step 1.}
Recall $N$ and $\neff$ from
\eqref{eq:multisource-total-sample-size} and
\eqref{eq:multisource-effective-sample-size}.
Let $\mathcal F$ be defined in \eqref{eq:localization-network-class},
and set
\begin{align}
\delta=\left(\frac{\sqrt{12\log(eN)}}{\neff}\right)^{p/(2p+K)}.
\label{eq:holder-approximation-tolerance}
\end{align}
By \eqref{eq:oracle-effective-sample-size} and
\eqref{eq:multisource-effective-sample-size}, $0<\delta\le1$.
Moreover,
\begin{align}
\log(2/\delta)
 \le C_5\log(eN),
\label{eq:holder-approximation-depth-scale}
\end{align}
where the inequality follows from
\eqref{eq:holder-approximation-tolerance},
\eqref{eq:multisource-effective-sample-size},
\eqref{eq:multisource-total-sample-size}, and
$\sigma_{\epsilon,j}^2\ge c_\epsilon$.
Also, \eqref{eq:holder-approximation-tolerance} gives
\begin{align}
\delta^{-K/p}
=\left(\frac{\neff}{\sqrt{12\log(eN)}}\right)^{K/(2p+K)}.
\label{eq:holder-approximation-width-scale}
\end{align}
The sparsity satisfies
\begin{align}
S\ge\frac{C_3C_2}{C_1+1}D\delta^{-K/p},
\label{eq:holder-approximation-sparsity-budget}
\end{align}
where the inequality follows from \eqref{eq:multisource-sparsity}, and
\eqref{eq:holder-approximation-width-scale}.
Choose $C_1\ge2$, $C_2\ge d$, and
$C_3\ge\max\{1,(C_1+1)/C_2\}$ sufficiently large.
Then $\sD\ge2$, $\sW\ge d$, and $\sS\ge\max\{\sW,\sD\}$, and
\eqref{eq:holder-approximation-depth-scale},
\eqref{eq:holder-approximation-width-scale}, and
\eqref{eq:holder-approximation-sparsity-budget} verify the budgets in
\Cref{lem:hierarchical-network-approximation}.
Choosing $\sC_A$ as required by that lemma yields $f_\delta\in\mathcal F$
satisfying
\begin{align}
\|f_\delta-f^*\|_\infty\le C_6\delta,
\label{eq:holder-hierarchical-uniform-approximation}
\end{align}
where $C_6>0$ depends only on $d,L,p,K,R,M$.
Consequently,
\begin{align}
\inf_{f\in\mathcal F}\|f-f^*\|_{\mathcal Q}^2
 \le\inf_{f\in\mathcal F}\|f-f^*\|_\infty^2
 \le\|f_\delta-f^*\|_\infty^2
 \le C_6^2\delta^2
 =C_6^2\left(\frac{\sqrt{12\log(eN)}}{\neff}\right)^{2p/(2p+K)},
\label{eq:holder-target-approximation}
\end{align}
where the first inequality uses $\mathcal Q([0,1]^d)=1$;
the second inequality uses $f_\delta\in\mathcal F$;
the third inequality uses \eqref{eq:holder-hierarchical-uniform-approximation};
and the equality uses \eqref{eq:holder-approximation-tolerance}.

{\bf Step 2.}
Applying \Cref{lem:capped-weighted-oracle-inequality} to
\eqref{eq:oracle-estimator} gives
\begin{align}
\mathbb E\|\widehat f-f^*\|_{\mathcal Q}^2
&\le3\inf_{f\in\mathcal F}\|f-f^*\|_{\mathcal Q}^2
\nonumber\\
&\quad+\frac{C_7\sS\sqrt{12\log(eN)}}{\neff}
\left\{\sD\log\bigl((\sW+1)\sC_A\bigr)
+\log(eN)\right\}
+4\sC_M^2N^{-3},
\label{eq:holder-estimated-oracle-application}
\end{align}
where $C_7>0$ depends only on $C_\epsilon,\sC_M,c_\epsilon$.
The complexity term satisfies
\begin{align}
&\frac{\sS\sqrt{12\log(eN)}}{\neff}
\left\{\sD\log\bigl((\sW+1)\sC_A\bigr)
+\log(eN)\right\}
\nonumber\\
&\qquad\le C_8
\left(\frac{\sqrt{12\log(eN)}}{\neff}\right)^{2p/(2p+K)}
\{\log(eN)\}^{3},
\label{eq:holder-complexity-bound}
\end{align}
where the inequality uses \eqref{eq:multisource-sparsity},
$1\le\neff/\sqrt{12\log(eN)}\le c_\epsilon^{-1}N$,
and a constant $C_8>0$ depending only on $\sC_A,C_1,C_2,C_3,p,K,c_\epsilon$.
Combining \eqref{eq:holder-estimated-oracle-application},
\eqref{eq:holder-target-approximation}, and
\eqref{eq:holder-complexity-bound} gives
\begin{align}
\mathbb E\|\widehat f-f^*\|_{\mathcal Q}^2
&\le C_9
\left(\frac{\sqrt{12\log(eN)}}{\neff}\right)^{2p/(2p+K)}
\{\log(eN)\}^{3}
+4\sC_M^2N^{-3},
\label{eq:holder-estimated-threshold-risk-bound}
\end{align}
where the inequality uses $\log(eN)\ge1$.

{\bf Step 3.}
The remainder satisfies
\begin{align}
N^{-3}
 \le c_\epsilon^{-2p/(2p+K)}
\left(\frac{\sqrt{12\log(eN)}}{\neff}\right)^{2p/(2p+K)},
\label{eq:holder-tail-approximation-comparison}
\end{align}
where the inequality uses
\eqref{eq:multisource-effective-sample-size},
\eqref{eq:multisource-total-sample-size},
$\sigma_{\epsilon,j}^2\ge c_\epsilon$, $N\ge2$, and
$\sqrt{12\log(eN)}\ge1$.
Combining \eqref{eq:holder-estimated-threshold-risk-bound} and
\eqref{eq:holder-tail-approximation-comparison} yields
\begin{align}
\mathbb E\|\widehat f-f^*\|_{\mathcal Q}^2
\le C_{10}\neff^{-2p/(2p+K)}
\{\log(eN)\}^{3+p/(2p+K)},
\nonumber
\end{align}
where the inequality uses $\log(eN)\ge1$.
This proves the theorem.
\end{proof}

\section{Technical results for \Cref{sec:multisource-adaptive}}

\begin{lemma}
\label{lem:adaptive-variance-weight-concentration}
Suppose
\Cref{ass:noise,ass:hierarchical-composition}, and \eqref{eq:multisource-adaptive-sample-size}
hold. Let $\widehat\nu_j$ and $\widehat\alpha_j$ be constructed as in
\Cref{sec:multisource-adaptive-algorithm}. For sufficiently large $N$,
with probability at least $1-N^{-3}$, the following bounds hold
simultaneously:
\begin{align}
\frac34\nu_j^*\le\widehat\nu_j\le\frac54\nu_j^*,
\qquad
\frac12\alpha_j^*\le\widehat\alpha_j\le2\alpha_j^*,
\qquad 
\frac{1}{\widehat\nu_j}\le\frac{4}{3c_\epsilon},
\qquad
\frac{\sigma_{\epsilon,j}^2}{\widehat\nu_j}\le\frac43,
\qquad j=1,\ldots,\rmJ.
\nonumber
\end{align}
\end{lemma}

\begin{proof}
{\bf Step 1.}
Fix $j\in\{1,\ldots,\rmJ\}$ and set
$Z_{ij}=Y_{ij}-\mathbb E Y_{ij}$ for $i\in\mathcal I_{j,0}$.
The response variance satisfies
\begin{align}
\nu_j^*
&=\operatorname{Var}(f^*(X_{ij}))+\sigma_{\epsilon,j}^2
\nonumber\\
&\le M^2+\sigma_{\epsilon,j}^2,
\nonumber
\end{align}
where the equality uses independence;
and the inequality uses $|f^*|\le M$ from \Cref{ass:hierarchical-composition}.
Moreover, \eqref{eq:multisource-adaptive-sample-size} and nonnegativity
of the signal variance give
\begin{align}
c_\epsilon\le\sigma_{\epsilon,j}^2\le\nu_j^*.
\label{eq:adaptive-response-variance-lower-bound}
\end{align}
Set $C_1=C_\epsilon+M^2/c_\epsilon$. For every $t\in\mathbb R$,
\begin{align}
\mathbb E\exp(tZ_{ij})
&=\mathbb E\exp\{t(f^*(X_{ij})-\mathbb E f^*(X_{ij}))\}
\mathbb E\exp(t\epsilon_{ij})
\nonumber\\
&\le\exp\left\{\frac{(M^2+C_\epsilon\sigma_{\epsilon,j}^2)t^2}{2}\right\}
 \le\exp\left(\frac{C_1\nu_j^*t^2}{2}\right),
\label{eq:adaptive-centered-response-mgf}
\end{align}
where the equality uses   independence and
zero mean in \Cref{ass:noise}; the first inequality uses Hoeffding's
lemma for $f^*(X_{ij})\in[-M,M]$ and   \Cref{ass:noise}; and the second inequality uses
\eqref{eq:adaptive-response-variance-lower-bound} and the definition of $C_1$.

{\bf Step 2.}
For each source $j$, define
\begin{align}
\overline Z_{j,0}=\frac{3}{n_j}\sum_{i\in\mathcal I_{j,0}}Z_{ij},
\qquad
T_j=\frac{3}{n_j}\sum_{i\in\mathcal I_{j,0}}Z_{ij}^2.
\nonumber
\end{align}
The definition of $\widehat V_j$ in \eqref{eq:adaptive-sample-variance-definition} gives
\begin{align}
\widehat V_j=\frac{n_j/3}{n_j/3-1}
\bigl(T_j-\overline Z_{j,0}^{\,2}\bigr).
\label{eq:adaptive-sample-variance-identity}
\end{align}
Chernoff's inequality applied to
\eqref{eq:adaptive-centered-response-mgf} gives, for every $u>0$,
\begin{align}
\mathbb P(|Z_{ij}|>u)
\le2\exp\left(-\frac{u^2}{2C_1\nu_j^*}\right).
\label{eq:adaptive-centered-response-tail}
\end{align}
Set $C_2=8C_1$. For every integer $q\ge2$,
\begin{align}
\mathbb E|Z_{ij}^2-\nu_j^*|^q
&\le2^q\mathbb E|Z_{ij}|^{2q}
\nonumber\\
&\le2^{q+1}(2C_1\nu_j^*)^q q!
\nonumber\\
&\le q!(C_2\nu_j^*)^q,
\label{eq:adaptive-centered-square-moments}
\end{align}
where the first inequality uses convexity and
$\nu_j^*=\mathbb E Z_{ij}^2$; the second inequality integrates
\eqref{eq:adaptive-centered-response-tail}; and the third inequality
uses $q\ge2$ and the definition of $C_2$.
Bernstein's inequality under the moment condition
\eqref{eq:adaptive-centered-square-moments} yields a constant
$C_3\ge48C_1$, depending only on $C_1$, such that
\begin{align}
\mathbb P\left(|T_j-\nu_j^*|>\frac{\nu_j^*}{8}\right)
\le2\exp(-n_j/C_3).
\label{eq:adaptive-second-moment-concentration}
\end{align}
Independence   and
\eqref{eq:adaptive-centered-response-mgf} also give
\begin{align}
\mathbb P\left(\overline Z_{j,0}^{\,2}>\frac{\nu_j^*}{8}\right)
&\le2\exp\left(-\frac{n_j}{48C_1}\right)
\nonumber\\
&\le2\exp(-n_j/C_3),
\label{eq:adaptive-mean-concentration}
\end{align}
where the first inequality applies the sub-Gaussian tail bound to the
average of $n_j/3$ independent centered responses; and the second
inequality uses $C_3\ge48C_1$.

{\bf Step 3.} 
Define
\begin{align}
\mathcal E=\bigcap_{j=1}^{\rmJ}
\left\{|T_j-\nu_j^*|\le\frac{\nu_j^*}{8},\quad
\overline Z_{j,0}^{\,2}\le\frac{\nu_j^*}{8}\right\}.
\nonumber
\end{align}
Choose the constant in
\eqref{eq:multisource-adaptive-sample-size} so that $C_0\ge5C_3$.
Then
\begin{align}
\mathbb P(\mathcal E^c)
&\le4\sum_{j=1}^{\rmJ}\exp(-n_j/C_3)
\nonumber\\
&\le4N(eN)^{-C_0/C_3}
\nonumber\\
&\le N^{-3},
\label{eq:adaptive-variance-failure-probability}
\end{align}
where the first inequality uses a union bound with
\eqref{eq:adaptive-second-moment-concentration} and
\eqref{eq:adaptive-mean-concentration}; the second inequality uses
$n_j\ge C_0\log(eN)$ and $\rmJ\le N$; and the third inequality
uses $C_0\ge5C_3$ and $N\ge2$.
For sufficiently large $N$, the sample-size assumption ensures
$n_j/3\ge10$ for every source. On $\mathcal E$,
\eqref{eq:adaptive-sample-variance-identity} therefore implies
\begin{align}
\frac34\nu_j^*\le \widehat V_j\le\frac54\nu_j^*,
\qquad j=1,\ldots,\rmJ.
\label{eq:adaptive-unclipped-variance-concentration}
\end{align}
For sufficiently large $N$, we also have $N^{-1}\le c_\epsilon$.
The definition $\widehat\nu_j=\max\{N^{-1},\widehat V_j\}$,
\eqref{eq:adaptive-response-variance-lower-bound}, and
\eqref{eq:adaptive-unclipped-variance-concentration} give, on $\mathcal E$,
\begin{align}
\frac34\nu_j^*\le\widehat\nu_j\le\frac54\nu_j^*,
\qquad j=1,\ldots,\rmJ.
\label{eq:adaptive-clipped-variance-concentration}
\end{align}

{\bf Step 4.}
On $\mathcal E$, the definitions of the mixture weights give
\begin{align}
\frac{\widehat\alpha_j}{\alpha_j^*}
=\frac{\nu_j^*/\widehat\nu_j}
{\sum_{k=1}^{\rmJ}\alpha_k^*(\nu_k^*/\widehat\nu_k)}.
\label{eq:adaptive-mixture-weight-ratio}
\end{align}
By \eqref{eq:adaptive-clipped-variance-concentration}, every ratio
$\nu_k^*/\widehat\nu_k$ lies in $[4/5,4/3]$.
Since $\alpha_k^*>0$ and $\sum_k\alpha_k^*=1$, the denominator in
\eqref{eq:adaptive-mixture-weight-ratio} also lies in $[4/5,4/3]$.
Therefore, \eqref{eq:adaptive-mixture-weight-ratio} gives
\begin{align}
\frac35\alpha_j^*\le\widehat\alpha_j\le\frac53\alpha_j^*,
\qquad j=1,\ldots,\rmJ.
\nonumber
\end{align}
Therefore $\frac12\alpha_j^*\le\widehat\alpha_j\le2\alpha_j^*$.
 Finally,
\eqref{eq:adaptive-response-variance-lower-bound} and
\eqref{eq:adaptive-clipped-variance-concentration} imply
\begin{align}
\frac{1}{\widehat\nu_j}\le\frac{4}{3c_\epsilon},
\qquad
\frac{\sigma_{\epsilon,j}^2}{\widehat\nu_j}\le\frac43,
\qquad j=1,\ldots,\rmJ.
\nonumber
\end{align}
All the asserted bounds hold on $\mathcal E$, whose probability is
at least $1-N^{-3}$ by
\eqref{eq:adaptive-variance-failure-probability}. This proves the lemma.
\end{proof}

\begin{lemma}
\label{lem:adaptive-pooled-ratio-error}
Suppose
\Cref{ass:noise,ass:holder-shift,ass:hierarchical-composition,ass:pooled-adaptive-source-ratios}, and \eqref{eq:multisource-adaptive-sample-size}
hold. Let $\widehat\alpha_j$ and $\widehat r$ be constructed as in
\Cref{sec:multisource-adaptive-algorithm}, and define
\begin{align}
\mathcal P_{\widehat\alpha}
=\sum_{j=1}^{\rmJ}\widehat\alpha_j\mathcal P_j,
\qquad
\rho_{\widehat\alpha}
=\frac{d\mathcal Q}{d\mathcal P_{\widehat\alpha}}.
\nonumber
\end{align}
Define the event
\begin{align}
\mathcal E=\bigcap_{j=1}^{\rmJ}
\left\{\frac34\nu_j^*\le\widehat\nu_j\le\frac54\nu_j^*\right\}.
\nonumber
\end{align}
Then
\begin{align}
\mathbb E\left[
\mathbf1_{\mathcal E}
\|\widehat r-\max\{1,\rho_{\widehat\alpha}\}\|_{\mathcal P_{\widehat\alpha}}^2
\right]
\le2\sum_{j=1}^{\rmJ}\frac{\varepsilon_j^2}{\alpha_j^*}.
\nonumber
\end{align}
\end{lemma}

\begin{proof}
{\bf Step 1.}
Fix the observations used to estimate the mixture weights and the
source-ratio estimators. The weights $\widehat\alpha_j$ are strictly
positive and sum to one, so $\mathcal P_j\ll\mathcal P_{\widehat\alpha}$
for every $j$. Choose nonnegative finite versions of the true ratios
$\rho_j$, and define
\begin{align}
\pi_j=\widehat\alpha_j
\frac{d\mathcal P_j}{d\mathcal P_{\widehat\alpha}},
\qquad j=1,\ldots,\rmJ.
\label{eq:adaptive-pointwise-mixture-weights}
\end{align}
The definition \eqref{eq:adaptive-pointwise-mixture-weights} and the
mixture identity imply, $\mathcal P_{\widehat\alpha}$-almost surely,
\begin{align}
\sum_{j=1}^{\rmJ}\pi_j=1.
\label{eq:adaptive-pointwise-weight-normalization}
\end{align}
Moreover, \eqref{eq:adaptive-pointwise-mixture-weights} gives
\begin{align}
\pi_j\,d\mathcal P_{\widehat\alpha}
=\widehat\alpha_j\,d\mathcal P_j.
\label{eq:adaptive-mixture-change-of-measure}
\end{align}
Multiplying \eqref{eq:adaptive-mixture-change-of-measure} by $\rho_j$
and using $\rho_j\,d\mathcal P_j=d\mathcal Q$ gives
\begin{align}
\pi_j\rho_j=\widehat\alpha_j\rho_{\widehat\alpha}
\qquad\mathcal P_{\widehat\alpha}\text{-almost surely}.
\label{eq:adaptive-mixture-ratio-product}
\end{align}
On $\{\rho_{\widehat\alpha}>0\}$,
\eqref{eq:adaptive-mixture-ratio-product} implies $\rho_j>0$ for every
$j$. On $\{\rho_{\widehat\alpha}=0\}$,
\eqref{eq:adaptive-pointwise-weight-normalization} and
\eqref{eq:adaptive-mixture-ratio-product} imply that at least one
$\rho_j$ is zero. Dividing \eqref{eq:adaptive-mixture-ratio-product}
by $\rho_j$ on $\{\rho_{\widehat\alpha}>0\}$ and summing using
\eqref{eq:adaptive-pointwise-weight-normalization} therefore yields
\begin{align}
\rho_{\widehat\alpha}
=\left(\sum_{j=1}^{\rmJ}
\frac{\widehat\alpha_j}{\rho_j}\right)^{-1}
\qquad\mathcal P_{\widehat\alpha}\text{-almost surely},
\label{eq:adaptive-true-pooled-ratio}
\end{align}
with the conventions $1/0=+\infty$ and $1/(+\infty)=0$.

{\bf Step 2.}
For $0<t\le\delta_N$, define the auxiliary clipped source estimates by
\begin{align}
\widetilde\rho_{j,t}=\max\{t,\widehat\rho_j\}
\label{eq:adaptive-auxiliary-source-ratios}
\end{align}
and their pooled estimate by
\begin{align}
\widehat\rho_t
=\left(\sum_{j=1}^{\rmJ}
\frac{\widehat\alpha_j}{\widetilde\rho_{j,t}}\right)^{-1}.
\label{eq:adaptive-auxiliary-pooled-ratio}
\end{align}
The positive summands in the denominator of
\eqref{eq:adaptive-auxiliary-pooled-ratio} give
\begin{align}
\frac{\widehat\rho_t}{\widetilde\rho_{j,t}}
\le\frac1{\widehat\alpha_j}.
\label{eq:adaptive-harmonic-coefficient-bound}
\end{align}
If $\widehat\rho_j(x)\ge\delta_N$ for every $j$, then
\eqref{eq:adaptive-auxiliary-source-ratios} and
\eqref{eq:adaptive-auxiliary-pooled-ratio} give
$\widehat\rho_t(x)=\widehat\rho(x)$ for every $0<t\le\delta_N$.
Otherwise, choose $j$ with $\widehat\rho_j(x)<\delta_N$. For every
$0<s\le\delta_N$, we have
\begin{align}
\widehat\rho_s(x)
&\le\frac{\widetilde\rho_{j,s}(x)}{\widehat\alpha_j}
\le\frac{\delta_N}{\widehat\alpha_j}
 \le1,
\nonumber
\end{align}
where the first inequality uses
\eqref{eq:adaptive-harmonic-coefficient-bound}; the second inequality
uses $s\le\delta_N$ and $\widehat\rho_j(x)<\delta_N$; and the third
inequality uses the choice $\delta_N\le\min_k\widehat\alpha_k$ in
\Cref{sec:multisource-adaptive-algorithm}.
Since $\widehat\rho_{\delta_N}=\widehat\rho$, both cases give
\begin{align}
\max\{1,\widehat\rho_t\}=\widehat r
\qquad\text{for every }0<t\le\delta_N.
\label{eq:adaptive-weight-floor-invariance}
\end{align}

{\bf Step 3.}
For every $0<t\le\delta_N$, the mixture identities give
\begin{align}
\sum_{j=1}^{\rmJ}\pi_j
\frac{\widehat\rho_t}{\widetilde\rho_{j,t}}
(\widetilde\rho_{j,t}-\rho_j)
&=\widehat\rho_t
-\widehat\rho_t\sum_{j=1}^{\rmJ}
\frac{\pi_j\rho_j}{\widetilde\rho_{j,t}}
 =\widehat\rho_t
-\widehat\rho_t\rho_{\widehat\alpha}
\sum_{j=1}^{\rmJ}
\frac{\widehat\alpha_j}{\widetilde\rho_{j,t}}
 =\widehat\rho_t-\rho_{\widehat\alpha},
\label{eq:adaptive-harmonic-error-identity}
\end{align}
where the first equality expands the sum and uses
\eqref{eq:adaptive-pointwise-weight-normalization}; the second equality
uses \eqref{eq:adaptive-mixture-ratio-product}; and the third equality
uses \eqref{eq:adaptive-auxiliary-pooled-ratio}.
Consequently,
\begin{align}
\|\widehat\rho_t-\rho_{\widehat\alpha}\|_{\mathcal P_{\widehat\alpha}}^2
&=\int\left\{\sum_{j=1}^{\rmJ}\pi_j
\frac{\widehat\rho_t}{\widetilde\rho_{j,t}}
(\widetilde\rho_{j,t}-\rho_j)\right\}^2
\,d\mathcal P_{\widehat\alpha}
\nonumber\\
&\le\sum_{j=1}^{\rmJ}\int\pi_j
\left(\frac{\widehat\rho_t}{\widetilde\rho_{j,t}}\right)^2
(\widetilde\rho_{j,t}-\rho_j)^2
\,d\mathcal P_{\widehat\alpha}
\nonumber\\
&\le\sum_{j=1}^{\rmJ}\frac1{\widehat\alpha_j}
\|\widetilde\rho_{j,t}-\rho_j\|_{\mathcal P_j}^2,
\label{eq:adaptive-auxiliary-pooled-error}
\end{align}
where the equality uses \eqref{eq:adaptive-harmonic-error-identity};
the first inequality uses Jensen's inequality and
\eqref{eq:adaptive-pointwise-weight-normalization}; and the second
inequality uses \eqref{eq:adaptive-harmonic-coefficient-bound} and
\eqref{eq:adaptive-mixture-change-of-measure}.

{\bf Step 4.}
For $a\in\mathbb R$, $b\ge0$, and $t>0$, lower clipping satisfies
\begin{align}
|\max\{t,a\}-b|^2
 \le|a-b|^2+(t-b)_+^2
 \le|a-b|^2+t^2,
\label{eq:adaptive-auxiliary-clipping-error}
\end{align}
where the first inequality follows by considering $a\ge t$, $a<t\le b$,
and $a<t$ with $b<t$; and the second inequality uses $b\ge0$.
Therefore, for every $0<t\le\delta_N$,
\begin{align}
\|\widehat r-\max\{1,\rho_{\widehat\alpha}\}\|_{\mathcal P_{\widehat\alpha}}^2
&\le\|\widehat\rho_t-\rho_{\widehat\alpha}\|_{\mathcal P_{\widehat\alpha}}^2
 \le\sum_{j=1}^{\rmJ}\frac1{\widehat\alpha_j}
\left\{\|\widehat\rho_j-\rho_j\|_{\mathcal P_j}^2+t^2\right\},
\label{eq:adaptive-weight-error-with-floor}
\end{align}
where the first inequality uses
\eqref{eq:adaptive-weight-floor-invariance} and the
$1$-Lipschitz property of $u\mapsto\max\{1,u\}$; and the second
inequality uses \eqref{eq:adaptive-auxiliary-pooled-error} and
\eqref{eq:adaptive-auxiliary-clipping-error}.
For the fixed realization, $\delta_N>0$ and
$\sum_j1/\widehat\alpha_j<\infty$. Letting $t\downarrow0$ in
\eqref{eq:adaptive-weight-error-with-floor} gives
\begin{align}
\|\widehat r-\max\{1,\rho_{\widehat\alpha}\}\|_{\mathcal P_{\widehat\alpha}}^2\le\sum_{j=1}^{\rmJ}\frac1{\widehat\alpha_j}
\|\widehat\rho_j-\rho_j\|_{\mathcal P_j}^2.
\label{eq:adaptive-pooled-weight-pathwise-error}
\end{align}

{\bf Step 5.}
On $\mathcal E$, we have
\begin{align}
\frac1{\widehat\alpha_j}
&=\frac1{\alpha_j^*}
\frac{\sum_{k=1}^{\rmJ}\alpha_k^*(\nu_k^*/\widehat\nu_k)}
{\nu_j^*/\widehat\nu_j}
 \le\frac5{3\alpha_j^*}
 \le\frac2{\alpha_j^*},
\label{eq:adaptive-inverse-mixture-weight-bound}
\end{align}
where the equality uses \eqref{eq:adaptive-mixture-weight-ratio};
the first inequality uses $4/5\le\nu_k^*/\widehat\nu_k\le4/3$
for every $k$ on $\mathcal E$ and $\sum_k\alpha_k^*=1$; and the
second inequality uses $5/3\le2$.
Taking expectations gives
\begin{align}
\mathbb E\left[
\mathbf1_{\mathcal E}
\|\widehat r-\max\{1,\rho_{\widehat\alpha}\}\|_{\mathcal P_{\widehat\alpha}}^2
\right]
&\le\sum_{j=1}^{\rmJ}
\mathbb E\left[
\frac{\mathbf1_{\mathcal E}}{\widehat\alpha_j}
\|\widehat\rho_j-\rho_j\|_{\mathcal P_j}^2
\right]
 \le2\sum_{j=1}^{\rmJ}\frac1{\alpha_j^*}
\mathbb E\|\widehat\rho_j-\rho_j\|_{\mathcal P_j}^2
 \le2\sum_{j=1}^{\rmJ}\frac{\varepsilon_j^2}{\alpha_j^*},
\nonumber
\end{align}
where the first inequality uses
\eqref{eq:adaptive-pooled-weight-pathwise-error}; the second inequality
uses \eqref{eq:adaptive-inverse-mixture-weight-bound} and
$\mathbf1_{\mathcal E}\le1$; and the third inequality uses
\eqref{eq:pooled-adaptive-source-errors}. This proves the lemma.
\end{proof}

\begin{lemma}
\label{lem:adaptive-pooled-tail-grid}
Suppose
\Cref{ass:noise,ass:holder-shift,ass:hierarchical-composition,ass:pooled-adaptive-source-ratios}, and \eqref{eq:multisource-adaptive-sample-size} 
hold. Let $\mathcal E$, $\mathcal P_{\widehat\alpha}$, and
$\rho_{\widehat\alpha}$ be as in \Cref{lem:adaptive-pooled-ratio-error},
and let $N_{\widehat\nu}$ and $\Lambda$ be as in
\Cref{sec:multisource-adaptive-algorithm}.
There exist constants $C_1,C_2,C_3>0$, depending only on
$M,c_\epsilon$, such that the following holds for sufficiently large
$N$. If $\neff>C_1\sqrt{\log(eN)}$, then on $\mathcal E$
we have $N_{\widehat\nu}\ge1$ and there exists $\lambda_0\in\Lambda$
satisfying
\begin{align}
\frac{\lambda_0}{N_{\widehat\nu}}
\le C_2\frac{\sqrt{\log(eN)}}{\neff}
\nonumber
\end{align}
and
\begin{align}
\int(\rho_{\widehat\alpha}-\lambda_0)_+
\,d\mathcal P_{\widehat\alpha}
+\int(\rho_{\widehat\alpha}-\lambda_0)_+^2
\,d\mathcal P_{\widehat\alpha}
\le C_3N^{-4}.
\nonumber
\end{align}
\end{lemma}

\begin{proof}
{\bf Step 1.}
Work on $\mathcal E$ and fix the estimated mixture weights.
Set $C_4=1+M^2/c_\epsilon$. For each $j$,
\begin{align}
\nu_j^*
 =\operatorname{Var}(f^*(X_{ij}))+\sigma_{\epsilon,j}^2
 \le M^2+\sigma_{\epsilon,j}^2
 \le C_4\sigma_{\epsilon,j}^2,
\label{eq:adaptive-response-noise-comparison}
\end{align}
where the equality uses independence;
the first inequality uses $|f^*|\le M$; and the second inequality
uses $\sigma_{\epsilon,j}^2\ge c_\epsilon$ and the definition of $C_4$.
Define
\begin{align}
H_{\widehat\nu}
=\sum_{j=1}^{\rmJ}
\frac{n_j}{3(1+\sigma_{\rho,j})\widehat\nu_j}.
\label{eq:adaptive-weighted-effective-sample-size}
\end{align}
On $\mathcal E$, we have
\begin{align}
H_{\widehat\nu}
&\ge\frac4{15}\sum_{j=1}^{\rmJ}
\frac{n_j}{(1+\sigma_{\rho,j})\nu_j^*}
 \ge\frac4{15C_4}\neff,
\label{eq:adaptive-weighted-effective-size-lower-bound}
\end{align}
where the first inequality uses
\eqref{eq:adaptive-weighted-effective-sample-size} and
$\widehat\nu_j\le5\nu_j^*/4$ on $\mathcal E$; and the second inequality
uses \eqref{eq:adaptive-response-noise-comparison} and
\eqref{eq:multisource-effective-sample-size}.
Similarly, on $\mathcal E$,
\begin{align}
H_{\widehat\nu}
&\le\frac49\sum_{j=1}^{\rmJ}
\frac{n_j}{(1+\sigma_{\rho,j})\nu_j^*}
 \le\frac49\neff,
\nonumber
\end{align}
where the first inequality uses
\eqref{eq:adaptive-weighted-effective-sample-size} and
$\widehat\nu_j\ge3\nu_j^*/4$ on $\mathcal E$; and the second inequality
uses \eqref{eq:adaptive-response-variance-lower-bound} and
\eqref{eq:multisource-effective-sample-size}.

{\bf Step 2.}
Fix $t\ge1$. For every $j$ with $\sigma_{\rho,j}>0$,
\begin{align}
\mathcal P_j\bigl(\rho_j>(1+\sigma_{\rho,j})t\bigr)
&\le\exp\left\{-\frac{((1+\sigma_{\rho,j})t-1)^2}
{2\sigma_{\rho,j}^2}\right\}
 \le\exp(-t^2/2),
\label{eq:adaptive-source-ratio-tail}
\end{align}
where the first inequality applies Chernoff's bound using
\Cref{ass:holder-shift}; and the second inequality uses
$(1+\sigma_{\rho,j})t-1\ge\sigma_{\rho,j}t$ for $t\ge1$.
If $\sigma_{\rho,j}=0$, then $\rho_j=1$ $\mathcal P_j$-almost surely,
so \eqref{eq:adaptive-source-ratio-tail} also holds in this case.
In addition,
\begin{align}
\mathcal Q\bigl(\rho_j>(1+\sigma_{\rho,j})t\bigr)
&=\int\rho_j\mathbf1_{\{\rho_j>(1+\sigma_{\rho,j})t\}}
\,d\mathcal P_j
 \le\left(\int\rho_j^2\,d\mathcal P_j\right)^{1/2}
\mathcal P_j\bigl(\rho_j>(1+\sigma_{\rho,j})t\bigr)^{1/2}
\nonumber\\
&\le\sqrt{1+\sigma_{\rho,j}^2}\exp(-t^2/4),
\label{eq:adaptive-target-source-ratio-tail}
\end{align}
where the equality uses $d\mathcal Q=\rho_j\,d\mathcal P_j$;
the first inequality uses Cauchy--Schwarz; and the second inequality
uses \eqref{eq:adaptive-source-ratio-tail} and
$\int\rho_j^2\,d\mathcal P_j\le1+\sigma_{\rho,j}^2$ from
\Cref{ass:holder-shift}.
The definitions of $\widehat\alpha_j$, $N_{\widehat\nu}$, and
\eqref{eq:adaptive-weighted-effective-sample-size} give
\begin{align}
\sum_{j=1}^{\rmJ}\frac{\widehat\alpha_j}{1+\sigma_{\rho,j}}
=\frac{H_{\widehat\nu}}{N_{\widehat\nu}}.
\label{eq:adaptive-harmonic-effective-size-identity}
\end{align}
For $\mathcal P_{\widehat\alpha}$-almost every point where
$\rho_j\le(1+\sigma_{\rho,j})t$ for all $j$,
\begin{align}
\rho_{\widehat\alpha}
&=\left(\sum_{j=1}^{\rmJ}\frac{\widehat\alpha_j}{\rho_j}\right)^{-1}
 \le t\left(\sum_{j=1}^{\rmJ}
\frac{\widehat\alpha_j}{1+\sigma_{\rho,j}}\right)^{-1}
 =\frac{N_{\widehat\nu}}{H_{\widehat\nu}}t,
\label{eq:adaptive-pooled-ratio-envelope}
\end{align}
where the first equality uses \eqref{eq:adaptive-true-pooled-ratio}
with its conventions for zero ratios; the inequality uses
$\rho_j\le(1+\sigma_{\rho,j})t$ for every $j$;
and the second equality uses
\eqref{eq:adaptive-harmonic-effective-size-identity}.
Consequently,
\begin{align}
\mathcal Q\left(\rho_{\widehat\alpha}>
\frac{N_{\widehat\nu}}{H_{\widehat\nu}}t\right)
 \le\sum_{j=1}^{\rmJ}
\mathcal Q\bigl(\rho_j>(1+\sigma_{\rho,j})t\bigr)
 \le\exp(-t^2/4)\sum_{j=1}^{\rmJ}\sqrt{1+\sigma_{\rho,j}^2}
 \le2N\exp(-t^2/4),
\label{eq:adaptive-pooled-target-tail}
\end{align}
where the first inequality uses \eqref{eq:adaptive-pooled-ratio-envelope}
and a union bound; the second inequality uses
\eqref{eq:adaptive-target-source-ratio-tail}; and the third inequality
uses $\sqrt{1+\sigma_{\rho,j}^2}\le1+\sigma_{\rho,j}<1+n_j$,
$\sum_jn_j=N$, and $\rmJ\le N$.

{\bf Step 3.}
Define
\begin{align}
t_N=\sqrt{24\log(eN)},
\qquad
\lambda_*=\frac{N_{\widehat\nu}}{H_{\widehat\nu}}t_N.
\label{eq:adaptive-tail-threshold}
\end{align}
Consequently \begin{align}
\frac{N_{\widehat\nu}}{H_{\widehat\nu}}
 =\left(\sum_{j=1}^{\rmJ}
\frac{\widehat\alpha_j}{1+\sigma_{\rho,j}}\right)^{-1}
 \le\max_{1\le j\le\rmJ}(1+\sigma_{\rho,j})
 \le2N,
\label{eq:adaptive-effective-size-ratio-upper-bound}
\end{align}
where the equality uses
\eqref{eq:adaptive-harmonic-effective-size-identity}; the first
inequality uses $\widehat\alpha_j\ge0$ and $\sum_j\widehat\alpha_j=1$;
and the second inequality uses $\sigma_{\rho,j}<n_j\le N$ and $N\ge2$.
The first truncation moment satisfies
\begin{align}
\int(\rho_{\widehat\alpha}-\lambda_*)_+
\,d\mathcal P_{\widehat\alpha}
&=\int_{\{\rho_{\widehat\alpha}>\lambda_*\}}
\frac{\rho_{\widehat\alpha}-\lambda_*}
{\rho_{\widehat\alpha}}\,d\mathcal Q
 \le\mathcal Q(\rho_{\widehat\alpha}>\lambda_*)
 \le2N\exp(-t_N^2/4),
\label{eq:adaptive-first-truncation-moment}
\end{align}
where the equality uses
$d\mathcal Q=\rho_{\widehat\alpha}\,d\mathcal P_{\widehat\alpha}$
on $\{\rho_{\widehat\alpha}>\lambda_*\}$;    the first inequality uses
$0\le(\rho_{\widehat\alpha}-\lambda_*)/\rho_{\widehat\alpha}\le1$
on the integration set; and the second
inequality uses \eqref{eq:adaptive-pooled-target-tail} and
\eqref{eq:adaptive-tail-threshold}.
Moreover,
\begin{align}
\int(\rho_{\widehat\alpha}-\lambda_*)_+\,d\mathcal Q
&=\int_{\lambda_*}^{\infty}
\mathcal Q(\rho_{\widehat\alpha}>v)\,dv
 =\frac{N_{\widehat\nu}}{H_{\widehat\nu}}
\int_{t_N}^{\infty}
\mathcal Q\left(\rho_{\widehat\alpha}>
\frac{N_{\widehat\nu}}{H_{\widehat\nu}}u\right)\,du
 \le2N\frac{N_{\widehat\nu}}{H_{\widehat\nu}}
\int_{t_N}^{\infty}\exp(-u^2/4)\,du
\nonumber\\
&\le\frac{2N}{t_N}\frac{N_{\widehat\nu}}{H_{\widehat\nu}}
\int_{t_N}^{\infty}u\exp(-u^2/4)\,du
 =\frac{4N}{t_N}\frac{N_{\widehat\nu}}{H_{\widehat\nu}}
\exp(-t_N^2/4),
\label{eq:adaptive-pooled-excess-moment}
\end{align}
where the first equality uses the tail-integral formula; the second
equality substitutes $v=(N_{\widehat\nu}/H_{\widehat\nu})u$ and uses
\eqref{eq:adaptive-tail-threshold}; the first inequality uses
\eqref{eq:adaptive-pooled-target-tail} and $t_N\ge1$; the second
inequality uses $u/t_N\ge1$ for $u\ge t_N$; and the third equality
evaluates the integral.
The second truncation moment satisfies
\begin{align}
\int(\rho_{\widehat\alpha}-\lambda_*)_+^2
\,d\mathcal P_{\widehat\alpha}
&=\int_{\{\rho_{\widehat\alpha}>\lambda_*\}}
\frac{(\rho_{\widehat\alpha}-\lambda_*)^2}
{\rho_{\widehat\alpha}}\,d\mathcal Q
 \nonumber \\
 &\le\int(\rho_{\widehat\alpha}-\lambda_*)_+\,d\mathcal Q
 \le\frac{4N}{t_N}\frac{N_{\widehat\nu}}{H_{\widehat\nu}}
\exp(-t_N^2/4)
 \le8N^2\exp(-t_N^2/4),
\label{eq:adaptive-second-truncation-moment}
\end{align}
where the equality uses
$d\mathcal Q=\rho_{\widehat\alpha}\,d\mathcal P_{\widehat\alpha}$
on $\{\rho_{\widehat\alpha}>\lambda_*\}$; the first inequality uses
$\rho_{\widehat\alpha}-\lambda_*\le\rho_{\widehat\alpha}$ on the
integration set;
the second inequality uses \eqref{eq:adaptive-pooled-excess-moment};
and the third inequality uses
\eqref{eq:adaptive-effective-size-ratio-upper-bound} and $t_N\ge1$.
For an absolute constant $C_5>0$, we obtain
\begin{align}
&\int(\rho_{\widehat\alpha}-\lambda_*)_+
\,d\mathcal P_{\widehat\alpha}
+\int(\rho_{\widehat\alpha}-\lambda_*)_+^2
\,d\mathcal P_{\widehat\alpha}
 \le(2N+8N^2)\exp(-t_N^2/4)
 \le C_5N^{-4},
\label{eq:adaptive-oracle-truncation-moments}
\end{align}
where the first inequality uses
\eqref{eq:adaptive-first-truncation-moment} and
\eqref{eq:adaptive-second-truncation-moment}; and the second inequality
uses $\exp(-t_N^2/4)=(eN)^{-6}$ from
\eqref{eq:adaptive-tail-threshold} and $N\ge2$.

{\bf Step 4.}
Set $C_6=15C_4\sqrt{24}/2$ and choose $C_1=\max\{1,C_6\}$.
The condition $\neff>C_1\sqrt{\log(eN)}$,
\eqref{eq:adaptive-weighted-effective-size-lower-bound}, and
\eqref{eq:adaptive-tail-threshold} imply
\begin{align}
H_{\widehat\nu}\ge2t_N\ge2.
\label{eq:adaptive-effective-size-grid-condition}
\end{align}
The definition \eqref{eq:adaptive-weighted-effective-sample-size}
and $1+\sigma_{\rho,j}\ge1$ give
$H_{\widehat\nu}\le N_{\widehat\nu}$. Hence
$N_{\widehat\nu}\ge1$ by
\eqref{eq:adaptive-effective-size-grid-condition}.
Moreover,
\begin{align}
1\le t_N\le\lambda_*
=\frac{N_{\widehat\nu}}{H_{\widehat\nu}}t_N
\le\frac{N_{\widehat\nu}}2,
\label{eq:adaptive-threshold-grid-range}
\end{align}
where the first inequality uses $N\ge2$ and
\eqref{eq:adaptive-tail-threshold}; the second inequality uses
$H_{\widehat\nu}\le N_{\widehat\nu}$; the equality uses
\eqref{eq:adaptive-tail-threshold}; and the third inequality uses
\eqref{eq:adaptive-effective-size-grid-condition}.
Choose $\lambda_0=2^{\lceil\log_2\lambda_*\rceil}$. Then
\begin{align}
\lambda_*\le\lambda_0\le2\lambda_*\le N_{\widehat\nu},
\label{eq:adaptive-rounded-grid-threshold}
\end{align}
where the first inequality uses
$\lceil\log_2\lambda_*\rceil\ge\log_2\lambda_*$; the second
inequality uses $\lceil\log_2\lambda_*\rceil<\log_2\lambda_*+1$;
and the third inequality uses \eqref{eq:adaptive-threshold-grid-range}.
By \eqref{eq:adaptive-threshold-grid-range} and
\eqref{eq:adaptive-rounded-grid-threshold}, $\lambda_0$ is a power
of two between $1$ and $N_{\widehat\nu}$, so $\lambda_0\in\Lambda$.
Moreover,
\begin{align}
\frac{\lambda_0}{N_{\widehat\nu}}
&\le\frac{2\lambda_*}{N_{\widehat\nu}}
 =\frac{2t_N}{H_{\widehat\nu}}
 \le C_6\frac{\sqrt{\log(eN)}}{\neff},
\label{eq:adaptive-oracle-normalized-threshold}
\end{align}
where the first inequality uses \eqref{eq:adaptive-rounded-grid-threshold};
the equality uses \eqref{eq:adaptive-tail-threshold}; and the second
inequality uses \eqref{eq:adaptive-weighted-effective-size-lower-bound}
and $t_N=\sqrt{24\log(eN)}$.
Finally,
\begin{align}
&\int(\rho_{\widehat\alpha}-\lambda_0)_+
\,d\mathcal P_{\widehat\alpha}
+\int(\rho_{\widehat\alpha}-\lambda_0)_+^2
\,d\mathcal P_{\widehat\alpha}
 \le\int(\rho_{\widehat\alpha}-\lambda_*)_+
\,d\mathcal P_{\widehat\alpha}
+\int(\rho_{\widehat\alpha}-\lambda_*)_+^2
\,d\mathcal P_{\widehat\alpha}
 \le C_5N^{-4},
\nonumber
\end{align}
where the first inequality uses $\lambda_0\ge\lambda_*$ from
\eqref{eq:adaptive-rounded-grid-threshold} and monotonicity of the
positive part; and the second inequality uses
\eqref{eq:adaptive-oracle-truncation-moments}.
For later use, if additionally
$\neff>C_1\{\log(eN)\}^{7/2}$, then
\begin{align}
\frac{\lambda_0\{\log(eN)\}^3}{N_{\widehat\nu}}
&\le C_6\frac{\{\log(eN)\}^{7/2}}{\neff}
 \le\frac{C_6}{C_1}
 \le1,
\label{eq:adaptive-oracle-approximation-scale}
\end{align}
where the first inequality uses
\eqref{eq:adaptive-oracle-normalized-threshold}; the second inequality
uses $\neff>C_1\{\log(eN)\}^{7/2}$; and the third inequality
uses $C_1\ge C_6$.
Taking $C_2=C_6$ and $C_3=C_5$ proves the lemma.
\end{proof}

\begin{lemma}
\label{lem:adaptive-training-error}
Suppose the assumptions in \Cref{lem:adaptive-pooled-tail-grid} hold
and $\neff>C_0\{\log(eN)\}^{7/2}$, where $C_0>0$ is sufficiently
large and depends only on $M,c_\epsilon$. Let $\lambda_0$ be chosen
as in \Cref{lem:adaptive-pooled-tail-grid}, and let
$\widehat f_{\lambda_0}$ be the corresponding training estimator in
\Cref{sec:multisource-adaptive-algorithm}. Define
\begin{align}
\mathcal R(f)
=\int\widehat r(x)(f(x)-f^*(x))^2
\,d\mathcal P_{\widehat\alpha}(x).
\nonumber
\end{align}
Write $\mathbb E_{\mathrm{train}}$ for expectation over the training
observations $\{(X_{ij},Y_{ij}):i\in\mathcal I_{j,1},\ j=1,\ldots,\rmJ\}$,
conditional on the data in $\mathcal I_{j,0}$ and on $\widehat\rho_j$,
for $j=1,\ldots,\rmJ$. Choose $\sC_A,C_\sD,C_\sW,C_\sS$ sufficiently large in
\eqref{eq:adaptive-network-class}, \eqref{eq:adaptive-network-depth},
\eqref{eq:adaptive-network-width}, and \eqref{eq:adaptive-network-sparsity}.
For sufficiently large $N$, on $\mathcal E$, we have
\begin{align}
\mathbb E_{\mathrm{train}}\mathcal R(\widehat f_{\lambda_0})
&\le C_1\left(\frac{\lambda_0\{\log(eN)\}^3}{N_{\widehat\nu}}\right)^{2p/(2p+K)}
 +C_2\|\widehat r-\max\{1,\rho_{\widehat\alpha}\}\|_{\mathcal P_{\widehat\alpha}}^2+C_3N^{-4}.
\nonumber
\end{align}
The constants $C_1,C_2,C_3>0$ depend only on
$d,L,p,K,R,M,c_\epsilon,C_\epsilon$.
\end{lemma}

\begin{proof}
{\bf Step 1.}
Work on $\mathcal E$. Define the clipped weight and population risk by
\begin{align}
w_0(x)=\min\{\widehat r(x),\lambda_0\},
\qquad
\mathcal R_{\lambda_0}(f)
=\int w_0(x)(f(x)-f^*(x))^2
\,d\mathcal P_{\widehat\alpha}(x).
\label{eq:adaptive-clipped-training-risk}
\end{align}
For $f\in\mathcal F_{\lambda_0}$, define the quadratic fluctuation by
\begin{align}
\phi(f)
=\frac1{N_{\widehat\nu}}
\sum_{j=1}^{\rmJ}\sum_{i\in\mathcal I_{j,1}}
\frac{w_0(X_{ij})}{\widehat\nu_j}
(f(X_{ij})-f^*(X_{ij}))^2
-\mathcal R_{\lambda_0}(f)
\label{eq:adaptive-training-quadratic-process}
\end{align}
and the noise term by
\begin{align}
\psi(f)
=\frac1{N_{\widehat\nu}}
\sum_{j=1}^{\rmJ}\sum_{i\in\mathcal I_{j,1}}
\frac{w_0(X_{ij})}{\widehat\nu_j}
\epsilon_{ij}(f(X_{ij})-f^*(X_{ij})).
\label{eq:adaptive-training-noise-process}
\end{align}
Then it holds that \begin{align}
&\frac1{N_{\widehat\nu}}
\sum_{j=1}^{\rmJ}\sum_{i\in\mathcal I_{j,1}}
\frac{w_0(X_{ij})}{\widehat\nu_j}
\bigl\{(Y_{ij}-f(X_{ij}))^2-(Y_{ij}-f^*(X_{ij}))^2\bigr\}
\nonumber\\
&\quad=\frac1{N_{\widehat\nu}}
\sum_{j=1}^{\rmJ}\sum_{i\in\mathcal I_{j,1}}
\frac{w_0(X_{ij})}{\widehat\nu_j}
\bigl\{(f(X_{ij})-f^*(X_{ij}))^2
-2\epsilon_{ij}(f(X_{ij})-f^*(X_{ij}))\bigr\}
\nonumber\\
&\quad=\mathcal R_{\lambda_0}(f)+\phi(f)-2\psi(f),
\label{eq:adaptive-training-excess-loss-decomposition}
\end{align}
where the first equality expands the squares using
\eqref{eq:multisource-model}; and the second equality uses
\eqref{eq:adaptive-training-quadratic-process} and
\eqref{eq:adaptive-training-noise-process}.
For every $f_0\in\mathcal F_{\lambda_0}$, the minimizing property
\eqref{eq:adaptive-training-estimator}, the definition of $w_0$ in
\eqref{eq:adaptive-clipped-training-risk}, and
\eqref{eq:adaptive-training-excess-loss-decomposition} give
\begin{align}
\mathcal R_{\lambda_0}(\widehat f_{\lambda_0})
&\le\mathcal R_{\lambda_0}(f_0)
+\phi(f_0)-\phi(\widehat f_{\lambda_0})
+2\{\psi(\widehat f_{\lambda_0})-\psi(f_0)\}.
\label{eq:adaptive-training-basic-inequality}
\end{align}

{\bf Step 2.}
Write $\mathbb P_{\mathrm{train}}$ for probability over the training
observations. Fix $f\in\mathcal F_{\lambda_0}$ independently of these
observations, and define
\begin{align}
U_{ij}(f)
=\frac{w_0(X_{ij})}{N_{\widehat\nu}\widehat\nu_j}
(f(X_{ij})-f^*(X_{ij}))^2,
\qquad i\in\mathcal I_{j,1}.
\nonumber
\end{align}
Since $|f-f^*|\le2M$, $w_0\le\lambda_0$, and
$1/\widehat\nu_j\le4/(3c_\epsilon)$ on $\mathcal E$, we have
\begin{align}
0\le U_{ij}(f)\le C_4\frac{\lambda_0}{N_{\widehat\nu}},
\qquad C_4=\frac{16M^2}{3c_\epsilon}.
\label{eq:adaptive-quadratic-summand-bound}
\end{align}
The definition of $\widehat\alpha_j$ gives
\begin{align}
\sum_{j=1}^{\rmJ}\sum_{i\in\mathcal I_{j,1}}
\mathbb E_{\mathrm{train}}U_{ij}(f)
&=\sum_{j=1}^{\rmJ}
\frac{n_j}{3N_{\widehat\nu}\widehat\nu_j}
\int w_0(x)(f(x)-f^*(x))^2\,d\mathcal P_j(x)
 =\mathcal R_{\lambda_0}(f),
\label{eq:adaptive-quadratic-summand-mean}
\end{align}
where the first equality uses $|\mathcal I_{j,1}|=n_j/3$ and
$X_{ij}\sim\mathcal P_j$; and the second equality uses
$n_j/(3N_{\widehat\nu}\widehat\nu_j)=\widehat\alpha_j$ and
\eqref{eq:adaptive-clipped-training-risk}.
Moreover,
\begin{align}
\sum_{j=1}^{\rmJ}\sum_{i\in\mathcal I_{j,1}}
\operatorname{Var}_{\mathrm{train}}(U_{ij}(f))
 \le\sum_{j=1}^{\rmJ}\sum_{i\in\mathcal I_{j,1}}
\mathbb E_{\mathrm{train}}U_{ij}(f)^2
 \le C_4\frac{\lambda_0}{N_{\widehat\nu}}
\sum_{j=1}^{\rmJ}\sum_{i\in\mathcal I_{j,1}}
\mathbb E_{\mathrm{train}}U_{ij}(f)
 =C_4\frac{\lambda_0}{N_{\widehat\nu}}
\mathcal R_{\lambda_0}(f),
\label{eq:adaptive-quadratic-variance-bound}
\end{align}
where  the second inequality uses \eqref{eq:adaptive-quadratic-summand-bound};
and the equality uses \eqref{eq:adaptive-quadratic-summand-mean}.
Apply \Cref{lem:bernstein-moment-bound} to the independent centered
variables $U_{ij}(f)-\mathbb E_{\mathrm{train}}U_{ij}(f)$, using
\eqref{eq:adaptive-quadratic-summand-bound} and
\eqref{eq:adaptive-quadratic-variance-bound}. For every $t>0$, this gives
\begin{align}
\mathbb P_{\mathrm{train}}\left(
|\phi(f)|>
\sqrt{2C_4\frac{\lambda_0t}{N_{\widehat\nu}}
\mathcal R_{\lambda_0}(f)}
+\frac{C_4\lambda_0t}{3N_{\widehat\nu}}
\right)\le2e^{-t}.
\label{eq:adaptive-quadratic-bernstein-tail}
\end{align}
By $\sqrt{2uv}\le u/8+4v$ for $u,v\ge0$,
\eqref{eq:adaptive-quadratic-bernstein-tail} implies that, for a constant
$C_5>0$ depending only on $M,c_\epsilon$,
\begin{align}
\mathbb P_{\mathrm{train}}\left(
|\phi(f)|>\frac18\mathcal R_{\lambda_0}(f)
+C_5\frac{\lambda_0t}{N_{\widehat\nu}}
\right)\le2e^{-t}.
\label{eq:adaptive-quadratic-localized-tail}
\end{align}
Let $\mathcal G\subseteq\mathcal F_{\lambda_0}$ be a nonempty finite
set chosen independently of the training observations. Applying
\eqref{eq:adaptive-quadratic-localized-tail} to each $g\in\mathcal G$
with $t$ replaced by $\log(2|\mathcal G|)+t$, and taking a union bound,
shows that, with $\mathbb P_{\mathrm{train}}$-probability at least
$1-e^{-t}$,
\begin{align}
|\phi(g)|\le\frac18\mathcal R_{\lambda_0}(g)
+C_5\frac{\lambda_0}{N_{\widehat\nu}}
\{\log(2|\mathcal G|)+t\}
\qquad\text{for every }g\in\mathcal G.
\nonumber
\end{align}

{\bf Step 3.}
Fix $f\in\mathcal F_{\lambda_0}$ independently of the training
observations. Write $\mathbb E_\epsilon$ and $\mathbb P_\epsilon$
for expectation and probability over the training errors
$\{\epsilon_{ij}:i\in\mathcal I_{j,1},\ j=1,\ldots,\rmJ\}$,
with the training covariates fixed.
For every $s\in\mathbb R$, we have
\begin{align}
\mathbb E_\epsilon\exp\{s\psi(f)\}
&\le\exp\left\{
\frac{C_\epsilon s^2}{2N_{\widehat\nu}^2}
\sum_{j=1}^{\rmJ}\sum_{i\in\mathcal I_{j,1}}
\frac{w_0(X_{ij})^2\sigma_{\epsilon,j}^2}{\widehat\nu_j^2}
(f(X_{ij})-f^*(X_{ij}))^2\right\}
\nonumber\\
&\le\exp\left\{
\frac{C_6s^2\lambda_0}{2N_{\widehat\nu}}
\bigl(\mathcal R_{\lambda_0}(f)+\phi(f)\bigr)\right\},
\label{eq:adaptive-noise-conditional-mgf}
\end{align}
where the first inequality uses the independence and exponential-moment
bound in \Cref{ass:noise} together with
\eqref{eq:adaptive-training-noise-process}; and the second inequality
uses $w_0^2\le\lambda_0w_0$,
$\sigma_{\epsilon,j}^2/\widehat\nu_j\le4/3$ on $\mathcal E$,
\eqref{eq:adaptive-training-quadratic-process}, and
$C_6=4C_\epsilon/3$.
Applying the Chernoff bound to
\eqref{eq:adaptive-noise-conditional-mgf} gives, for every $t>0$,
\begin{align}
\mathbb P_\epsilon\left(
|\psi(f)|>
\sqrt{\frac{2C_6\lambda_0t}{N_{\widehat\nu}}
\bigl(\mathcal R_{\lambda_0}(f)+\phi(f)\bigr)}
\right)
\le2e^{-t}.
\label{eq:adaptive-noise-conditional-tail}
\end{align}
By \eqref{eq:adaptive-quadratic-localized-tail}, with
$\mathbb P_{\mathrm{train}}$-probability at least $1-2e^{-t}$,
\begin{align}
\mathcal R_{\lambda_0}(f)+\phi(f)
\le\frac98\mathcal R_{\lambda_0}(f)
+C_5\frac{\lambda_0t}{N_{\widehat\nu}}.
\label{eq:adaptive-empirical-quadratic-upper-bound}
\end{align}
For a constant $C_7>0$ depending only on $M,c_\epsilon,C_\epsilon$,
we also have
\begin{align}
\sqrt{\frac{2C_6\lambda_0t}{N_{\widehat\nu}}
\left(\frac98\mathcal R_{\lambda_0}(f)
+C_5\frac{\lambda_0t}{N_{\widehat\nu}}\right)}
&\le\sqrt{\frac{9C_6\lambda_0t}{4N_{\widehat\nu}}
\mathcal R_{\lambda_0}(f)}
+\sqrt{2C_5C_6}\frac{\lambda_0t}{N_{\widehat\nu}}
\nonumber\\
&\le\frac1{16}\mathcal R_{\lambda_0}(f)
+C_7\frac{\lambda_0t}{N_{\widehat\nu}},
\label{eq:adaptive-noise-threshold-bound}
\end{align}
where the first inequality uses $\sqrt{u+v}\le\sqrt u+\sqrt v$ for
$u,v\ge0$; and the second inequality uses
$\sqrt{uv}\le u/16+4v$ and combines the coefficients of
$\lambda_0t/N_{\widehat\nu}$ into $C_7$.
Taking expectation over the training covariates in
\eqref{eq:adaptive-noise-conditional-tail} and using
\eqref{eq:adaptive-empirical-quadratic-upper-bound} and
\eqref{eq:adaptive-noise-threshold-bound}, a union bound gives
\begin{align}
\mathbb P_{\mathrm{train}}\left(
|\psi(f)|>\frac1{16}\mathcal R_{\lambda_0}(f)
+C_7\frac{\lambda_0t}{N_{\widehat\nu}}
\right)\le4e^{-t}.
\label{eq:adaptive-noise-localized-tail}
\end{align}
For any nonempty finite $\mathcal G\subseteq\mathcal F_{\lambda_0}$
chosen independently of the training observations, apply
\eqref{eq:adaptive-noise-localized-tail} with $t$ replaced by
$\log(4|\mathcal G|)+t$. A union bound shows that, with
$\mathbb P_{\mathrm{train}}$-probability at least $1-e^{-t}$,
\begin{align}
|\psi(g)|\le\frac1{16}\mathcal R_{\lambda_0}(g)
+C_7\frac{\lambda_0}{N_{\widehat\nu}}
\{\log(4|\mathcal G|)+t\}
\qquad\text{for every }g\in\mathcal G.
\nonumber
\end{align}

{\bf Step 4.}
Set $\eta=N^{-5}$. Since $\lambda_0\in\Lambda$ and
$1/\widehat\nu_j\le4/(3c_\epsilon)$ on $\mathcal E$, we have
\begin{align}
1\le\lambda_0\le N_{\widehat\nu}\le\frac{4N}{9c_\epsilon}.
\label{eq:adaptive-training-grid-size-bound}
\end{align}
By \Cref{lem:network-covering}, $\mathcal F_{\lambda_0}$ is totally
bounded and hence separable. Choose a finite $\eta$-net
$\mathcal G\subseteq\mathcal F_{\lambda_0}$ in the uniform norm,
independently of the training observations. The architecture choices
\eqref{eq:adaptive-network-depth}, \eqref{eq:adaptive-network-width}, and
\eqref{eq:adaptive-network-sparsity}, together with
\eqref{eq:adaptive-training-grid-size-bound} and
\eqref{eq:adaptive-oracle-approximation-scale}, give
\begin{align}
\log(6|\mathcal G|)+1
\le C_8\left(\frac{N_{\widehat\nu}}
{\lambda_0\{\log(eN)\}^3}\right)^{K/(2p+K)}
\{\log(eN)\}^3.
\label{eq:adaptive-training-net-entropy}
\end{align}
Define
\begin{align}
H(f)=|\phi(f)|+2|\psi(f)|-\frac14\mathcal R_{\lambda_0}(f),
\qquad f\in\mathcal F_{\lambda_0}.
\label{eq:adaptive-training-combined-process}
\end{align}
By \eqref{eq:adaptive-quadratic-localized-tail} and
\eqref{eq:adaptive-noise-localized-tail}, a union bound gives,
for every fixed $f\in\mathcal F_{\lambda_0}$ and $t>0$,
\begin{align}
\mathbb P_{\mathrm{train}}\left(
H(f)>C_9\frac{\lambda_0t}{N_{\widehat\nu}}\right)
\le6e^{-t},
\qquad C_9=C_5+2C_7.
\label{eq:adaptive-training-combined-tail}
\end{align}
To bound the Lipschitz constant of $H$, define
\begin{align}
T_\epsilon
=\frac1{N_{\widehat\nu}}
\sum_{j=1}^{\rmJ}\sum_{i\in\mathcal I_{j,1}}
\frac{w_0(X_{ij})}{\widehat\nu_j}|\epsilon_{ij}|.
\nonumber
\end{align}
Its expectation satisfies
\begin{align}
\mathbb E_{\mathrm{train}}T_\epsilon
&\le\frac{\lambda_0}{N_{\widehat\nu}}
\sum_{j=1}^{\rmJ}\frac{n_j\sigma_{\epsilon,j}}{3\widehat\nu_j}
 \le\frac{4N}{9\sqrt{c_\epsilon}},
\label{eq:adaptive-training-noise-modulus-mean}
\end{align}
where the first inequality uses $w_0\le\lambda_0$ and
$\mathbb E|\epsilon_{ij}|\le\sigma_{\epsilon,j}$;
and the second inequality uses
$\sigma_{\epsilon,j}/\widehat\nu_j
\le\sqrt{(\sigma_{\epsilon,j}^2/\widehat\nu_j)(1/\widehat\nu_j)}
\le4/(3\sqrt{c_\epsilon})$ on $\mathcal E$ and
$\lambda_0\le N_{\widehat\nu}$.
For every $f,g\in\mathcal F_{\lambda_0}$, we have
\begin{align}
|H(f)-H(g)|
&\le|\phi(f)-\phi(g)|+2|\psi(f)-\psi(g)|
+\frac14|\mathcal R_{\lambda_0}(f)-\mathcal R_{\lambda_0}(g)|
\nonumber\\
&\le(9M\lambda_0+2T_\epsilon)\|f-g\|_\infty,
\label{eq:adaptive-training-combined-continuity}
\end{align}
where the first inequality uses the triangle inequality in
\eqref{eq:adaptive-training-combined-process}; and the second inequality
uses $|(f-f^*)^2-(g-f^*)^2|\le4M\|f-g\|_\infty$,
$w_0\le\lambda_0$, and
$\sum_j n_j/(3N_{\widehat\nu}\widehat\nu_j)=1$ in
\eqref{eq:adaptive-clipped-training-risk},
\eqref{eq:adaptive-training-quadratic-process}, and
\eqref{eq:adaptive-training-noise-process}.
Writing $Z=\sup_{f\in\mathcal F_{\lambda_0}}(H(f))_+$,
we obtain
\begin{align}
\mathbb E_{\mathrm{train}}Z
&\le C_9\frac{\lambda_0}{N_{\widehat\nu}}
\{\log(6|\mathcal G|)+1\}
+\eta\{9M\lambda_0+2\mathbb E_{\mathrm{train}}T_\epsilon\}
\nonumber\\
&\le C_{10}\left(\frac{\lambda_0\{\log(eN)\}^3}{N_{\widehat\nu}}\right)^{2p/(2p+K)}+C_{11}N^{-4},
\label{eq:adaptive-training-uniform-expectation}
\end{align}
where the first inequality applies
\Cref{lem:random-lipschitz-supremum-bound} with
$a=C_9\lambda_0/N_{\widehat\nu}$, $b=6$, and
$L=9M\lambda_0+2T_\epsilon$, using
\eqref{eq:adaptive-training-combined-tail},
\eqref{eq:adaptive-training-combined-continuity}, and
\eqref{eq:adaptive-training-noise-modulus-mean}; and the second
inequality uses \eqref{eq:adaptive-training-net-entropy},
\eqref{eq:adaptive-training-noise-modulus-mean},
\eqref{eq:adaptive-training-grid-size-bound}, and $\eta=N^{-5}$.

{\bf Step 5.}
For every $f_0\in\mathcal F_{\lambda_0}$, we have
\begin{align}
\mathcal R_{\lambda_0}(\widehat f_{\lambda_0})
&\le\mathcal R_{\lambda_0}(f_0)
+|\phi(f_0)|+2|\psi(f_0)|
+|\phi(\widehat f_{\lambda_0})|+2|\psi(\widehat f_{\lambda_0})|
\nonumber\\
&\le\frac54\mathcal R_{\lambda_0}(f_0)
+\frac14\mathcal R_{\lambda_0}(\widehat f_{\lambda_0})+2Z,
\label{eq:adaptive-training-localized-basic-inequality}
\end{align}
where the first inequality applies the triangle inequality to
\eqref{eq:adaptive-training-basic-inequality}; and the second inequality
uses \eqref{eq:adaptive-training-combined-process} and
$H(f)\le Z$ for every $f\in\mathcal F_{\lambda_0}$.
Set
\begin{align}
\delta=\left(\frac{\lambda_0\{\log(eN)\}^3}{N_{\widehat\nu}}\right)^{p/(2p+K)}.
\nonumber
\end{align}
By \eqref{eq:adaptive-oracle-approximation-scale}, $0<\delta\le1$.
The choices \eqref{eq:adaptive-network-depth},
\eqref{eq:adaptive-network-width}, and
\eqref{eq:adaptive-network-sparsity}, together with
\eqref{eq:adaptive-training-grid-size-bound} and
\eqref{eq:adaptive-oracle-approximation-scale}, satisfy the architecture
requirements of \Cref{lem:hierarchical-network-approximation} for
sufficiently large $\sC_A,C_\sD,C_\sW,C_\sS$.
Thus, there is $f_0\in\mathcal F_{\lambda_0}$ such that
\begin{align}
\|f_0-f^*\|_\infty\le C_{12}\delta.
\label{eq:adaptive-training-network-approximation}
\end{align}
Moreover,
\begin{align}
\int w_0\,d\mathcal P_{\widehat\alpha}
&\le\int\widehat r\,d\mathcal P_{\widehat\alpha}
 \le2+\int|\widehat r-\max\{1,\rho_{\widehat\alpha}\}|
\,d\mathcal P_{\widehat\alpha}
\nonumber\\
&\le2+\|\widehat r-\max\{1,\rho_{\widehat\alpha}\}\|_{\mathcal P_{\widehat\alpha}},
\label{eq:adaptive-training-weight-mass}
\end{align}
where the first inequality uses $w_0\le\widehat r$;
the second inequality uses
$\widehat r\le1+\rho_{\widehat\alpha}
+|\widehat r-\max\{1,\rho_{\widehat\alpha}\}|$ and
$\int\rho_{\widehat\alpha}\,d\mathcal P_{\widehat\alpha}=1$;
and the third inequality is Cauchy--Schwarz.
Therefore,
\begin{align}
\mathcal R_{\lambda_0}(f_0)
&\le C_{12}^2\delta^2(2+\|\widehat r-\max\{1,\rho_{\widehat\alpha}\}\|_{\mathcal P_{\widehat\alpha}})
\nonumber\\
&\le C_{13}\delta^2+\|\widehat r-\max\{1,\rho_{\widehat\alpha}\}\|_{\mathcal P_{\widehat\alpha}}^2,
\label{eq:adaptive-training-approximation-risk}
\end{align}
where the first inequality uses
\eqref{eq:adaptive-clipped-training-risk},
\eqref{eq:adaptive-training-network-approximation}, and
\eqref{eq:adaptive-training-weight-mass}; and the second inequality
uses $uv\le u^2/4+v^2$ with
$u=C_{12}^2\delta^2$, $v=\|\widehat r-\max\{1,\rho_{\widehat\alpha}\}\|_{\mathcal P_{\widehat\alpha}}$, and $\delta^4\le\delta^2$.
Consequently,
\begin{align}
\mathbb E_{\mathrm{train}}\mathcal R_{\lambda_0}(\widehat f_{\lambda_0})
&\le\frac53\mathcal R_{\lambda_0}(f_0)
+\frac83\mathbb E_{\mathrm{train}}Z
\nonumber\\
&\le C_{14}\left(\frac{\lambda_0\{\log(eN)\}^3}{N_{\widehat\nu}}\right)^{2p/(2p+K)}
\nonumber\\
&\quad+C_{15}\|\widehat r-\max\{1,\rho_{\widehat\alpha}\}\|_{\mathcal P_{\widehat\alpha}}^2+C_{16}N^{-4},
\label{eq:adaptive-clipped-training-risk-bound}
\end{align}
where the first inequality rearranges
\eqref{eq:adaptive-training-localized-basic-inequality} and takes
expectation over the training observations; and the second inequality
uses \eqref{eq:adaptive-training-approximation-risk},
\eqref{eq:adaptive-training-uniform-expectation}, and
$\delta^2=(\lambda_0\{\log(eN)\}^3/N_{\widehat\nu})^{2p/(2p+K)}$.

{\bf Step 6.}
For $\mathcal P_{\widehat\alpha}$-almost every $x$, we have
\begin{align}
\widehat r(x)
&=w_0(x)+(\widehat r(x)-\lambda_0)_+
\nonumber\\
&\le w_0(x)+|\widehat r(x)-\max\{1,\rho_{\widehat\alpha}(x)\}|
+(\rho_{\widehat\alpha}(x)-\lambda_0)_+
\nonumber\\
&\le2w_0(x)+|\widehat r(x)-\max\{1,\rho_{\widehat\alpha}(x)\}|^2
+(\rho_{\widehat\alpha}(x)-\lambda_0)_+,
\label{eq:adaptive-training-weight-transfer}
\end{align}
where the equality uses $w_0=\min\{\widehat r,\lambda_0\}$;
the first inequality uses $(u+v)_+\le|u|+(v)_+$ and
$(\max\{1,\rho_{\widehat\alpha}\}-\lambda_0)_+
=(\rho_{\widehat\alpha}-\lambda_0)_+$ for $\lambda_0\ge1$;
and the second inequality uses $|u|\le1+u^2$ and $w_0\ge1$.
Consequently, for every $f\in\mathcal F_{\lambda_0}$,
\begin{align}
\mathcal R(f)
&\le2\mathcal R_{\lambda_0}(f)+4M^2\|\widehat r-\max\{1,\rho_{\widehat\alpha}\}\|_{\mathcal P_{\widehat\alpha}}^2
\nonumber\\
&\quad+4M^2\int(\rho_{\widehat\alpha}-\lambda_0)_+
\,d\mathcal P_{\widehat\alpha}
\nonumber\\
&\le2\mathcal R_{\lambda_0}(f)+4M^2\|\widehat r-\max\{1,\rho_{\widehat\alpha}\}\|_{\mathcal P_{\widehat\alpha}}^2+C_{17}N^{-4},
\label{eq:adaptive-training-risk-transfer}
\end{align}
where the first inequality integrates
\eqref{eq:adaptive-training-weight-transfer} against
$(f-f^*)^2\,d\mathcal P_{\widehat\alpha}$ and uses $|f-f^*|\le2M$;
and the second inequality applies \Cref{lem:adaptive-pooled-tail-grid}.
It follows that
\begin{align}
\mathbb E_{\mathrm{train}}\mathcal R(\widehat f_{\lambda_0})
&\le2\mathbb E_{\mathrm{train}}
\mathcal R_{\lambda_0}(\widehat f_{\lambda_0})
\nonumber\\
&\quad+4M^2\|\widehat r-\max\{1,\rho_{\widehat\alpha}\}\|_{\mathcal P_{\widehat\alpha}}^2+C_{17}N^{-4}
\nonumber\\
&\le C_{18}\left(\frac{\lambda_0\{\log(eN)\}^3}{N_{\widehat\nu}}\right)^{2p/(2p+K)}
\nonumber\\
&\quad+C_{19}\|\widehat r-\max\{1,\rho_{\widehat\alpha}\}\|_{\mathcal P_{\widehat\alpha}}^2+C_{20}N^{-4},
\nonumber
\end{align}
where the first inequality applies
\eqref{eq:adaptive-training-risk-transfer} to $\widehat f_{\lambda_0}$
and takes expectation over the training observations; and the second
inequality uses \eqref{eq:adaptive-clipped-training-risk-bound}.
 
\end{proof}

\begin{lemma}
\label{lem:adaptive-validation-oracle}
Suppose the assumptions and the condition on $\neff$ in
\Cref{lem:adaptive-pooled-tail-grid} hold. Let $\lambda_0$ be chosen
as in \Cref{lem:adaptive-pooled-tail-grid}, and let
$\{\widehat f_\lambda:\lambda\in\Lambda\}$ and $\widehat f$ be the
training candidates and final estimator in
\Cref{sec:multisource-adaptive-algorithm}. Use $\mathcal R$ as defined in \Cref{lem:adaptive-training-error}.
Write $\mathbb E_{\mathrm{val}}$ for expectation over the validation
observations $\{(X_{ij},Y_{ij}):i\in\mathcal I_{j,2},\ j=1,\ldots,\rmJ\}$,
conditional on the data in $\mathcal I_{j,0}$ and $\mathcal I_{j,1}$,
and on $\widehat\rho_j$, for $j=1,\ldots,\rmJ$.
For sufficiently large $N$, on $\mathcal E$, we have
\begin{align}
\mathbb E_{\mathrm{val}}\|\widehat f-f^*\|_{\mathcal Q}^2
&\le C_1\mathcal R(\widehat f_{\lambda_0})
 +C_2\|\widehat r-\max\{1,\rho_{\widehat\alpha}\}\|_{\mathcal P_{\widehat\alpha}}^2
+C_3\frac{\lambda_0}{N_{\widehat\nu}}\log(eN)
+C_4N^{-4}.
\nonumber
\end{align}
The constants $C_1,C_2,C_3,C_4>0$ depend only on
$M,c_\epsilon,C_\epsilon$.
\end{lemma}

\begin{proof}
{\bf Step 1.}
Work on $\mathcal E$. Define
\begin{align}
w_0(x)=\min\{\widehat r(x),\lambda_0\},
\qquad
v_0(x)=\widehat r(x)-w_0(x).
\label{eq:adaptive-validation-weight-decomposition}
\end{align}
Then $1\le w_0\le\lambda_0$ and $v_0\ge0$.
For $\mathcal P_{\widehat\alpha}$-almost every $x$, we have
\begin{align}
v_0(x)
&=(\widehat r(x)-\lambda_0)_+
 \le|\widehat r(x)-\max\{1,\rho_{\widehat\alpha}(x)\}|
+(\rho_{\widehat\alpha}(x)-\lambda_0)_+,
\label{eq:adaptive-validation-remainder-pointwise}
\end{align}
where the equality uses
\eqref{eq:adaptive-validation-weight-decomposition}; and the inequality
uses $(u+v)_+\le|u|+(v)_+$ and
$(\max\{1,\rho_{\widehat\alpha}\}-\lambda_0)_+
=(\rho_{\widehat\alpha}-\lambda_0)_+$ for $\lambda_0\ge1$.
Therefore,
\begin{align}
\int v_0^2\,d\mathcal P_{\widehat\alpha}
&\le2\|\widehat r-\max\{1,\rho_{\widehat\alpha}\}\|_{\mathcal P_{\widehat\alpha}}^2
\nonumber\\
&\quad+2\int(\rho_{\widehat\alpha}-\lambda_0)_+^2
\,d\mathcal P_{\widehat\alpha}
\nonumber\\
&\le2\|\widehat r-\max\{1,\rho_{\widehat\alpha}\}\|_{\mathcal P_{\widehat\alpha}}^2+C_5N^{-4},
\label{eq:adaptive-validation-remainder-second-moment}
\end{align}
where the first inequality squares
\eqref{eq:adaptive-validation-remainder-pointwise}, uses
$(u+v)^2\le2u^2+2v^2$, and integrates with respect to
$\mathcal P_{\widehat\alpha}$; and the second inequality applies
\Cref{lem:adaptive-pooled-tail-grid}.

{\bf Step 2.}
For each $j$, the observations with $i\in\mathcal I_{j,1}$ are used
for training, and those with $i\in\mathcal I_{j,2}$ are used for
validation. These two collections of observations are independent.
The data in $\mathcal I_{j,0}$ and $\widehat\rho_j$, for
$j=1,\ldots,\rmJ$, are also independent of the data in
$\mathcal I_{j,2}$, for $j=1,\ldots,\rmJ$. Therefore, the training candidates
$\{\widehat f_\lambda:\lambda\in\Lambda\}$ are independent of the
validation observations.
Use $\mathcal R_{\lambda_0}$ from
\eqref{eq:adaptive-clipped-training-risk}. For each candidate
$f\in\{\widehat f_\lambda:\lambda\in\Lambda\}$, define
\begin{align}
\phi_{\mathrm{val}}(f)
=\frac1{N_{\widehat\nu}}
\sum_{j=1}^{\rmJ}\sum_{i\in\mathcal I_{j,2}}
\frac{w_0(X_{ij})}{\widehat\nu_j}
(f(X_{ij})-f^*(X_{ij}))^2
-\mathcal R_{\lambda_0}(f)
\label{eq:adaptive-validation-quadratic-process}
\end{align}
and
\begin{align}
\psi_{\mathrm{val}}(f)
=\frac1{N_{\widehat\nu}}
\sum_{j=1}^{\rmJ}\sum_{i\in\mathcal I_{j,2}}
\frac{w_0(X_{ij})}{\widehat\nu_j}
\epsilon_{ij}(f(X_{ij})-f^*(X_{ij})).
\label{eq:adaptive-validation-noise-process}
\end{align}
Define  
\begin{align}
\xi_{\mathrm{val}}(f)
&=\frac1{N_{\widehat\nu}}
\sum_{j=1}^{\rmJ}\sum_{i\in\mathcal I_{j,2}}
\frac{v_0(X_{ij})}{\widehat\nu_j}
\bigl\{(f(X_{ij})-f^*(X_{ij}))^2
-2\epsilon_{ij}(f(X_{ij})-f^*(X_{ij}))\bigr\}
\nonumber\\
&\quad-\int v_0(x)(f(x)-f^*(x))^2
\,d\mathcal P_{\widehat\alpha}(x).
\label{eq:adaptive-validation-remainder-process}
\end{align}
Expanding the validation excess loss gives
\begin{align}
&\frac1{N_{\widehat\nu}}
\sum_{j=1}^{\rmJ}\sum_{i\in\mathcal I_{j,2}}
\frac{\widehat r(X_{ij})}{\widehat\nu_j}
\bigl\{(Y_{ij}-f(X_{ij}))^2-(Y_{ij}-f^*(X_{ij}))^2\bigr\}
\nonumber\\
&\quad=\frac1{N_{\widehat\nu}}
\sum_{j=1}^{\rmJ}\sum_{i\in\mathcal I_{j,2}}
\frac{w_0(X_{ij})+v_0(X_{ij})}{\widehat\nu_j}
\bigl\{(f(X_{ij})-f^*(X_{ij}))^2
-2\epsilon_{ij}(f(X_{ij})-f^*(X_{ij}))\bigr\}
\nonumber\\
&\quad=\mathcal R(f)+\phi_{\mathrm{val}}(f)
-2\psi_{\mathrm{val}}(f)+\xi_{\mathrm{val}}(f),
\label{eq:adaptive-validation-excess-loss-decomposition}
\end{align}
where the first equality uses \eqref{eq:multisource-model} and
\eqref{eq:adaptive-validation-weight-decomposition}; and the second
equality uses \eqref{eq:adaptive-validation-quadratic-process},
\eqref{eq:adaptive-validation-noise-process}, and
\eqref{eq:adaptive-validation-remainder-process}.
The minimizing property \eqref{eq:adaptive-validation-selector} and
\eqref{eq:adaptive-validation-excess-loss-decomposition} imply
\begin{align}
\mathcal R(\widehat f)
&\le\mathcal R(\widehat f_{\lambda_0})
+\phi_{\mathrm{val}}(\widehat f_{\lambda_0})
-\phi_{\mathrm{val}}(\widehat f)
+2\{\psi_{\mathrm{val}}(\widehat f)
-\psi_{\mathrm{val}}(\widehat f_{\lambda_0})\}
\nonumber\\
&\quad+\xi_{\mathrm{val}}(\widehat f_{\lambda_0})
-\xi_{\mathrm{val}}(\widehat f).
\label{eq:adaptive-validation-basic-inequality}
\end{align}

{\bf Step 3.}
Write $\mathbb P_{\mathrm{val}}$ for probability over the validation
observations. The arguments establishing
\eqref{eq:adaptive-quadratic-localized-tail} and
\eqref{eq:adaptive-noise-localized-tail} use only $|f|,|f^*|\le M$,
$w_0\le\lambda_0$, and the bounds
$1/\widehat\nu_j\le4/(3c_\epsilon)$ and
$\sigma_{\epsilon,j}^2/\widehat\nu_j\le4/3$ on $\mathcal E$,
in addition to independence and \Cref{ass:noise}.
Thus, they apply to every training candidate when
$\mathcal I_{j,1}$ is replaced by $\mathcal I_{j,2}$, since both sets
have $n_j/3$ elements. Combining the two tail bounds by a union bound
gives, for every $\lambda\in\Lambda$ and $t>0$,
\begin{align}
\mathbb P_{\mathrm{val}}\left(
|\phi_{\mathrm{val}}(\widehat f_\lambda)|
+2|\psi_{\mathrm{val}}(\widehat f_\lambda)|
-\frac14\mathcal R_{\lambda_0}(\widehat f_\lambda)
>C_6\frac{\lambda_0t}{N_{\widehat\nu}}\right)
\le6e^{-t},
\label{eq:adaptive-validation-bounded-tail}
\end{align}
for a constant $C_6>0$ depending only on $M,c_\epsilon,C_\epsilon$.
Define
\begin{align}
Z_{\mathrm{val}}
=\max_{\lambda\in\Lambda}
\left(
|\phi_{\mathrm{val}}(\widehat f_\lambda)|
+2|\psi_{\mathrm{val}}(\widehat f_\lambda)|
-\frac14\mathcal R_{\lambda_0}(\widehat f_\lambda)
\right)_+.
\label{eq:adaptive-validation-bounded-maximum}
\end{align}
Applying \eqref{eq:adaptive-validation-bounded-tail} with $t$ replaced
by $\log(6|\Lambda|)+t$ and taking a union bound over $\lambda\in\Lambda$
yields
\begin{align}
\mathbb P_{\mathrm{val}}\left(
Z_{\mathrm{val}}>C_6\frac{\lambda_0}{N_{\widehat\nu}}
\{\log(6|\Lambda|)+t\}\right)\le e^{-t}.
\label{eq:adaptive-validation-bounded-maximum-tail}
\end{align}
The definition of $\Lambda$ and
\eqref{eq:adaptive-training-grid-size-bound} give
\begin{align}
|\Lambda|=1+\lfloor\log_2N_{\widehat\nu}\rfloor
\le C_7\log(eN).
\label{eq:adaptive-validation-grid-cardinality}
\end{align}
Consequently,
\begin{align}
\mathbb E_{\mathrm{val}}Z_{\mathrm{val}}
&\le C_6\frac{\lambda_0}{N_{\widehat\nu}}
\{\log(6|\Lambda|)+1\}
 \le C_8\frac{\lambda_0}{N_{\widehat\nu}}\log(eN),
\label{eq:adaptive-validation-bounded-expectation}
\end{align}
where the first inequality integrates
\eqref{eq:adaptive-validation-bounded-maximum-tail}, bounding the
probability by one below
$C_6\lambda_0\log(6|\Lambda|)/N_{\widehat\nu}$; and the second
inequality uses \eqref{eq:adaptive-validation-grid-cardinality}.

{\bf Step 4.}
Fix a candidate $f\in\{\widehat f_\lambda:\lambda\in\Lambda\}$ and
write $h=f-f^*$. Define
\begin{align}
V_{ij}(f)=\frac{v_0(X_{ij})}{\widehat\nu_j}
\{h(X_{ij})^2-2\epsilon_{ij}h(X_{ij})\},
\qquad i\in\mathcal I_{j,2}.
\nonumber
\end{align}
By \eqref{eq:adaptive-validation-remainder-process},
$\mathbb E\epsilon_{ij}=0$, independence of the errors and covariates,
and $n_j/(3N_{\widehat\nu}\widehat\nu_j)=\widehat\alpha_j$, we have
\begin{align}
\xi_{\mathrm{val}}(f)
=\frac1{N_{\widehat\nu}}
\sum_{j=1}^{\rmJ}\sum_{i\in\mathcal I_{j,2}}
\{V_{ij}(f)-\mathbb E_{\mathrm{val}}V_{ij}(f)\}.
\label{eq:adaptive-validation-centered-remainder-sum}
\end{align}
For each fixed $x\in[0,1]^d$, expanding the square gives
\begin{align}
\mathbb E_{\mathrm{val}}\{h(x)^2-2\epsilon_{ij}h(x)\}^2
&=h(x)^4-4h(x)^3\mathbb E\epsilon_{ij}
+4h(x)^2\mathbb E\epsilon_{ij}^2
\nonumber\\
&=h(x)^4+4\sigma_{\epsilon,j}^2h(x)^2,
\label{eq:adaptive-validation-noise-square-expansion}
\end{align}
where the first equality expands the square and uses linearity of
expectation; and the second equality uses
$\mathbb E\epsilon_{ij}=0$ and
$\mathbb E\epsilon_{ij}^2=\sigma_{\epsilon,j}^2$.
It follows that
\begin{align}
\mathbb E_{\mathrm{val}}V_{ij}(f)^2
=\frac1{\widehat\nu_j^2}
\int v_0(x)^2\{h(x)^4+4\sigma_{\epsilon,j}^2h(x)^2\}
\,d\mathcal P_j(x),
\label{eq:adaptive-validation-summand-second-moment}
\end{align}
where the equality uses independence of the errors and covariates to
integrate \eqref{eq:adaptive-validation-noise-square-expansion}
against $v_0(x)^2\,d\mathcal P_j(x)/\widehat\nu_j^2$.
For a constant $C_9>0$ depending only on $M,c_\epsilon$,
\begin{align}
\mathbb E_{\mathrm{val}}\xi_{\mathrm{val}}(f)^2
&=\frac1{N_{\widehat\nu}^2}
\sum_{j=1}^{\rmJ}\sum_{i\in\mathcal I_{j,2}}
\operatorname{Var}_{\mathrm{val}}(V_{ij}(f))
 \le\frac1{N_{\widehat\nu}^2}
\sum_{j=1}^{\rmJ}\sum_{i\in\mathcal I_{j,2}}
\mathbb E_{\mathrm{val}}V_{ij}(f)^2
\nonumber\\
&=\frac1{N_{\widehat\nu}^2}
\sum_{j=1}^{\rmJ}\frac{n_j}{3\widehat\nu_j^2}
\int v_0(x)^2\{h(x)^4+4\sigma_{\epsilon,j}^2h(x)^2\}
\,d\mathcal P_j(x)
 \le\frac{C_9}{N_{\widehat\nu}}
\int v_0^2\,d\mathcal P_{\widehat\alpha},
\label{eq:adaptive-validation-remainder-variance}
\end{align}
where the first equality uses
\eqref{eq:adaptive-validation-centered-remainder-sum} and independence
of the validation observations; the first inequality bounds each
variance by its second moment; the second equality uses
\eqref{eq:adaptive-validation-summand-second-moment} and
$|\mathcal I_{j,2}|=n_j/3$; and the second inequality uses $|h|\le2M$,
$1/\widehat\nu_j\le4/(3c_\epsilon)$,
$\sigma_{\epsilon,j}^2/\widehat\nu_j\le4/3$, and
$n_j/(3N_{\widehat\nu}\widehat\nu_j)=\widehat\alpha_j$.
Consequently,
\begin{align}
\mathbb E_{\mathrm{val}}
\max_{\lambda\in\Lambda}|\xi_{\mathrm{val}}(\widehat f_\lambda)|
&\le\left\{\sum_{\lambda\in\Lambda}
\mathbb E_{\mathrm{val}}\xi_{\mathrm{val}}(\widehat f_\lambda)^2
\right\}^{1/2}
 \le\left\{\frac{C_9|\Lambda|}{N_{\widehat\nu}}
\int v_0^2\,d\mathcal P_{\widehat\alpha}\right\}^{1/2}
 \le\frac12\int v_0^2\,d\mathcal P_{\widehat\alpha}
+\frac{C_9|\Lambda|}{2N_{\widehat\nu}}
\nonumber\\
&\le C_{10}\|\widehat r-\max\{1,\rho_{\widehat\alpha}\}\|_{\mathcal P_{\widehat\alpha}}^2
+C_{11}\frac{\lambda_0}{N_{\widehat\nu}}\log(eN)
+C_{12}N^{-4},
\label{eq:adaptive-validation-remainder-expectation}
\end{align}
where the first inequality uses Cauchy--Schwarz and
$\max_\lambda|u_\lambda|^2\le\sum_\lambda|u_\lambda|^2$;
the second inequality uses
\eqref{eq:adaptive-validation-remainder-variance};
the third inequality uses $\sqrt{uv}\le(u+v)/2$;
and the fourth inequality uses
\eqref{eq:adaptive-validation-remainder-second-moment},
\eqref{eq:adaptive-validation-grid-cardinality}, and $\lambda_0\ge1$.

{\bf Step 5.}
Combining the validation fluctuations gives
\begin{align}
\mathcal R(\widehat f)
&\le\mathcal R(\widehat f_{\lambda_0})
+\frac14\mathcal R_{\lambda_0}(\widehat f_{\lambda_0})
+\frac14\mathcal R_{\lambda_0}(\widehat f)
+2Z_{\mathrm{val}}
+2\max_{\lambda\in\Lambda}|\xi_{\mathrm{val}}(\widehat f_\lambda)|
\nonumber\\
&\le\frac54\mathcal R(\widehat f_{\lambda_0})
+\frac14\mathcal R(\widehat f)
+2Z_{\mathrm{val}}
+2\max_{\lambda\in\Lambda}|\xi_{\mathrm{val}}(\widehat f_\lambda)|,
\label{eq:adaptive-validation-localized-basic-inequality}
\end{align}
where the first inequality applies the triangle inequality to
\eqref{eq:adaptive-validation-basic-inequality} and uses
\eqref{eq:adaptive-validation-bounded-maximum}; and the second
inequality uses $w_0\le\widehat r$, which gives
$\mathcal R_{\lambda_0}(f)\le\mathcal R(f)$.
Therefore,
\begin{align}
\mathbb E_{\mathrm{val}}\mathcal R(\widehat f)
&\le\frac53\mathcal R(\widehat f_{\lambda_0})
+\frac83\mathbb E_{\mathrm{val}}Z_{\mathrm{val}}
+\frac83\mathbb E_{\mathrm{val}}
\max_{\lambda\in\Lambda}|\xi_{\mathrm{val}}(\widehat f_\lambda)|
\nonumber\\
&\le C_{13}\mathcal R(\widehat f_{\lambda_0})
\nonumber\\
&\quad+C_{14}\|\widehat r-\max\{1,\rho_{\widehat\alpha}\}\|_{\mathcal P_{\widehat\alpha}}^2
+C_{15}\frac{\lambda_0}{N_{\widehat\nu}}\log(eN)
+C_{16}N^{-4},
\label{eq:adaptive-validation-weighted-risk-bound}
\end{align}
where the first inequality rearranges
\eqref{eq:adaptive-validation-localized-basic-inequality} and takes
expectation over the validation observations; and the second inequality
uses \eqref{eq:adaptive-validation-bounded-expectation} and
\eqref{eq:adaptive-validation-remainder-expectation}.
For every candidate $f$, the target risk satisfies
\begin{align}
\|f-f^*\|_{\mathcal Q}^2
&=\int\rho_{\widehat\alpha}(x)(f(x)-f^*(x))^2
\,d\mathcal P_{\widehat\alpha}(x)
\nonumber\\
&\le\mathcal R(f)
+\int|\widehat r(x)-\max\{1,\rho_{\widehat\alpha}(x)\}|
(f(x)-f^*(x))^2\,d\mathcal P_{\widehat\alpha}(x)
\nonumber\\
&\le2\mathcal R(f)+4M^2\|\widehat r-\max\{1,\rho_{\widehat\alpha}\}\|_{\mathcal P_{\widehat\alpha}}^2,
\label{eq:adaptive-validation-target-risk-transfer}
\end{align}
where the equality uses
$d\mathcal Q=\rho_{\widehat\alpha}\,d\mathcal P_{\widehat\alpha}$;
the first inequality uses
$\rho_{\widehat\alpha}\le\widehat r+|\widehat r-\max\{1,\rho_{\widehat\alpha}\}|$;
and the second inequality uses $|u|\le1+u^2$,
$\widehat r\ge1$, and $|f-f^*|\le2M$.
Applying \eqref{eq:adaptive-validation-target-risk-transfer} to
$\widehat f$, taking expectation over the validation observations,
and using \eqref{eq:adaptive-validation-weighted-risk-bound}
proves the lemma.
For later use, under the condition on $\neff$ in
\Cref{lem:adaptive-training-error}, we have
\begin{align}
\frac{\lambda_0}{N_{\widehat\nu}}\log(eN)
&\le\frac{\lambda_0\{\log(eN)\}^3}{N_{\widehat\nu}}
\nonumber\\
&\le\left(\frac{\lambda_0\{\log(eN)\}^3}
{N_{\widehat\nu}}\right)^{2p/(2p+K)},
\label{eq:adaptive-validation-cost-comparison}
\end{align}
where the first inequality uses $\log(eN)\ge1$; and the second
inequality uses \eqref{eq:adaptive-oracle-approximation-scale} and
$0<2p/(2p+K)<1$.
\end{proof}

\subsection{Proof of \Cref{thm:multisource-adaptive-rate}}

\begin{proof}[Proof of \Cref{thm:multisource-adaptive-rate}]
{\bf Step 1.}
Let $C_4\ge1$ be at least the constants in the conditions on
$\neff$ in
\Cref{lem:adaptive-pooled-tail-grid,lem:adaptive-training-error}.
If $\neff\le C_4\{\log(eN)\}^{7/2}$, then
\begin{align}
\mathbb E\|\widehat f-f^*\|_{\mathcal Q}^2
&\le4M^2
 \le C_5\neff^{-2p/(2p+K)}
\{\log(eN)\}^{7p/(2p+K)},
\nonumber
\end{align}
where the first inequality uses $\|f^*\|_\infty\le M$ from
\Cref{ass:hierarchical-composition} and $\|\widehat f\|_\infty\le M$; and the second inequality
uses $\neff\le C_4\{\log(eN)\}^{7/2}$,
$0<2p/(2p+K)<1$, and $C_4\ge1$, with $C_5=4M^2C_4$.
This proves the desired bound in this case. Henceforth, assume
\begin{align}
\neff>C_4\{\log(eN)\}^{7/2}.
\nonumber
\end{align}

{\bf Step 2.}
Let $\mathcal E$ be defined in \Cref{lem:adaptive-pooled-ratio-error},
and use $\mathcal R$ from \Cref{lem:adaptive-training-error}. On $\mathcal E$,
\Cref{lem:adaptive-pooled-tail-grid} gives $N_{\widehat\nu}\ge1$ and
$\lambda_0\in\Lambda$ satisfying
\begin{align}
0<\frac{\lambda_0}{N_{\widehat\nu}}
\le\min\left\{1,\frac{C_6\sqrt{\log(eN)}}{\neff}\right\}.
\label{eq:adaptive-theorem-threshold-ratio}
\end{align}
Taking expectation in \Cref{lem:adaptive-validation-oracle} over
the data in $\mathcal I_{j,1}$, for $j=1,\ldots,\rmJ$, and
applying \Cref{lem:adaptive-training-error} yields
\begin{align}
\mathbb E_{\mathrm{train}}\mathbb E_{\mathrm{val}}
\|\widehat f-f^*\|_{\mathcal Q}^2
&\le C_7\left\{
\left(\frac{\lambda_0\{\log(eN)\}^3}
{N_{\widehat\nu}}\right)^{2p/(2p+K)}
+\frac{\lambda_0}{N_{\widehat\nu}}\log(eN)\right\}
\nonumber\\
&\quad+C_8\|\widehat r-\max\{1,\rho_{\widehat\alpha}\}\|_{\mathcal P_{\widehat\alpha}}^2+C_9N^{-4}.
\label{eq:adaptive-theorem-training-validation-bound}
\end{align}
Here $\lambda_0$, $N_{\widehat\nu}$, $\widehat r$, and
$\mathcal P_{\widehat\alpha}$ depend only on the data in
$\mathcal I_{j,0}$ and on $\widehat\rho_j$, for
$j=1,\ldots,\rmJ$.
Consequently, on $\mathcal E$,
\begin{align}
\mathbb E_{\mathrm{train}}\mathbb E_{\mathrm{val}}
\|\widehat f-f^*\|_{\mathcal Q}^2
&\le C_{10}\left(\frac{\lambda_0\{\log(eN)\}^3}
{N_{\widehat\nu}}\right)^{2p/(2p+K)}
\nonumber\\
&\quad+C_{11}\|\widehat r-\max\{1,\rho_{\widehat\alpha}\}\|_{\mathcal P_{\widehat\alpha}}^2+C_{12}N^{-4}
\nonumber\\
&\le C_{13}\neff^{-2p/(2p+K)}
\{\log(eN)\}^{7p/(2p+K)}
\nonumber\\
&\quad+C_{11}\|\widehat r-\max\{1,\rho_{\widehat\alpha}\}\|_{\mathcal P_{\widehat\alpha}}^2+C_{12}N^{-4},
\label{eq:adaptive-theorem-good-event-conditional-risk}
\end{align}
where the first inequality combines
\eqref{eq:adaptive-theorem-training-validation-bound} and
\eqref{eq:adaptive-validation-cost-comparison}; and the second
inequality uses \eqref{eq:adaptive-theorem-threshold-ratio}.

{\bf Step 3.}
Taking expectation over all data gives
\begin{align}
\mathbb E\left[\mathbf1_{\mathcal E}
\|\widehat f-f^*\|_{\mathcal Q}^2\right]
&=\mathbb E\left[\mathbf1_{\mathcal E}
\mathbb E_{\mathrm{train}}\mathbb E_{\mathrm{val}}
\|\widehat f-f^*\|_{\mathcal Q}^2\right]
\nonumber\\
&\le C_{13}\neff^{-2p/(2p+K)}
\{\log(eN)\}^{7p/(2p+K)}
\nonumber\\
&\quad+C_{11}\mathbb E\left[\mathbf1_{\mathcal E}\|\widehat r-\max\{1,\rho_{\widehat\alpha}\}\|_{\mathcal P_{\widehat\alpha}}^2\right]+C_{12}N^{-4}
\nonumber\\
&\le C_{13}\neff^{-2p/(2p+K)}
\{\log(eN)\}^{7p/(2p+K)}
+C_{14}\sum_{j=1}^{\rmJ}\frac{\varepsilon_j^2}{\alpha_j^*}
+C_{12}N^{-4},
\label{eq:adaptive-theorem-good-event-risk}
\end{align}
where the equality uses the law of total expectation and the fact
that $\mathcal E$ depends only on observations with
$i\in\mathcal I_{j,0}$; the first inequality uses
\eqref{eq:adaptive-theorem-good-event-conditional-risk} and
$\mathbb P(\mathcal E)\le1$; and the second inequality applies
\Cref{lem:adaptive-pooled-ratio-error}.
On the complementary event, we have
\begin{align}
\mathbb E\left[\mathbf1_{\mathcal E^c}
\|\widehat f-f^*\|_{\mathcal Q}^2\right]
&\le4M^2\mathbb P(\mathcal E^c)
 \le4M^2N^{-3},
\label{eq:adaptive-theorem-bad-event-risk}
\end{align}
where the first inequality uses
$\|\widehat f\|_\infty,\|f^*\|_\infty\le M$; and the second inequality
uses \Cref{lem:adaptive-variance-weight-concentration} and the
definition of $\mathcal E$ in \Cref{lem:adaptive-pooled-ratio-error}.
Consequently,
\begin{align}
\mathbb E\|\widehat f-f^*\|_{\mathcal Q}^2
&=\mathbb E\left[\mathbf1_{\mathcal E}
\|\widehat f-f^*\|_{\mathcal Q}^2\right]
+\mathbb E\left[\mathbf1_{\mathcal E^c}
\|\widehat f-f^*\|_{\mathcal Q}^2\right]
\nonumber\\
&\le C_{15}\neff^{-2p/(2p+K)}
\{\log(eN)\}^{7p/(2p+K)}
+C_{16}\sum_{j=1}^{\rmJ}\frac{\varepsilon_j^2}{\alpha_j^*}
+C_{17}N^{-3},
\nonumber
\end{align}
where the equality uses $\mathbf1_{\mathcal E}+\mathbf1_{\mathcal E^c}=1$;
and the inequality combines \eqref{eq:adaptive-theorem-good-event-risk}
and \eqref{eq:adaptive-theorem-bad-event-risk}, using $N^{-4}\le N^{-3}$.
Together with the case $\neff\le C_4\{\log(eN)\}^{7/2}$, this proves
the theorem.
\end{proof}

\section{Auxiliary results}

\begin{lemma}
\label{lem:subgaussian-tail-moment}
Let $Z$ be a real-valued random variable. Suppose that, for constants
$a,B,s>0$,
\begin{align}
\mathbb P(Z>su)\le B e^{-au^2}
\qquad\text{for every }u\ge1.
\nonumber
\end{align}
Then, for every $t\ge1$, we have
\begin{align}
\mathbb E(Z-st)_+\le\frac{Bs}{2at}e^{-at^2}.
\nonumber
\end{align}
\end{lemma}

\begin{proof}
This is a standard result.
\end{proof}

\begin{lemma}
\label{lem:cutoff-contraction}
For every $u,v,a\in\mathbb R$, it holds that 
\begin{align}
|\max\{u,a\}-\max\{v,a\}|\le|u-v|,
\nonumber
\end{align}
and that \begin{align}
|\min\{u,a\}-\min\{v,a\}|\le|u-v|.
\nonumber
\end{align}
\end{lemma}
\begin{proof}
The results follow from triangle inequalities.
\end{proof}
\begin{lemma}
\label{lem:hierarchical-network-approximation}
Suppose \Cref{ass:hierarchical-composition} holds, and fix $\sC_M\ge M$.
There exist constants $\sC_A\ge1$ and $C_1,C_2,C_3,C_4>0$,
depending only on $d,L,p,K,R,M$, such that the following holds.
For every $0<\delta\le1$ and integers $\sD,\sW,\sS$ satisfying
\begin{align}
\sD&\ge\max\{2,C_1\log(2/\delta)\}, \quad 
\sW\ge\max\{d,C_2\delta^{-K/p}\}, \quad 
\sS\ge C_3\sD\delta^{-K/p},
\end{align}
there is a network $f_\delta\in\mathcal F(\sD,\sW,\sS,\sC_A,\sC_M)$ such that
\begin{align}
\|f_\delta-f^*\|_\infty\le C_4\delta,
\nonumber
\end{align}
and that  $\|f^*\|_\infty\le M$.
\end{lemma}

\begin{proof}
By Theorem~5 and Lemma~3 of \citet{schmidthieber2020}, choosing
the network budgets $D,W,S$ as in
\Cref{lem:hierarchical-network-approximation} gives an
approximating network $g_\delta$ satisfying
\begin{align}
\|g_\delta-f^*\|_\infty\le C_5\delta,
\label{eq:hierarchical-approximation-unclipped}
\end{align}
where $C_5>0$ is independent of $\delta,\sD,\sW,\sS$ and $f^*$.
Define $f_\delta=\max\{-M,\min\{M,g_\delta\}\}$.
Clipping gives
\begin{align}
\|f_\delta-f^*\|_\infty
&\le\|g_\delta-f^*\|_\infty
 \le C_5\delta,
\nonumber
\end{align}
where the first inequality uses $\|f^*\|_\infty\le M$ and
clipping to $[-M,M]$; and the second inequality uses
\eqref{eq:hierarchical-approximation-unclipped}.
For sufficiently large $\sC_A,C_1,C_2,C_3$, we have
$f_\delta\in\mathcal F(\sD,\sW,\sS,\sC_A,\sC_M)$, proving the lemma.
\end{proof}

\begin{lemma}
\label{lem:network-covering}
Let $\sD\ge2$, $\sW\ge d$, and $\sS\ge\max\{\sW,\sD\}$ be integers, and
let $\sC_A\ge1$ and $\sC_M>0$. There exists a universal constant
$C_1>0$ such that, for every $0<\delta\le1$,
\begin{align}
&\log\mathcal N\bigl(
\delta,\mathcal F(\sD,\sW,\sS,\sC_A,\sC_M),\|\cdot\|_\infty\bigr)
\nonumber\\
&\qquad\le C_1\sS\left\{
\sD\log\bigl((\sW+1)\sC_A\bigr)+\log(4/\delta)
\right\}.
\nonumber
\end{align}
Here $\|\cdot\|_\infty$ is the uniform norm on $[0,1]^d$, and
$\mathcal N$ denotes the covering number with centers in the covered
class.
\end{lemma}

\begin{proof}
Write $\mathcal F=\mathcal F(\sD,\sW,\sS,\sC_A,\sC_M)$ and
$\mathcal G=\mathcal R(d,\sW,\sD,\sC_A,\sS)$, using the notation of
\citet{oubolcskei2026} and restricting the functions to $[0,1]^d$.
Then $\mathcal F\subseteq\mathcal G$, and continuity identifies
the uniform norm with the $L^\infty([0,1]^d)$ norm. Therefore,
\begin{align}
\log\mathcal N(\delta,\mathcal F,\|\cdot\|_\infty)
&\le\log\mathcal N(\delta/4,\mathcal G,\|\cdot\|_\infty)
\nonumber\\
&\le C_1\min\{\sS,\sW^2\sD\}
\log\left(\frac{4(\sW+1)^\sD\sC_A^\sD}{\delta}\right)
\nonumber\\
&\le C_1\sS\left\{
\sD\log\bigl((\sW+1)\sC_A\bigr)+\log(4/\delta)
\right\},
\nonumber
\end{align}
where the first inequality uses Lemma~F.1 of \citet{oubolcskei2026}
and monotonicity in the covering radius;
the second inequality applies their Theorem~5.1 with norm exponent
$\infty$ and radius $\delta/4<1/2$;
and the third inequality uses $\min\{\sS,\sW^2\sD\}\le \sS$ and expands the logarithm.
\end{proof}

\begin{lemma}
\label{lem:conditional-noise-moments}
Under \Cref{ass:noise}, for every integer $m\ge2$, there exists a
constant $C>0$, depending only on $m$ and $C_\epsilon$, such that
\begin{align}
\mathbb E|\epsilon_{ij}|^m\le C\sigma_{\epsilon,j}^m,
\qquad i=1,\ldots,n_j,\quad j=1,\ldots,\rmJ.
\nonumber
\end{align}
\end{lemma}

\begin{proof}
If $\sigma_{\epsilon,j}=0$, then \Cref{ass:noise} implies
$\mathbb E|\epsilon_{ij}|^m=0$.
For $\sigma_{\epsilon,j}>0$, we have
\begin{align}
\mathbb E|\epsilon_{ij}|^m
&\le m!\sigma_{\epsilon,j}^m
\mathbb E\exp\left(\frac{|\epsilon_{ij}|}{\sigma_{\epsilon,j}}\right)
\nonumber\\
&\le2e^{C_\epsilon/2}m!\sigma_{\epsilon,j}^m,
\nonumber
\end{align}
where the first inequality uses $u^m\le m!e^u$ for $u\ge0$;
and the second inequality uses $e^{|u|}\le e^u+e^{-u}$ and
\Cref{ass:noise} at $t=\pm1/\sigma_{\epsilon,j}$.
Taking $C=2e^{C_\epsilon/2}m!$ proves the assertion.
\end{proof}

\begin{lemma}
\label{lem:bernstein-moment-bound}
Let $\{Z_{ij}:1\le i\le n_j,\ 1\le j\le\rmJ\}$ be
independent, centered real-valued random variables. Suppose that,
for some $b>0$ and $v_{ij}\ge0$ and all
$i=1,\ldots,n_j$, $j=1,\ldots,\rmJ$,
\begin{align}
\operatorname{Var}(Z_{ij})\le v_{ij}^2
\label{eq:bernstein-variance-bound}
\end{align}
and
\begin{align}
\|Z_{ij}\|_{L^\infty}\le b.
\label{eq:bernstein-uniform-bound}
\end{align}
Then, for every $t>0$,
\begin{align}
\mathbb P\left(
\left|\sum_{j=1}^{\rmJ}\sum_{i=1}^{n_j}Z_{ij}\right|
>\sqrt{2t\sum_{j=1}^{\rmJ}\sum_{i=1}^{n_j}v_{ij}^2}+\frac{bt}{3}
\right)\le2e^{-t}.
\nonumber
\end{align}
\end{lemma}

\begin{proof}
The result follows from the standard exponential-moment proof of
Bernstein's inequality applied to \eqref{eq:bernstein-variance-bound}
and \eqref{eq:bernstein-uniform-bound},
with total variance bound $\sum_{j=1}^{\rmJ}\sum_{i=1}^{n_j}v_{ij}^2$
and common uniform bound $b$.
\end{proof}

\begin{lemma}
\label{lem:finite-maximum-exponential-tail}
Let $m\ge1$ be an integer and let $X_1,\ldots,X_m$ be real-valued
random variables satisfying
\begin{align}
\max_{1\le\ell\le m}X_\ell\ge0.
\nonumber
\end{align}
Suppose that, for some $a>0$, every $1\le\ell\le m$ and every
$t>0$ satisfy
\begin{align}
\mathbb P(X_\ell>at)\le2e^{-t}.
\label{eq:finite-maximum-individual-tail}
\end{align}
Then
\begin{align}
\mathbb E\left[\max_{1\le\ell\le m}X_\ell\right]
\le a\{\log(2m)+1\}.
\nonumber
\end{align}
No independence is required.
\end{lemma}
\begin{proof}
For every $t>0$, we have
\begin{align}
\mathbb P\left(\max_{1\le\ell\le m}X_\ell
>a\{\log(2m)+t\}\right)
&\le\sum_{\ell=1}^m
\mathbb P\left(X_\ell>a\{\log(2m)+t\}\right)
\nonumber\\
&\le2m\exp\{-\log(2m)-t\}
\nonumber\\
&=e^{-t},
\label{eq:finite-maximum-shifted-tail}
\end{align}
where the first inequality uses the union bound;
the second inequality uses \eqref{eq:finite-maximum-individual-tail};
and the equality cancels the factor $2m$.
Consequently,
\begin{align}
\mathbb E\left[\max_{1\le\ell\le m}X_\ell\right]
&\le a\log(2m)+a\int_0^\infty e^{-t}\,dt
\nonumber\\
&=a\{\log(2m)+1\},
\nonumber
\end{align}
where the inequality uses the tail-integral formula for the
nonnegative maximum, bounds probabilities by $1$ below $a\log(2m)$, and applies
\eqref{eq:finite-maximum-shifted-tail} after the substitution
$x=a\{\log(2m)+t\}$; and the equality evaluates the integral.
\end{proof}

\begin{lemma}
\label{lem:localized-supremum-bound}
Let $\nu$ be a finite nonnegative measure on $[0,1]^d$ satisfying
$\nu([0,1]^d)\le1$, and let $\mathcal F$ be a 
class of   functions from $[0,1]^d$ to $\mathbb R$
 such that 
$\sup_{f\in\mathcal F}\|f\|_\infty\le \sC_M$ for some $\sC_M>0$.
Let $\mathcal X(f)$, $f\in\mathcal F$, be a real-valued random process.
Suppose that, for deterministic constants $a,L,\tau>0$, every
$f,g\in\mathcal F$ and $t>0$ satisfy
\begin{align}
\mathbb P\left(
|\mathcal X(f)-\mathcal X(g)|-\tau\|f-g\|_\nu^2>at
\right)\le2e^{-t},
\label{eq:localized-supremum-tail}
\end{align}
and, simultaneously for all $f,g\in\mathcal F$,
\begin{align}
|\mathcal X(f)-\mathcal X(g)|\le L\|f-g\|_\infty.
\label{eq:localized-supremum-continuity}
\end{align}
Then, for every $\delta>0$,
\begin{align}
&\mathbb E\left[
\sup_{f,g\in\mathcal F}
\left\{|\mathcal X(f)-\mathcal X(g)|-\tau\|f-g\|_\nu^2\right\}
\right]
 \le
a\bigl\{2\log \mathcal N(\delta,\mathcal F,\|\cdot\|_\infty)
+\log2+1\bigr\}
+(2L+8\tau \sC_M)\delta,
\nonumber
\end{align}
where the covering centers belong to $\mathcal F$.
\end{lemma}

\begin{proof}
{\bf Step 1.}
Fix $\delta>0$ and choose a deterministic $\delta$-net
$\mathcal F_0\subset\mathcal F$ with
$|\mathcal F_0|=\mathcal N(\delta,\mathcal F,\|\cdot\|_\infty)$.
For $f,g\in\mathcal F_0$, define
\begin{align}
X_{f,g}=|\mathcal X(f)-\mathcal X(g)|
-\tau\|f-g\|_\nu^2.
\label{eq:localized-supremum-net-variables}
\end{align}
Then
\begin{align}
&\mathbb E\left[
\max_{f,g\in\mathcal F_0}
\left\{|\mathcal X(f)-\mathcal X(g)|-\tau\|f-g\|_\nu^2\right\}
\right]
 \le a\bigl\{\log(2|\mathcal F_0|^2)+1\bigr\}
=a\bigl\{2\log \mathcal N(\delta,\mathcal F,\|\cdot\|_\infty)
+\log2+1\bigr\},
\label{eq:localized-supremum-net-maximum}
\end{align}
where the inequality applies
\Cref{lem:finite-maximum-exponential-tail} to the
$|\mathcal F_0|^2$ random variables in
\eqref{eq:localized-supremum-net-variables}, using
\eqref{eq:localized-supremum-tail} and $X_{f,f}=0$ to ensure
that their maximum is nonnegative; and the equality uses the
chosen cardinality of $\mathcal F_0$.

{\bf Step 2.}
For each $f,g\in\mathcal F$, choose $f_0,g_0\in\mathcal F_0$
with $\|f-f_0\|_\infty,\|g-g_0\|_\infty\le\delta$.
Simultaneously for all $f,g\in\mathcal F$,
\begin{align}
|\mathcal X(f)-\mathcal X(g)|
&\le |\mathcal X(f_0)-\mathcal X(g_0)|+|\mathcal X(f)-\mathcal X(f_0)|+|\mathcal X(g)-\mathcal X(g_0)|
\nonumber\\
&\le |\mathcal X(f_0)-\mathcal X(g_0)|+2L\delta,
\label{eq:localized-supremum-process-comparison}
\end{align}
where the first inequality uses the triangle inequality; and the
second inequality uses \eqref{eq:localized-supremum-continuity}
and the net approximation.
For every bounded measurable $h:[0,1]^d\to\mathbb R$,
\begin{align}
\|h\|_\nu
&\le\|h\|_\infty\sqrt{\nu([0,1]^d)}
\nonumber\\
&\le\|h\|_\infty,
\label{eq:localized-supremum-norm-comparison}
\end{align}
where the first inequality integrates $|h(x)|^2\le\|h\|_\infty^2$;
and the second inequality uses $\nu([0,1]^d)\le1$.
Therefore,
\begin{align}
\left|\|f-g\|_\nu^2-\|f_0-g_0\|_\nu^2\right|
&\le\bigl(\|f-g\|_\nu+\|f_0-g_0\|_\nu\bigr)
\bigl(\|f-f_0\|_\nu+\|g-g_0\|_\nu\bigr)
\nonumber\\
&\le8\sC_M\delta,
\label{eq:localized-supremum-penalty-comparison}
\end{align}
where   the second inequality
uses the uniform bound $\sC_M$, the net approximation, and
\eqref{eq:localized-supremum-norm-comparison}.
Consequently,
\begin{align}
&\left\{|\mathcal X(f)-\mathcal X(g)|-\tau\|f-g\|_\nu^2\right\}
\nonumber\\
&\quad\le
\left\{|\mathcal X(f_0)-\mathcal X(g_0)|-\tau\|f_0-g_0\|_\nu^2\right\}
+(2L+8\tau \sC_M)\delta
\nonumber\\
&\quad\le
\max_{u,v\in\mathcal F_0}
\left\{|\mathcal X(u)-\mathcal X(v)|-\tau\|u-v\|_\nu^2\right\}
+(2L+8\tau \sC_M)\delta,
\label{eq:localized-supremum-net-comparison}
\end{align}
where the first inequality uses
\eqref{eq:localized-supremum-process-comparison} and
\eqref{eq:localized-supremum-penalty-comparison}; and the second inequality bounds
the selected pair by the maximum over $\mathcal F_0$.
 
Taking the supremum over $f,g\in\mathcal F$ and expectations in
\eqref{eq:localized-supremum-net-comparison}, then applying
\eqref{eq:localized-supremum-net-maximum}, proves the assertion.
\end{proof}

\begin{lemma}
\label{lem:random-lipschitz-supremum-bound}
Let $\mathcal F$ be a nonempty separable class of bounded real-valued
functions on $[0,1]^d$, equipped with the uniform norm, and let
$H(f)$, $f\in\mathcal F$, be a real-valued random process.
Suppose there is a nonnegative random variable $L$ with
$\mathbb E L<\infty$ such that, almost surely, simultaneously for all
$f,g\in\mathcal F$,
\begin{align}
|H(f)-H(g)|\le L\|f-g\|_\infty.
\label{eq:random-lipschitz-supremum-continuity}
\end{align}
Suppose also that, for deterministic constants $a>0$ and $b\ge1$,
every $f\in\mathcal F$ and $t>0$ satisfy
\begin{align}
\mathbb P\{H(f)>at\}\le b e^{-t}.
\label{eq:random-lipschitz-supremum-tail}
\end{align}
Then, for every $\eta>0$ and every nonempty finite $\eta$-net
$\mathcal G\subseteq\mathcal F$ in the uniform norm,
\begin{align}
\mathbb E\left[\sup_{f\in\mathcal F}(H(f))_+\right]
\le a\{\log(b|\mathcal G|)+1\}+\eta\,\mathbb E L.
\nonumber
\end{align}
Here $\mathcal G$ is deterministic and $(u)_+=\max\{u,0\}$.
\end{lemma}

\begin{proof}
{\bf Step 1.}
Write $m=|\mathcal G|$ and $V=\max_{g\in\mathcal G}(H(g))_+$.
For every $t>0$,
\begin{align}
\mathbb P(V>at)
&\le\sum_{g\in\mathcal G}\mathbb P\{H(g)>at\}
\nonumber\\
&\le bm e^{-t},
\label{eq:random-lipschitz-net-tail}
\end{align}
where the first inequality is the union bound; and the second
inequality uses \eqref{eq:random-lipschitz-supremum-tail}.
Therefore,
\begin{align}
\mathbb E V
&\le a\log(bm)+\int_{a\log(bm)}^\infty bm e^{-u/a}\,du
\nonumber\\
&=a\{\log(bm)+1\},
\label{eq:random-lipschitz-net-expectation}
\end{align}
where the inequality integrates the tail probability, using the bound
one on $[0,a\log(bm)]$ and \eqref{eq:random-lipschitz-net-tail} on
$[a\log(bm),\infty)$; and the equality evaluates the integral.

{\bf Step 2.}
For every $f\in\mathcal F$, choose $g\in\mathcal G$ with
$\|f-g\|_\infty\le\eta$. Simultaneously for all $f\in\mathcal F$,
\begin{align}
(H(f))_+
&\le(H(g))_++L\eta
 \le V+L\eta,
\nonumber
\end{align}
where the first inequality uses
\eqref{eq:random-lipschitz-supremum-continuity}, the net approximation,
and $L\ge0$; and the second inequality uses the definition of $V$.
Separability of $\mathcal F$ and
\eqref{eq:random-lipschitz-supremum-continuity} allow the supremum to
be taken over a countable dense subset, so it is measurable.
Taking the supremum over $f\in\mathcal F$ and then expectations gives
\begin{align}
\mathbb E\left[\sup_{f\in\mathcal F}(H(f))_+\right]
&\le\mathbb E V+\eta\,\mathbb E L
 \le a\{\log(b|\mathcal G|)+1\}+\eta\,\mathbb E L,
\nonumber
\end{align}
where the first inequality uses the common bound $V+L\eta$; and the
second inequality uses \eqref{eq:random-lipschitz-net-expectation}.
\end{proof}

\section{Proof related to \Cref{thm:holder-covariate-shift-lower-bound}}

\begin{proof}[Proof of \Cref{thm:holder-covariate-shift-lower-bound}]

{\bf Step 1.}
Let $\mathcal Q$ be the uniform distribution on $[0,1]^d$,
and let $\delta_0$ be the point mass at the origin. For each
$j=1,\ldots,\rmJ$, set
\begin{align}
\mathcal P_j
=\frac{1}{1\vee\sigma_{\rho,j}}\mathcal Q
+\left(1-\frac{1}{1\vee\sigma_{\rho,j}}\right)\delta_0.
\nonumber
\end{align}
Since $\mathcal Q(\{0\})=0$, we have $\mathcal Q\ll\mathcal P_j$
and may take
\begin{align}
\rho_j(x)=(1\vee\sigma_{\rho,j})\mathbf 1_{\{x\ne0\}}.
\nonumber
\end{align}
If $0\le\sigma_{\rho,j}\le1$, then $\mathcal P_j=\mathcal Q$
and $\mathbb E_{X\sim\mathcal P_j}\exp\{t(\rho_j(X)-1)\}=1$
for every $t\in\mathbb R$.
If $\sigma_{\rho,j}>1$, then $0\le\rho_j\le\sigma_{\rho,j}$ and
$\mathbb E_{X\sim\mathcal P_j}[\rho_j(X)]=1$. Thus, for every
$t\in\mathbb R$,
\begin{align}
\mathbb E_{X\sim\mathcal P_j}
\left[\exp\{t(\rho_j(X)-1)\}\right]
\le\exp\left(\frac{\sigma_{\rho,j}^2t^2}{8}\right)
\le\exp\left(\frac{\sigma_{\rho,j}^2t^2}{2}\right),
\nonumber
\end{align}
where the first inequality uses Hoeffding's lemma; and the second
inequality uses $\sigma_{\rho,j}^2t^2\ge0$.
Hence \Cref{ass:holder-shift} holds for every source.
Take mutually independent $X_{ij}\sim\mathcal P_j$ and mutually
independent $\epsilon_{ij}\sim N(0,\sigma_{\epsilon,j}^2)$,
independent of all covariates. Then \Cref{ass:noise} holds since
$C_\epsilon\ge1$.

{\bf Step 2.}
Set $a=1\wedge R\wedge M$. By
\Cref{lem:holder-packing-lower-bound}, there exist constants
$c_1,c_2,\varepsilon_0>0$, depending only on $p,K,a$, such that,
for every $0<\varepsilon\le\varepsilon_0$, we can choose
$(p,a)$-smooth functions $h_1,\ldots,h_J:[0,1]^K\to[-a,a]$,
supported in $[1/4,3/4]^K$, satisfying
\begin{align}
\log J\ge c_1\varepsilon^{-K/p}.
\label{eq:holder-packing-cardinality}
\end{align}
Their pairwise distances satisfy
\begin{align}
\|h_j-h_k\|_{L^2([0,1]^K)}\ge\varepsilon
\qquad\text{for all }j\ne k.
\nonumber
\end{align}
Moreover,
\begin{align}
\|h_j\|_{L^2([0,1]^K)}\le c_2\varepsilon
\qquad\text{for all }1\le j\le J.
\nonumber
\end{align}
We will choose $\varepsilon$ in Step 4.

{\bf Step 3.}
Set $\eta=R/(1+R)$ and define $u:[-R,R]^K\to[-R,R]$ by
\begin{align}
u(z)=\eta z_1.
\nonumber
\end{align}
Then $0<\eta<1$, and $u$ is $(p,R)$-smooth because the expression
in \Cref{def:holder} equals $\eta(R+1)=R$.
For $j=1,\ldots,J$ and $\ell=1,\ldots,L$, define
\begin{align}
g_{j,\ell}(x)=\eta^{\ell-1}h_j(x_1,\ldots,x_K).
\nonumber
\end{align}
Since $a\le R$, we have $g_{j,1}\in\mathcal H(1,p,K,R)$.
For $\ell\ge2$, the definitions give
\begin{align}
g_{j,\ell}(x)=u\bigl(g_{j,\ell-1}(x),0,\ldots,0\bigr).
\nonumber
\end{align}
The zero function belongs to every level of the hierarchical class,
so \Cref{def:hierarchical} gives
$g_{j,\ell}\in\mathcal H(\ell,p,K,R)$ by induction.
Set
\begin{align}
f_j(x)=g_{j,L}(x)=\eta^{L-1}h_j(x_1,\ldots,x_K).
\nonumber
\end{align}
Thus $f_j\in\mathcal H(L,p,K,R)$ and $\|f_j\|_\infty\le a\le M$,
so \Cref{ass:hierarchical-composition} holds. The support of $h_j$ gives $f_j(0)=0$.
Since the first $K$ coordinates under $\mathcal Q$ are uniform on
$[0,1]^K$, the packing bounds give
\begin{align}
\|f_j-f_k\|_{\mathcal Q}\ge\eta^{L-1}\varepsilon
\qquad\text{for all }1\le j\ne k\le J.
\label{eq:holder-lower-separation}
\end{align}
Moreover,
\begin{align}
\|f_j\|_{\mathcal Q}\le\eta^{L-1}c_2\varepsilon
\qquad\text{for all }1\le j\le J.
\label{eq:holder-lower-function-norm}
\end{align}

{\bf Step 4.}
Recall from \eqref{eq:multisource-effective-sample-size} that
\begin{align}
\neff=\sum_{j=1}^{\rmJ}
\frac{n_j}{(1+\sigma_{\rho,j})\sigma_{\epsilon,j}^2}.
\nonumber
\end{align}
For $f\in\{0,f_1,\ldots,f_J\}$, let $\mathbb P_f$ denote the joint
law of all source observations under regression function $f$.
For any $f,g\in\{0,f_1,\ldots,f_J\}$,
\begin{align}
\operatorname{KL}(\mathbb P_f,\mathbb P_g)
&=\frac12\sum_{j=1}^{\rmJ}
\frac{n_j}{\sigma_{\epsilon,j}^2}\|f-g\|_{\mathcal P_j}^2
\nonumber\\
&=\frac12\sum_{j=1}^{\rmJ}
\frac{n_j}{(1\vee\sigma_{\rho,j})\sigma_{\epsilon,j}^2}
\|f-g\|_{\mathcal Q}^2
\nonumber\\
&\le\neff\|f-g\|_{\mathcal Q}^2,
\label{eq:holder-lower-kl}
\end{align}
where the first equality uses \Cref{lem:gaussian-regression-kl};
the second equality uses the definition of $\mathcal P_j$ and
$f(0)=g(0)=0$; and the inequality uses
$1\vee\sigma_{\rho,j}\ge(1+\sigma_{\rho,j})/2$ and
\eqref{eq:multisource-effective-sample-size}.
Choose $c_0>0$ satisfying
\begin{align}
c_0\le\min\left\{
\varepsilon_0,
\left(\frac{c_1}{\log4}\right)^{p/K},
\left(\frac{c_1}{8\eta^{2(L-1)}c_2^2}\right)^{p/(2p+K)}
\right\},
\label{eq:holder-lower-scale-constant}
\end{align}
and set
\begin{align}
\varepsilon=c_0\neff^{-p/(2p+K)}.
\label{eq:holder-lower-packing-scale}
\end{align}
By \eqref{eq:oracle-effective-sample-size} and
\eqref{eq:multisource-effective-sample-size}, $\neff\ge1$,
so $0<\varepsilon\le c_0\le\varepsilon_0$. Moreover,
\begin{align}
\log J
\ge c_1\varepsilon^{-K/p}
\ge c_1c_0^{-K/p}
\ge\log4,
\nonumber
\end{align}
where the first inequality uses \eqref{eq:holder-packing-cardinality};
the second inequality uses $\varepsilon\le c_0$; and the third
inequality uses \eqref{eq:holder-lower-scale-constant}.
In particular, $J\ge4$. Also,
\begin{align}
\neff\varepsilon^{2+K/p}
=c_0^{2+K/p}
\le\frac{c_1}{8\eta^{2(L-1)}c_2^2},
\label{eq:holder-lower-scale-check}
\end{align}
where the equality uses \eqref{eq:holder-lower-packing-scale};
and the inequality uses \eqref{eq:holder-lower-scale-constant}.
The average Kullback--Leibler divergence satisfies
\begin{align}
\frac1J\sum_{k=1}^J\operatorname{KL}(\mathbb P_{f_k},\mathbb P_0)
&\le\frac{\neff}{J}\sum_{k=1}^J\|f_k\|_{\mathcal Q}^2
\nonumber\\
&\le\neff\eta^{2(L-1)}c_2^2\varepsilon^2
\nonumber\\
&\le\frac{\eta^{2(L-1)}c_2^2}{c_1}
\neff\varepsilon^{2+K/p}\log J
\nonumber\\
&\le\frac18\log J,
\label{eq:holder-fano-kl-condition}
\end{align}
where the first inequality uses \eqref{eq:holder-lower-kl} with $g=0$;
the second inequality uses \eqref{eq:holder-lower-function-norm};
the third inequality uses \eqref{eq:holder-packing-cardinality};
and the fourth inequality uses \eqref{eq:holder-lower-scale-check}.

{\bf Step 5.}
The source distributions and likelihood ratios are identical for all
$f_1,\ldots,f_J$. The estimation form of Fano's inequality gives
\begin{align}
&\inf_{\widetilde f}
\sup_{\{\mathcal P_j\}_{j=1}^{\rmJ},\mathcal Q,f^*,\epsilon}
\mathbb E\|\widetilde f-f^*\|_{\mathcal Q}^2
\nonumber\\
&\quad\ge\frac{\eta^{2(L-1)}\varepsilon^2}{4}
\left[1-\frac{
J^{-1}\sum_{k=1}^J\operatorname{KL}(\mathbb P_{f_k},\mathbb P_0)
+\log2}{\log J}\right]
\nonumber\\
&\quad\ge\frac{3\eta^{2(L-1)}\varepsilon^2}{32},
\label{eq:holder-fano-risk-bound}
\end{align}
where the first inequality uses the admissibility of the constructed
models, Fano's inequality, \eqref{eq:holder-lower-separation}, and
reference distribution $\mathbb P_0$; and the second inequality uses
\eqref{eq:holder-fano-kl-condition} and $J\ge4$.
Substituting \eqref{eq:holder-lower-packing-scale} into
\eqref{eq:holder-fano-risk-bound} and using
\eqref{eq:multisource-effective-sample-size} proves the theorem with
$C_1=3\eta^{2(L-1)}c_0^2/32>0$, depending only on $L,p,K,d,R,M$.

\end{proof}

 \begin{lemma}
\label{lem:holder-packing-lower-bound}
Let $p\ge1$, let $K\ge1$ be an integer, and let $a>0$.
Let $\mathcal Q$ be the uniform distribution on $[0,1]^K$.
Using \Cref{def:holder}, define
\begin{align}
\mathcal G
=\Bigl\{h:[0,1]^K\to\mathbb R:\;&
h\text{ is $(p,a)$-smooth},
\nonumber\\
&\operatorname{supp}(h)\subseteq[1/4,3/4]^K\Bigr\}.
\nonumber
\end{align}
There exist constants $c_1,c_2,\varepsilon_0>0$, depending only on
$p,K,a$, such that, for every $0<\varepsilon\le\varepsilon_0$,
there are $h_1,\ldots,h_J\in\mathcal G$ satisfying
\begin{align}
\log J\ge c_1\varepsilon^{-K/p}.
\nonumber
\end{align}
These functions are separated in $L^2(\mathcal Q)$:
\begin{align}
\|h_j-h_k\|_{L^2(\mathcal Q)}\ge\varepsilon
\qquad\text{for all }j\ne k.
\nonumber
\end{align}
Moreover,
\begin{align}
\|h_j\|_{L^2(\mathcal Q)}\le c_2\varepsilon
\qquad\text{for all }1\le j\le J.
\nonumber
\end{align}
\end{lemma}

\begin{proof}
The result follows from the H\"older packing construction in the proof
of Theorem~3 of \citet{schmidthieber2020}, after rescaling the support
to $[1/4,3/4]^K$ and the smoothness bound to $a$.
\end{proof}

\begin{lemma}
\label{lem:gaussian-regression-kl}
Let $\mathcal P_1,\ldots,\mathcal P_{\rmJ}$ be probability
measures on $[0,1]^d$, let $n_j\ge1$ be integers, and let
$\sigma_{\epsilon,j}>0$ for $j=1,\ldots,\rmJ$.
For $f\in\bigcap_{j=1}^{\rmJ}L^2(\mathcal P_j)$, let
$\mathbb P_f$ denote the joint law of the observations
\begin{align}
X_{ij}\sim\mathcal P_j,
\qquad Y_{ij}=f(X_{ij})+\epsilon_{ij},
\qquad i=1,\ldots,n_j,\quad j=1,\ldots,\rmJ,
\nonumber
\end{align}
where the covariates are mutually independent, and the errors
$\epsilon_{ij}\sim N(0,\sigma_{\epsilon,j}^2)$ are mutually independent
and independent of all covariates.
Then, for any $f,g\in\bigcap_{j=1}^{\rmJ}L^2(\mathcal P_j)$,
\begin{align}
\operatorname{KL}(\mathbb P_f,\mathbb P_g)
=\frac12\sum_{j=1}^{\rmJ}
\frac{n_j}{\sigma_{\epsilon,j}^2}\|f-g\|_{\mathcal P_j}^2.
\nonumber
\end{align}
\end{lemma}

\begin{proof}
 The conditional Gaussian densities and independence give
\begin{align}
\log\frac{d\mathbb P_f}{d\mathbb P_g}
=\frac12\sum_{j=1}^{\rmJ}\frac{1}{\sigma_{\epsilon,j}^2}
\sum_{i=1}^{n_j}
\left\{(Y_{ij}-g(X_{ij}))^2-(Y_{ij}-f(X_{ij}))^2\right\}.
\label{eq:gaussian-regression-log-likelihood}
\end{align}

Write $\mathbb E_f$ for expectation under $\mathbb P_f$. Then
\begin{align}
\operatorname{KL}(\mathbb P_f,\mathbb P_g)
&=\mathbb E_f\left[\log\frac{d\mathbb P_f}{d\mathbb P_g}\right]
\nonumber\\
&=\frac12\sum_{j=1}^{\rmJ}\frac{1}{\sigma_{\epsilon,j}^2}
\sum_{i=1}^{n_j}\mathbb E_f\left[
(f(X_{ij})-g(X_{ij}))^2
+2\epsilon_{ij}(f(X_{ij})-g(X_{ij}))\right]
\nonumber\\
&=\frac12\sum_{j=1}^{\rmJ}\frac{1}{\sigma_{\epsilon,j}^2}
\sum_{i=1}^{n_j}\mathbb E_f\left[(f(X_{ij})-g(X_{ij}))^2\right]
\nonumber\\
&=\frac12\sum_{j=1}^{\rmJ}
\frac{n_j}{\sigma_{\epsilon,j}^2}\|f-g\|_{\mathcal P_j}^2,
\nonumber
\end{align}
where the first equality is the definition of Kullback--Leibler
divergence; the second equality uses
\eqref{eq:gaussian-regression-log-likelihood} and
$Y_{ij}=f(X_{ij})+\epsilon_{ij}$; the third equality uses the centered
errors independent of the covariates; and the fourth equality uses
$X_{ij}\sim\mathcal P_j$.
\end{proof}

\end{document}